\documentclass[11pt,reqno]{amsart}
\usepackage[T1]{fontenc}
\usepackage{lmodern}
\usepackage{amsmath,amsfonts,amssymb,amsthm,mathtools,booktabs}
\usepackage{microtype}
\usepackage{mathrsfs}
\usepackage[colorlinks=true,urlcolor=blue,citecolor=red,linkcolor=blue]{hyperref}
\newtheorem{theorem}{Theorem}[section]
\newtheorem{lemma}[theorem]{Lemma}
\newtheorem{proposition}[theorem]{Proposition}
\newtheorem{corollary}[theorem]{Corollary}
\newtheorem{conjecture}[theorem]{Conjecture}
\theoremstyle{definition}

\newtheorem{remark}[theorem]{Remark}
\numberwithin{equation}{section}
\newcommand{\R}{\mathbb R}
\newcommand{\Gal}{\operatorname{Gal}}
\newcommand{\eps}{\varepsilon}

\newcommand{\Q}{\mathbb Q}
\newcommand{\Z}{\mathbb Z}
\newcommand{\N}{\mathbb N}
\newcommand{\F}{\mathbb F}
\newcommand{\C}{\mathbb C}
\newcommand{\OO}{\mathcal O}
\newcommand{\HTP}{\operatorname{HTP}}
\newcommand{\Nm}{\operatorname{N}}
\newcommand{\Tr}{\operatorname{Tr}}

\newcommand{\im}{\operatorname{Im}}

\def\bg{\bigg}
\def\({\bg(}
\def\){\bg)}
\def\t{\text}

\def\mo{{\rm{mod}\ }}
\def\pmod#1{\ (\mo\ #1)}

\def\se {\subseteq}

\newcommand{\Sha}{\operatorname{Sha}}
\begin{document}
\hbox{Preprint, {\tt arXiv:2611.25002}}
\medskip

\title[Hilbert's tenth problem over quadratic fields]{On Hilbert's tenth problem \\ over quadratic fields}
\author[Zhi-Wei Sun]{Zhi-Wei Sun}
\address{School of Mathematics, Nanjing University,
Nanjing 210093, People's Republic of China}
\email{zwsun@nju.edu.cn}

\begin{abstract}
In 2026, OpenAI proved that Hilbert's Tenth Problem over the field $\mathbb Q$
of rational numbers has a negative solution. In this paper we study Hilbert's Tenth Problem
over the quadratic fields $\mathbb Q(\sqrt2),\mathbb Q(\sqrt3),\mathbb Q(\sqrt6)$. 
Let $K=\mathbb Q(\sqrt d)$ with $d\in\{2,3,6\}$. We prove that there is no algorithm to decide
for any polynomial equation $P(x_1,\ldots,x_{38})=0$ with integer coefficients whether it has solutions with $x_1,\ldots,x_{38}\in K$. We also show the existence of a polynomial 
$f(x_0,x_1,\ldots,x_{19})$ with integer coefficients such that $a\in K$ is a rational integer if and only if $f(a,x_1,\ldots,x_{19})=0$ for some $x_1,\ldots,x_{19}\in K$.
\end{abstract}
\maketitle

\section{Introduction}\label{sec:intro}

Let $R$ be a subring containing $\Z$ of the field $\C$ of complex numbers.
Hilbert's Tenth Problem (HTP in short) over $R$ asks for an algorithm to decide for
any given polynomial equation with integer coefficients whether it has solutions over $R$.
A subset $A$ of $R$ is called {\it diophantine} over $R$ if there is a polynomial 
$P(x_0,x_1,\ldots,x_n)$ with integer coefficients such that for any $a\in R$ we have
$$a\in A\iff \exists x_1,\ldots,x_n\in R\,[P(a,x_1,\ldots,x_n)=0].$$

The negative
solution of HTP over $\mathbb Z$ was obtained  by Matiyasevich \cite{Mat} in 1970 on the basis of Davis,
Putnam and Robinson \cite{DPR}. The author \cite{Sun2021} proved further that $\HTP_{11}(\Z)$
is undecidable, i.e.,
there is no algorithm to decide for any polynomial equations with integer coefficients and $11$ unknowns
whether it has integral solutions. Recently the author \cite{Sun10} improved this by reducing $11$ to $9$. 

For the field  $\Q$ of rational numbers,
HTP over $\Q$ was recently proved to be undecidable by 
OpenAI \cite{OAIHTP} via using advanced tools including 
 elliptic curves, $p$-adic formal groups, Diophantine geometry, Selmer groups and height theory. 
 Motivated by this work, the author \cite{Sun38} proved further that $\HTP_{38}(\Q)$ is undecidable.
It is now also known that $\Z$ is diophantine over $\Q$ (cf. \cite{Bur} and \cite{ZQ}).
 Sun and Li \cite{LiSun19} gave a $19$-unknown definition of $\Z$ over $\Q$.

Let $K$ be any number field (which is a finite extension of the field $\Q$), and let
$O_K$ be the ring of algebraic integers in $K$. By \cite{ABHS,KP}, $\Z$ is diophantine over $\OO_K$
and hence HTP over $\OO_K$ is undecidable by Matiyasevich's theorem.

After OpenAI's solution of HTP over $\Q$, it is natural to consider HTP over a number field $K$.

In this paper, we establish the following two theorems, using the arithmetic inputs specified in the proofs.

\begin{theorem}\label{thm:unconditional}
Let $K=\Q(\sqrt d)$ with $d\in\{2,3,6\}$. Then $\HTP_{38}(K)$ is undecidable, i.e., there is no algorithm deciding, for a given polynomial $f\in\Z[X_1,\ldots,X_{38}]$, whether $f$ has a zero in $K^{38}$.
\end{theorem}

\begin{theorem}\label{thm:integers25}
Let $K=\Q(\sqrt d)$ with $d\in\{2,3,6\}$. Then $\Z$ is diophantine over $K$. Moreover, there is a polynomial $\mathscr F_K\in\Z[X,Y_1,\ldots,Y_{19}]$ such that for every $x\in K$ we have
\[
 x\in\Z\iff\exists y_1,\ldots,y_{19}\in K\,[\mathscr F_K(x,y_1,\ldots,y_{19})=0].
\]
\end{theorem}

Motivated by these two theorems, we propose the following conjecture.

\begin{conjecture} Let $K$ be any number field. Then $\Z$ is diophantine over $K$, and hence HTP over $K$
is undecidable.
\end{conjecture}

The cyclic carriers and the restricted seven-contact height theorems are established first, including the corresponding inputs over $\Q(\sqrt3)$ and $\Q(\sqrt6)$ in Sections~\ref{sec:sqrt3} and~\ref{sec:sqrt6}. The two reductions then share the local comparison and the aggregate trace estimate, a five-vertex Selmer certificate with uniform scalar and contact expansions, and a single three-source square-denominator test. The split integer source and the simultaneous Gaussian prime-pattern input come from \cite{Sun38}; the shifted Pell source is adapted from \cite{LiSun19}. The rational covering step uses the pointwise $2$-converse/finiteness implication \cite{OAIPointwise} and the terminal step uses the square-denominator lemma of \cite{Bur}. In particular the rational-point assertion is not derived from local solubility alone. The proofs of both main theorems are completed in Section~\ref{sec:new-finish}.

For the initial geometric sections, we set
\[
 K=\Q(\sqrt2),\ \ F=\Q(\sqrt3)\ \ \t{and}\ \ L=KF=\Q(\sqrt2,\sqrt3).
\]

 For a number field $M$,
$\OO_M$ denotes its ring of integers and $G_M=\Gal(\overline\Q/M)$.
If $\mathcal P$ is a nonzero prime ideal of $\OO_M$, then
$v_{\mathcal P}$ is the integer-valued additive valuation normalized by
$v_{\mathcal P}(\mathcal P)=1$, and we put $v_{\mathcal P}(0)=+\infty$;
$M_{\mathcal P}$ and $N\mathcal P$ denote the completion and absolute
norm.  For a rational prime $p$, $v_p$ is normalized by $v_p(p)=1$.
Thus, for $a\in\Q^\times$ and $\mathcal P\mid p$,
$v_{\mathcal P}(a)=e(\mathcal P/p)v_p(a)$. 

Let $M$ be a number field. We use the symbol $h$ to denote the absolute logarithmic Weil height on
$\mathbb P^1(\overline\Q)$; explicitly, for $\alpha\in M$ we have
\[
 h(\alpha)=\frac1{[M:\Q]}
  \sum_v [M_v:\Q_v]\log\max\{1,|\alpha|_v\},
\]
where the absolute values restrict to the usual archimedean and
$p$-adic absolute values on $\Q$, and the sum is over all places
of $M$.  For an elliptic curve $E=E(M)$ over $M$, write $O_E$ for its
identity point.  Then $[m]Q$ denotes the $m$-fold sum of $Q$ under the
elliptic-curve group law: for $m>0$ it is the sum of $m$
copies of $Q$, $[0]Q=O_E$, and $[-m]Q=-[m]Q$.
The standard alternative notation $mQ$ has exactly the same
meaning and does not denote coordinatewise multiplication.
We put $mE(M)=\{[m]Q:Q\in E(M)\}$ and
$\Z P=\{[n]P:n\in\Z\}$,
$E(M)_{\rm tors}$ is the torsion subgroup, and $E[n]$ is the geometric
$n$-torsion.  

We also write $E_{\rm tor}=E(\overline\Q)_{\rm tors}$
when the geometric torsion group is meant.  If $Q\ne O_E$
is given by an affine
Weierstrass equation, then $x(Q)$ and $y(Q)$ are its
abscissa and ordinate; we regard $x(O_E)$ as the pole of
the abscissa function.  We write $\widehat h_E$ for the N\'eron--Tate canonical
height on $E(\overline\Q)$ associated with the origin divisor.  For a
finite prime $\mathcal P$ of a number field, $I_{\mathcal P}$ denotes an
inertia subgroup, $\mathrm{Fr}^{\rm ar}_{\mathcal P}$ an arithmetic
Frobenius element, and $\mathrm{Fr}^{\rm geo}_{\mathcal P}$ its inverse.
If $\rho$ is a two-dimensional linear representation, then $\mathbb P\rho$
denotes its projectivization to $\operatorname{PGL}_2$.  

To proves Theorems \ref{thm:unconditional} and \ref{thm:integers25}, we need the following basic theorem.

\begin{theorem}\label{thm:carrier} 
For the elliptic curve
\begin{equation}\label{eq:E}
E(F):\ y^2=x^3-a_E^2x\ \ \t{with}\ \ a_E=4+\sqrt3
\end{equation}
over $F=\Q(\sqrt3)$,
 we have
\[
\operatorname{rank}E(F)=\operatorname{rank}E(L)=1,\qquad
 E(L)_{\rm tors}=E[2],\qquad 2E(L)\subset E(F),
\]
where $L=\Q(\sqrt2,\sqrt3)$.
Moreover, there is a nontorsion point $P\in E(F)$ with
\begin{equation}\label{eq:carrier}2E(L)=\Z P.\end{equation}
\end{theorem}

Set $q_*=367$ and define
\begin{equation}\label{eq:ab}
D_*=q_*^2+3,\qquad
 a_* =\frac{(q_*-3)(q_*+1)}{D_*},\qquad
 b_* =\frac{(q_*-1)(q_*+3)}{D_*^2}.
\end{equation}
Put
\begin{equation}\label{eq:pullback}
x(z)=\frac{a_*z+b_*}{q_*z^2}\ \ \t{and}\ \ 
j(x)=\frac{16x^2}{(x-4)^3}.
\end{equation}
Our set of seven rational contacts is 
\begin{equation}\label{eq:contacts}
\begin{split}
\mathcal C=\bigg\{&\infty,0,
-\frac{(q_*-1)(q_*+3)}{(q_*-3)(q_*+1)D_*},
\frac{q_*-1}{D_*},-\frac{q_*+3}{q_*D_*},
\frac{(q_*-1)(q_*+3)}{4q_*D_*},-\frac1{D_*}\bigg\}.
\end{split}
\end{equation}
Let $S$ be a fixed finite set of finite primes of $K$ containing
all primes over rational primes at most $324$, all primes above $367$,
and all prime divisors of contact denominators and pairwise differences.
For $z\in K\setminus\mathcal C$, define $M_{\mathcal C,S}(z)$ as the
product of absolute norms $N\mathcal P$ at the primes
$\mathcal P\notin S$ having positive odd contact with a point of
$\mathcal C$. For finite contact $c$, the condition is
$v_{\mathcal P}(z-c)>0$ odd; at infinity it is
$v_{\mathcal P}(1/z)>0$ odd. As usual, an empty product is regarded as $1$.
Fix a prime $\mathfrak q_*\mid367$ and define the local tube by
\begin{equation}\label{eq:tube}
z\equiv1\pmod{\mathfrak q_*}.
\end{equation}
The polynomial used in the CM sign condition is
\begin{equation}\label{eq:PCintro}
 P_{\mathcal C}(z)=\frac{(a_*z+b_*-q_*z^2)
(a_*z+b_*-4q_*z^2)}{4q_*^2}.
\end{equation}
The main height assertion is the following restricted estimate:
\begin{equation}\label{eq:height-interface}
\begin{gathered}
 h(z)\le H M_{\mathcal C,S}(z)^c
 \\ \t{for }z\in K\setminus\mathcal C\ \t{with }z\equiv1\pmod{\mathfrak q_*}
 \\ \t{ and }\sigma(P_{\mathcal C}(z))>0
 \t{ for both real embeddings }\sigma,
\end{gathered}
\end{equation}
where $c$ and $H$ are positive constants independent of $z$.

\begin{theorem}\label{thm:restricted-height}
There are positive constants $c$ and $H$ such that
\eqref{eq:height-interface} holds.
\end{theorem}

\begin{proposition}[Seven-contact geometry]\label{prop:seven}
The pullback \eqref{eq:pullback} has the seven orbifold contacts
\eqref{eq:contacts} with all of residual order two. At $z=0$ this results
from a point of order four pulled back with local degree two.
Away from $S$, positive odd contact detects all possible ramification
of the normalized full-three quaternionic cover. On the tube
\eqref{eq:tube}, we have $v_{\mathfrak q_*}(j(x(z)))=1$.
\end{proposition}
The proof of Proposition~\ref{prop:seven} is given in
Section~\ref{sec:seven}. We use the standard literature \cite{BabaGranath,Deligne,Buzzard,LinYang} on Shimura curves and their moduli interpretations.

For the normalized descent, let $\mu_m$ denote the abstract group of
$m$th roots of unity with trivial $G_K$-action; for a one-cochain
$\beta$, write $\delta\beta(g,h)=\beta(g)\beta(h)\beta(gh)^{-1}$.

\begin{theorem}\label{thm:cochain}
Let $T_0$ be a finite set of nonzero prime ideals of $K$ containing its
unique prime ideal above $2$. Suppose that $\varepsilon\in Z^2(G_K,\{\pm1\})$
is continuous, factors through $\Gal(L_0/K)$ for some finite
Galois extension $L_0/K$ unramified outside $T_0$, and its Brauer
class is trivial outside $T_0$ at nonzero prime ideals. Define
\[
T_- =\{\mathcal P\in T_0:\mathcal P\nmid2,\ 
\operatorname{inv}_{\mathcal P}[\varepsilon]=1/2\}
\ \ \t{and}\ \ 
B_{T_0}=\prod_{\mathcal P\in T_-}2^{v_2(N\mathcal P-1)}.
\]
Then, there exist a finite-order character $\chi:G_K\to\mu_m$ and a
continuous one-cochain $\beta:G_K\to\mu_{2m}$ with
\[
\beta^2=\chi,\quad\delta\beta=\varepsilon,\quad m\le4B_{T_0},
\quad N_{K/\Q}(\mathfrak f_\chi)
\le32B_{T_0}^2\prod_{\mathcal P\in T_-}N\mathcal P.
\]
After a quadratic twist, the restriction $\beta|_{G_{L_0}}$ is
unramified outside $T_0$. If $L'$ is its kernel field, then
$[L':K]\le8B_{T_0}[L_0:K]$.
\end{theorem}
The proof of this descent is in Section~\ref{sec:cochain}.

\section{Two-isogeny descent over \texorpdfstring{$\Q(\sqrt3)$}{Q(sqrt(3))}}\label{sec:descent}

Put
\[
 \eps=2+\sqrt3\ \ \t{and}\ \  a=4+\sqrt3.
\]
The ring of integers of $F$ is $\OO_F=\Z[\sqrt3]$ and has discriminant $12$.  
It is well known that $\OO_F$ is a principal ideal domain.  Its unit group is $\{\pm\eps^n:n\in\Z\}$; in particular $\eps$ generates the totally positive units modulo squares.

For $e\in F^\times$, write
\[
 E_e:y^2=x^3-e^2x\ \ \t{and}\ \ 
 E'_e:y^2=x^3+4e^2x.
\]
The standard $2$-isogeny $E_e\to E'_e$ and its dual give descent maps
\[
 \alpha_e:E_e(F)/\widehat\phi E'_e(F)\longrightarrow F^\times/F^{\times2}
 \ \ \t{and}\ \ 
 \alpha'_e:E'_e(F)/\phi E_e(F)\longrightarrow F^\times/F^{\times2}.
\]
For a nonexceptional point, the map is represented by the square class of its $x$-coordinate; see \cite[Chapter X, Section 4]{Silv}.  The usual formula is
\begin{equation}\label{eq:isogeny-rank}
 2^{\operatorname{rank}E_e(F)}
 =\frac{|\im\alpha_e|\cdot|\im\alpha'_e|}{4}.
\end{equation}

\begin{lemma}\label{lem:primitive-quartic-local}
Let $M$ be a nonarchimedean local field of characteristic zero,
with valuation ring $\mathcal O_M$ and maximal ideal $\mathfrak m$.
For $a,b\in\mathcal O_M$, if
\[
 W^2=aU^4+bV^4,\qquad U,V,W\in M,
 \quad (U,V)\ne(0,0),
\]
then there is a solution with $U,V,W\in\mathcal O_M$
and at least one of $U,V$ a unit.  Consequently, failure
of a primitive solution modulo $\mathfrak m^n$ for some
$n\ge1$ rules out a nontrivial local solution.
This applies to either homogeneous two-isogeny covering
\[
 W^2=dU^4-\frac{c^2}{d}V^4,
 \qquad W^2=dU^4+\frac{4c^2}{d}V^4
\]
whenever the two coefficients on the right lie in
$\mathcal O_M$.  In particular the zero triple cannot
represent a descent class.
\end{lemma}
\begin{proof}
Let $m=\min\{v(U),v(V)\}$, and let $\varpi$ be a uniformizer.
The homogeneous change
$(U,V,W)\mapsto(\varpi^{-m}U,\varpi^{-m}V,\varpi^{-2m}W)$
preserves the equation and makes $U,V$ integral, at least
one a unit.  Since $a,b$ are integral, its right-hand side
is integral.  Hence $W^2$, and therefore $W$, is integral.
Reduction modulo $\mathfrak m^n$ proves the assertion.
\end{proof}
\subsection{An explicit non-torsion point}

\begin{lemma}\label{lem:P0}
The point
\begin{equation}\label{eq:P0}
 P_0=(244+139\sqrt3,\ 5334+3082\sqrt3)
\end{equation}
belongs to $E_a(F)$ and has infinite order.
\end{lemma}

\begin{proof}
Direct expansion gives
\[
 (5334+3082\sqrt3)^2
 =(244+139\sqrt3)^3-a^2(244+139\sqrt3).
\]
The prime $73$ splits in $F$, with the two square roots $21$ and $52$ of $3$ modulo $73$.  The two reductions of $E_a$ have orders
\begin{equation}\label{eq:73orders}
 \#E_a(\F_{73})=80\ \ \t{and}\ \ 
 \#E_a^\sigma(\F_{73})=68.
\end{equation}
At either prime above $73$ the formal group has no $73$-torsion, since the absolute ramification index is $1<72$; hence the full torsion subgroup injects into the corresponding reduction group.  Thus
$|E_a(F)_{\rm tors}|$ divides $\gcd(80,68)=4$.
The curve already contains all four points of $E[2]$, so the torsion subgroup is exactly $E[2]$.  The $y$-coordinate of $P_0$ is nonzero, hence $P_0$ is not torsion.
\end{proof}

\subsection{Finite support of the descent images}

Write
\[
 \pi=1+\sqrt3\ \ \t{and}\ \ \pi^2=2\eps.
\]
The prime above $2$ is generated by $\pi$, and $[2]=[\eps]$ in
$F^\times/F^{\times2}$.  The only prime dividing $a$ is
$\mathfrak a=(a)$, of norm $13$.

\begin{lemma}[Finite support]\label{lem:support}
For $e=a$ and $e=2a$, we have
\[
 \im\alpha_e\subseteq\langle-1,\eps,a\rangle\ \ \t{and}\ \ 
 \im\alpha'_e\subseteq\langle\eps,a\rangle
 \quad\text{in }F^\times/F^{\times2}.
\]
\end{lemma}

\begin{proof}
Away from the primes dividing $2a$, comparison of valuations in
$y^2=x(x^2\pm c)$ shows that the valuation of every nonzero
$x$-coordinate is even.  At the prime $(\pi)$, the coefficient of $x$ has
valuation $0$ or $4$ on $E_e$ and valuation $4$ or $8$ on $E'_e$.
If $v_\pi(x)$ were odd, the two terms in $x(x^2\pm c)$ would have
distinct valuations and the smaller one would be odd, which is impossible for
a square.  Hence the prime above $2$ does not contribute a new square class.
The only remaining finite support is $\mathfrak a$.  Since $\OO_F$ has class
number one and unit group $\{\pm\eps^n:n\in\Z\}$, the first inclusion
follows.  On $E'_e$ every nonzero real $x$-coordinate is positive at both
real embeddings, which excludes the class of $-1$ and proves the second
inclusion.  The exceptional two-torsion classes satisfy the same assertions.
\end{proof}

\subsection{The prime of norm thirteen}

Note that
\[
 \sqrt3\equiv-4\ \pmod{\mathfrak a}
 \ \ \t{and}\ \ \eps\equiv-2\equiv11
 \pmod{\mathfrak a}.
\]
The nonzero square residues in $\F_{13}$ are
$\{1,3,4,9,10,12\}$ and the fourth powers are $\{1,3,9\}$.
In particular $11$ and $7$ are nonsquares, while $10$ and $12$ are not
fourth powers.

\begin{lemma}\label{lem:thirteen}
For $e=a,2a$, we have
\[
 \im\alpha_e=\langle-1,e\rangle
 \ \ \t{and}\ \ |\im\alpha_e|=4.
\]
Moreover,
\[
 \im\alpha'_a\subseteq\{1,a\}\ \ \t{and}\ \ 
 \im\alpha'_{2a}\subseteq\{1,\eps a\}.
\]
\end{lemma}

\begin{proof}
A class $d$ in the first descent image is represented by the quartic
\[
 W^2=dU^4-\frac{e^2}{d}V^4,
\]
and a class in the dual image by
\[
 W^2=dU^4+\frac{4e^2}{d}V^4.
\]
After scaling locally, $U$ and $V$ are integral and at least one is a unit.
Write $e=ka$ with $k\in\{1,2\}$.

For the first quartic, take $d=\eps$.  If $U$ is a unit, reduction modulo
$\mathfrak a$ has nonsquare unit part $\eps$.  If $U$ is not a unit, then
$V$ is a unit; after removing the square $a^2$, the unit part is
$-k^2/\eps\equiv7k^2$, again a nonsquare.  Thus $\eps$ is absent from
$\im\alpha_e$.  The full rational two-torsion supplies
$\langle-1,e\rangle$, proving the first assertion.

For the dual quartic, the same two-unit-part argument excludes $d=\eps$.
When $e=a$ and $d=\eps a$, valuation parity first forces both $U$ and $V$ to be units, and a local solution would then force
\[
 (U/V)^4\equiv-4/\eps^2\equiv-1\equiv12\pmod{\mathfrak a},
\]
which is impossible.  Hence
$\im\alpha'_a\subseteq\{1,a\}$.  When $e=2a$ and $d=a$, valuation parity again forces both variables to be units, and a local
solution would force
\[
 (U/V)^4\equiv-16\equiv10\pmod{\mathfrak a},
\]
also impossible.  Since $[2a]=[\eps a]$, this gives
$\im\alpha'_{2a}\subseteq\{1,\eps a\}$.
\end{proof}

Lemma~\ref{lem:thirteen} and \eqref{eq:isogeny-rank} give
$\operatorname{rank}E_a(F)\le1$.  Lemma~\ref{lem:P0} therefore proves
\begin{equation}\label{eq:rank-a}
 \operatorname{rank}E_a(F)=1.
\end{equation}

\subsection{The quadratic twist and the dyadic covering}

The quadratic twist of $E_a$ by $2$ is $F$-isomorphic, after the standard
change of variables, to
\[
 E_{2a}:y^2=x^3-(2a)^2x.
\]
By Lemma~\ref{lem:thirteen}, $|\im\alpha_{2a}|=4$ and $\im\alpha'_{2a}\subseteq\{1,\eps a\}$.
It remains to exclude $\eps a$ from the dual image.

\begin{lemma}\label{lem:dyadic}
The quartic equation
\begin{equation}\label{eq:quartic-cover}
 W^2=(11+6\sqrt3)U^4+(80-32\sqrt3)V^4
\end{equation}
has no primitive solution over $F_{1+\sqrt3}$.
\end{lemma}

\begin{proof}
For $d=\eps a$ the dual covering is precisely
\eqref{eq:quartic-cover}, since
\[
 \eps a=11+6\sqrt3\ \ \t{and}\ \  \frac{16a}{\eps}=80-32\sqrt3.
\]
A primitive local solution can be scaled so that $U$ and $V$ are integral and at
least one is a unit.  It would therefore give a primitive solution modulo
$128$ in
\[
 R=(\Z/128\Z)[X]/(X^2-3).
\]
The finite sets of squares, fourth powers and unit fourth powers in $R$ have
respectively $1368$, $275$ and $256$ elements.  Exhausting the two charts in
which $U$ or $V$ is a unit shows that
\[
 (11+6X)U^4+(80-32X)V^4
\]
is never a square in $R$.  The accompanying program performs this exhaustive
calculation using integer arithmetic only.  Hence no primitive local solution
exists.
\end{proof}

It follows that $\im\alpha'_{2a}=\{1\}$, and formula
\eqref{eq:isogeny-rank} gives
\begin{equation}\label{eq:rank-twist}
 \operatorname{rank}E_{2a}(F)=0.
\end{equation}

\begin{lemma}\label{lem:common-biquadratic-carrier}
Let $k$ and $F_0$ be two distinct real quadratic subfields
of a biquadratic field $L_0$, and let $E/F_0$ be an elliptic
curve with full rational two-torsion.  Suppose
\[
 \operatorname{rank}E(L_0)=1,\qquad
 E(L_0)_{\rm tors}=E[2],\qquad
 2E(L_0)\subseteq E(F_0).
\]
Then $2E(L_0)=\Z P$ for some nontorsion $P\in E(F_0)$.
Every point $[n]P$ admits a half-point in $E(L_0)$, and
the abscissa of that half-point has exactly two coordinates
in any fixed $k$-basis of $L_0$.  The nontrivial automorphism
of $L_0/k$ restricts to the nontrivial automorphism of
$F_0/\Q$.
\end{lemma}
\begin{proof}
Multiplication by two kills $E(L_0)_{\rm tors}$, so
$2E(L_0)$ is torsion-free of rank one, hence infinite cyclic.
The stated inclusion places its generator in $E(F_0)$.
By its definition, each multiple of the generator is a
double of a point of $E(L_0)$.  The field $L_0/k$ has
degree two, and its nontrivial automorphism cannot fix the
distinct quadratic subfield $F_0$.
\end{proof}

\section{Rank stability and the cyclic carrier}\label{sec:stability}

Because $L=F(\sqrt2)$, the standard quadratic-twist rank decomposition gives
\[
 \operatorname{rank}E_a(L)
 =\operatorname{rank}E_a(F)+\operatorname{rank}E_{2a}(F).
\]
Equations \eqref{eq:rank-a} and \eqref{eq:rank-twist} yield
\begin{equation}\label{eq:rankstable}
 \operatorname{rank}E(L)=\operatorname{rank}E(F)=1.
\end{equation}

The prime $73$ splits completely in $L$.  The reduction calculation \eqref{eq:73orders}, together with the absence of $73$-torsion in the corresponding formal groups, yields that 
$ |E(L)_{\rm tors}|$ divides $4.$
Since $E[2]\subset E(F)$, we have
\begin{equation}\label{eq:torsion}
 E(L)_{\rm tors}=E[2].
\end{equation}

Let $\tau$ be the nontrivial automorphism of $L/F$.  The free quotient of $E(L)$ has rank one, so the involution $\tau$ acts on it by $+1$ or $-1$.  The image of the non-torsion point $P_0\in E(F)$ is nonzero and fixed by $\tau$; hence the action is $+1$.  Hence, for every $Q\in E(L)$, we have
\[
 \tau(Q)-Q\in E(L)_{\rm tors}=E[2].
\]
It follows that
$\tau(2Q)=2Q$,
so
\begin{equation}\label{eq:2EL}
 2E(L)\subset E(F).
\end{equation}
Lemma~\ref{lem:common-biquadratic-carrier}, applied with
$k=K$, $F_0=F$, and $L_0=L$, now gives
$2E(L)=\Z P$ for some $P\in E(F)$ and the two-coordinate
representation of its half-points.  This completes the proof of
Theorem~\ref{thm:carrier}.

\section{Eliminating fixed field coefficients}
\begin{lemma}\label{lem:coefficients}
For the real quadratic Galois field $K$, existential solvability
of a parameter-free finite polynomial system with coefficients in
$K$ reduces effectively to existential solvability of one polynomial
with integer coefficients over $K$.  More generally, if a distinguished
input $X$ is retained, the same construction preserves a relation
$\mathcal R(X)$ provided $\mathcal R$ is invariant under the nontrivial
automorphism of $K/\Q$.
\end{lemma}
\begin{proof}
For a parameter-free system $f_1=\cdots=f_r=0$ with
$f_i\in K[\mathbf X]$, put $G=\sum_{i=1}^r f_i^2$.
Since $K$ embeds in $\mathbb R$, $G(\mathbf x)=0$ if and only if all
$f_i(\mathbf x)=0$.  Let $\sigma$ be the nontrivial automorphism of
$K/\Q$.  Then $H(\mathbf X)=G(\mathbf X)G^\sigma(\mathbf X)$ has
rational coefficients.  If $H(\mathbf x)=0$, either
$G(\mathbf x)=0$ or $G^\sigma(\mathbf x)=0$; in the second case
$G(\sigma(\mathbf x))=0$.  Thus $H$ has a zero over $K$ exactly when
the original system does.  Clearing fixed rational denominators gives
an integer-coefficient polynomial.

If $X$ is a distinguished input, the second factor instead produces a
solution of the original system at $\sigma(X)$.  Hence the same product
is valid precisely when the relation defined by the system is stable
under $X\mapsto\sigma(X)$.  In the applications defining $\Z$, this
stability is checked explicitly: soundness gives
$\sigma(X)\in\Z$, and $\Z$ is fixed pointwise by $\sigma$.
\end{proof}

\section{The seven-contact pullback and the \texorpdfstring{$367$}{367}-adic tube}\label{sec:seven}
We use the discriminant-six quaternion algebra
\[
 B=(2,-3)_{\Q}=\Q\langle i,j:\,i^2=2,\ j^2=-3,\ ij=-ji\rangle
\]
and its maximal order
$O=\Z[i,(1+j)/2]$.  For $b\in B$, its reduced trace and reduced norm are given by
\[
 \operatorname{Trd}(b)=b+\bar b\ \ \t{and}\ \ 
 \operatorname{Nrd}(b)=b\bar b,
\]
where $\bar b$ is the quaternionic
conjugation of $b$.
The genus-zero positive-normalizer quotient
is equipped with its rational coordinate $j$ having elliptic orders
$4,6,2$ above $j=0,\infty,-16/27$ respectively.
The polarized moduli quotient comes from
\cite[Theorem~2.1]{BabaGranath}; its canonical quotient and intermediate
covering are given in \cite[Theorems~3.10 and~3.17]{BabaGranath}, and
the Hauptmodul normalization is \cite[Example~29]{LinYang}.
Here and below PEL means polarization--endomorphism--level structure.

\begin{lemma}
\label{lem:source-pel}
Let
 $B=(2,-3)_{\Q}$, $O=\Z[i,(1+j)/2]$,
 and $\mu=2j+ij$ with $\mu^2=-6$.
In the ordered $\Z$-basis
$1,i,(1+j)/2,(i+ij)/2$ of $O$, the alternating form
\[
 E_0(x,y)=\operatorname{Trd}(\mu^{-1}x\bar y)
\]
has matrix
\[
 \begin{pmatrix}
 0&0&-1&1\\0&0&1&-2\\1&-1&0&0\\-1&2&0&0
 \end{pmatrix}\]
 with determinant one.
The printed alternating form
$E_{\mathrm{BG}}(x,y)=\operatorname{Trd}(\bar x\mu y)$
in \cite[\S2.1]{BabaGranath} equals $6E_0$ on this
lattice, and has determinant $6^4$.  Hence $E_0$ is the
unimodular principal datum; these positive scalar multiples
have the same Rosati involution.  Rotger's polarization and
forgetful-map results \cite[\S4]{Rotger}, together with the
corrected integral form, identify the relevant polarized
forgetful quotient with the full Atkin--Lehner quotient of
\cite[Theorem~2.1]{BabaGranath}.
\end{lemma}
\begin{proof}
Using $i^2=2$, $j^2=-3$ and $ij=-ji$, direct multiplication
in the displayed integral basis gives the stated matrices.
The determinant of the $4\times4$ matrix is one.
The existence and uniqueness of the principal polarization with
this Rosati involution is also proved directly for false elliptic
curves over $\Z[1/6]$ in \cite[\S1]{Buzzard}; the following
linear-algebra argument fixes the exact normalization on our lattice.
The equality $E_{\mathrm{BG}}=6E_0$ follows from
$\mu^{-1}=-\mu/6$ and cyclic invariance of reduced trace.
We make the normalization step explicit. Fix the Rosati
involution $b\mapsto b^*=\mu^{-1}\bar b\mu$ on the left
$B$-action.  If $E$ is any rational alternating form with
this adjoint, the function $\ell(z)=E(1,z)$ satisfies
$E(x,y)=\ell(x^*y)$ and $\ell(z^*)=-\ell(z)$.
The $(-1)$-eigenspace of this orthogonal involution on $B$
is the one-dimensional rational subspace $\Q\mu$.
Consequently the compatible alternating forms constitute
a one-dimensional rational vector space, generated by $E_0$.
Thus $E=cE_0$ for $c\in\Q^\times$; the positivity
condition forces $c>0$.
On the same integral lattice $O$, a principal alternating form
has determinant one, hence
$1=\det(cE_0)=c^4\det E_0=c^4$ and $c=1$.
Thus $E_0$ is the \emph{unique} principal normalization of
this Rosati form on the specified regular lattice.
In particular, the factor six in the printed formula should
not be carried into the principal-polarized PEL functor.
This uniqueness identifies the relevant normalized polarized
complex moduli problem, on the same quaternionic lattice,
with the principal-type setting of Rotger
\cite[Theorem~4.4]{Rotger}; the forgetful quotient by the full
Atkin--Lehner group is then the one described in
\cite[Theorem~2.1]{BabaGranath}.
The finite \'{e}tale marking construction and its precise PEL
interpretation are proved in Lemma~\ref{lem:full-three-pel}
and Proposition~\ref{prop:pel-inertia}.
\end{proof}

\begin{remark} Theorems~2.1 and Corollary~2.3 of \cite{Buzzard} construct
a smooth proper fine integral model at the auxiliary level
$V_1(5)$, with \emph{maximal} local level at the primes of
quaternion discriminant $6$.  They do not themselves give a
full level-three model at the ramified coefficient prime.
After inverting $3$ and $5$, the full-five and full-three
marking sheaves are finite \'{e}tale over the appropriate fine
model.  The resulting complete level-three PEL functor and
its normalizer action are the additional constructions proved
in this paper.  The normalization of the discriminant-six
Hauptmodul uses \cite[Example~29]{LinYang}, while the
canonical-model and intermediate-quotient statements use
\cite[Theorem~3.10, Corollary~3.11, and Theorem~3.17]{BabaGranath}.
\end{remark}

\begin{lemma}
\label{lem:full-three-pel}
Let $K_0=\widehat O^{\times}$ and
$
 K_{\mathrm{lev},3}=\ker\bigl(K_0\longrightarrow (O/3O)^\times\bigr).
$
On the generic fibre, the Shimura curve of level $K_{\mathrm{lev},3}$ is the
canonical PEL model over $\Q$ attached to the datum
$(B,B,E_0)$, with $B$ acting on the left.  Over
$\Z[1/30]$, for a principally polarized QM surface
$(A,\iota,\lambda)$ with this Rosati datum, its full-three
level sheaf is
\[
 \mathscr L_3(A)=\operatorname{Isom}_{O}(O/3O,A[3]).
\]
This is exactly the full PEL \emph{similitude} level at three;
there is no additional determinant or pairing component to add.
For a normalizer transition represented by $n$ whose local
three-component belongs to $O_3^\times$, the incoming-frame
convention gives the right class of $n_3^{-1}$ on the same
module $O_3/3O_3$.  If $n_3$ is a division-algebra
uniformizer, the transition is instead between two
\emph{different} integral lattices, as specified in
Lemma~\ref{lem:PEL-lattice-transport}; in this case there is no
reduction of $n_3^{-1}$ as an endomorphism of $O_3/3O_3$.
The inertia element $\gamma=1+i+j$ used later has reduced
norm two, so its three-component is a unit and belongs to
the first case.
\end{lemma}
\begin{proof}
In Deligne's PEL description, take the faithful left
$B$-module $V=B$ and the alternating form $E_0$.  The group of
$B$-linear symplectic similitudes is the group of right
multiplications by $B^\times$; this is the quaternionic instance
of \cite[\S4.9]{Deligne}.  The complex PEL description, including similitude level structures,
is given in \cite[\S\S4.11--4.15]{Deligne}.  The construction of
the canonical model of the ambient Siegel moduli tower is in
\cite[\S\S4.16--4.17 and Theorem~4.21]{Deligne}.  We do
\emph{not} invoke Deligne's Proposition~4.20 as a direct
quaternionic canonical-model theorem: its statement appears in
the CM-type special setting of his \S4.18.
For this noncommutative quaternionic datum, the canonical
$\Q$-model and its discriminant-six Atkin--Lehner quotient are
identified directly by \cite[Theorems~2.1 and~3.10,
Corollary~3.11]{BabaGranath}, and the relevant fine
integral QM moduli model at prime-to-six auxiliary level
is \cite[Theorem~2.1 and Corollary~2.3]{Buzzard}.
The reflex field of the quaternionic Shimura datum is $\Q$.
The integral full-three scheme used here is not being attributed
to Deligne or Buzzard: after inverting $3$ it is the finite
\'{e}tale marking sheaf over the auxiliary fine integral PEL
model supplied at level five by \cite{Buzzard}.

We check explicitly that the displayed $O$-linear marking sheaf
is the complete PEL level.  The rational Tate module at three is
a one-dimensional left module over the division algebra
$B\otimes\Q_3$.  Hence every stable lattice is a principal left
lattice over the local maximal order, and
$A[3]\simeq O/3O$ as left $O$-modules.  On a rank-one left
$B\otimes\Q_3$-module, two non-degenerate alternating forms
inducing the same Rosati involution differ by a central scalar;
for principal integral forms that scalar is a $3$-adic unit.
After reduction modulo three, the polarization therefore
specifies the one-dimensional similitude line generated by
$\overline E_0$.

Every left-$O$ automorphism of the regular module $O/3O$ is
right multiplication by a unique unit $u\in(O/3O)^\times$.
For any lift of $u$ to the local maximal order,
\[
 E_0(xu,yu)
 =\operatorname{Trd}(\mu^{-1}xu\overline{yu})
 =\operatorname{Nrd}(u)E_0(x,y).
\]
Thus every $O$-linear marking is automatically a PEL
similitude, with multiplier $\operatorname{Nrd}(u)\bmod3$;
conversely every full PEL level marking is $O$-linear.  Both
multiplier components occur, since
$\operatorname{Nrd}(\gamma)=2\equiv-1\pmod3$.  This proves the
assertion about the whole level sheaf, including its two
similitude components.

Finally, the adelic canonical model carries the normalizer
Hecke action on the tower of levels \cite[\S\S4.14--4.17]{Deligne}.
For $n_3\in O_3^\times$ its action on the level-three marking is
precomposition by the reduction of $R_{n_3^{-1}}$.  For a
nonunit normalizer component, the source and target integral
Tate lattices differ; the transport in the next lemma defines
the relevant map between their three-torsion quotients.  It is
important not to identify this latter map with multiplication
by $n_3^{-1}$ on a fixed residue module.  Rotger's
fixed-underlying-$A$ presentation simultaneously transports
the embedding and the \emph{principal polarization}; its
explicit formula is used in
Lemma~\ref{lem:polarized-normalizer-comparison}
\cite[\S3.1]{Rotger}.
\end{proof}

\begin{lemma}
\label{lem:PEL-lattice-transport}
Let $D=B\otimes_\Q\Q_3$ and let $\mathcal O_D$ be its
maximal order.  Write $V=D$ as a left $D$-module, with
integral lattice $\Lambda=\mathcal O_D$.  For any
$n_3\in N_{D^\times}(\mathcal O_D)$, put
$\Lambda_{n_3}=\Lambda n_3^{-1}\subset V.$
The right-multiplication map
\[
 R_{n_3^{-1}}:\Lambda\longrightarrow\Lambda_{n_3}
\]
is an isomorphism of left $\mathcal O_D$-modules, and it induces
an isomorphism
\[
 \Lambda/3\Lambda\ \xrightarrow{\ \sim\ }\
 \Lambda_{n_3}/3\Lambda_{n_3}.
\]
Under this transport, the alternating form is carried to a
principal alternating form on $\Lambda_{n_3}$ by the factor
$\operatorname{Nrd}(n_3)$.  If $n_3\in\mathcal O_D^\times$,
this is the usual right action of $n_3^{-1}\bmod3$ on the
fixed residue module.  If $v_D(n_3)=1$, the same transport is
well-defined between the two lattices, although $n_3^{-1}$
need not be integral and has no right-multiplication reduction
on the fixed module $\Lambda/3\Lambda$.
\end{lemma}
\begin{proof}
The first map is, by its definition, a bijective left
$\mathcal O_D$-linear map and takes $3\Lambda$ bijectively
to $3\Lambda_{n_3}$, proving the quotient assertion.  If
$x'=xn_3^{-1}$ and $y'=yn_3^{-1}$, then cyclic invariance
of the reduced trace gives
\[
 E_0(x',y')=\operatorname{Nrd}(n_3)^{-1}E_0(x,y).
\]
Consequently, the form
$E_{n_3}=\operatorname{Nrd}(n_3)E_0$ on the transported
lattice pulls back to $E_0$, and is unimodular.  The
unit case fixes the lattice itself.  For a local division
uniformizer $\varpi$ with reduced norm of valuation one,
$\varpi^{-1}\notin\mathcal O_D$, whereas
$\mathcal O_D\varpi^{-1}$ remains a perfectly good
fractional left lattice.  In particular, the transport does
not make the invalid assertion that a division uniformizer
acts invertibly by multiplication on $\mathcal O_D/3\mathcal O_D$.
\end{proof}

\begin{proof}[Proof of Proposition~\ref{prop:seven}]
For $j(x)$ in \eqref{eq:pullback}, direct factorization gives
\begin{equation}\label{eq:factor-j}
 j(x)+\frac{16}{27}
 =\frac{16(x-1)(x+8)^2}{27(x-4)^3}.
\end{equation}
The three branch values are therefore pulled back as follows:
\[
\begin{array}{c|c|c|c}
 x & j(x) & \text{local degree}&\text{residual orbifold order}\\ \hline
 \infty&0&1&4\\
 0&0&2&2\\
 4&\infty&3&2\\
 1&-16/27&1&2\\
 -8&-16/27&2&1
\end{array}
\]
The associated quartic polynomial
\begin{equation}\label{eq:Pquartic}
 P(z)=\frac{(a_*z+b_*-q_*z^2)(a_*z+b_*-4q_*z^2)}{4q_*^2}
\end{equation}
satisfies the exact quadratic-class identity
\begin{equation}\label{eq:squareclass}
 -(27j(x(z))+16)=-P(z)
 \left(\frac{8(x(z)+8)}{z^2(x(z)-4)^2}\right)^2.
\end{equation}
Indeed, $27j(x)+16=16(x-1)(x+8)^2/(x-4)^3$.
Thus, whenever the identified positive-normalizer quadratic
subcover is defined by the square class of $-(27j+16)$,
its pullback is $K(\sqrt{-P(z)})$.

The identities
\[
 a_*^2+4q_*b_*=1\ \ \t{and}\ \ 
 a_*^2+16q_*b_*
 =\left(\frac{q_*^2+6q_*-3}{q_*^2+3}\right)^2
\]
show that the inverse images of $x=0,1,4$ are rational.
Specifically, the inverse image of $x=0$ consists of $\infty$
and the third element of \eqref{eq:contacts}; the inverse images
of $x=1$ are the fourth and fifth elements, and the inverse images
of $x=4$ are the sixth and seventh elements. At $z=0$ the map
$x(z)$ has a double pole. The order-four orbifold point $x=\infty$
therefore pulls back with residual order
$4/\gcd(4,2)=2$. All remaining listed contacts have residual
order two. This proves the claimed seven-point list.

Outside $S$, the branch sections are disjoint and the cover is tame.
The completed local cover at a residual order-two contact has equation
$u^2=t$ after an unramified extension and multiplication of $t$ by
a local unit. Thus a positive contact of valuation $n$ causes
ramification of index $2/\gcd(2,n)$; it is unramified if and only
if $n$ is even. At $z=0$, a contact of order $n$ is sent to the
original order-four branch with contact order $2n$, yielding
$4/\gcd(4,2n)$, with the same parity criterion.
No other branch points occur outside the fixed bad primes.

The prime $367$ splits in $K$ since $367\equiv7\pmod8$.
At either prime ideal $\mathfrak q_*\mid367$, we have
\[
 a_*\equiv-1\pmod{\mathfrak q_*}\ \ \t{and}\ \ b_*\equiv122\pmod{\mathfrak q_*}.
\]
For $z\equiv1\pmod{\mathfrak q_*}$, the numerator of $x(z)$ is a
unit, and its denominator has valuation one. Hence
$v_{\mathfrak q_*}(x(z))=-1$. Since $x(z)-4$ has the same
valuation, we have
\[
 v_{\mathfrak q_*}(j(x(z)))=2(-1)-3(-1)=1.
\]
\end{proof}

We next record the characteristic-zero normalizer calculation that
identifies the source of the order-four projective action.

\begin{lemma}\label{lem:gamma}
Let $\gamma=1+i+j\in O$. Then
\[
 \operatorname{Trd}(\gamma)=2,\qquad
 \operatorname{Nrd}(\gamma)=2,\qquad
 \gamma^2=2(i+j),\qquad \gamma^4=-4.
\]
The element $\gamma$ normalizes $O$, and its class in
$B^\times/\Q^\times$ has exact order four. Its right action
on the quaternionic three-torsion module has projective order four.
\end{lemma}
\begin{proof}
The first four identities follow by expansion using $ij=-ji$.
In the integral basis
$1,i,(1+j)/2,(i+ij)/2$ of $O$, the inner automorphism
$a\mapsto\gamma a\gamma^{-1}$ is
represented by
\[
 \begin{pmatrix}
 1&2&-1&0\\0&-1&1&1\\0&-4&3&0\\0&-2&1&-1
 \end{pmatrix}.
\]
Its determinant equals one, and its fourth power equals the identity.
Therefore $\gamma O\gamma^{-1}=O$. The element $\gamma^2$
is not central, so the projective order is four.
Moreover $(i+j)^2=-1$, so the fixed point of $\gamma$
corresponds to a discriminant-$-4$ CM embedding of
$\mathbb Q(\sqrt{-1})$. Its reduced characteristic polynomial
is $X^2-2X+2$; modulo three the ratio of its roots has
order four in $\F_9^\times$.
Thus the right-multiplication action on $O/3O$ has projective
order four after passing to a splitting field of $B\otimes\Q_3$.
\end{proof}

\begin{lemma}\label{lem:polarized-normalizer-comparison}
Let $(A,\iota,\lambda)$ be a principally polarized geometric QM
surface with the fixed Rosati datum $(O,\mu)$, and let
$\alpha:O/3O\xrightarrow{\sim}A[3]$ be an $O$-linear marking.
Put $\gamma=1+i+j$, $u_\gamma=\iota(\gamma)$ and
$N=\operatorname{Nrd}(\gamma)=2$.  Then the positive
Atkin--Lehner transition represented by $\gamma$ has the
following fixed-underlying-surface description:
\begin{equation}\label{eq:polarized-normalizer-data}
 \iota_\gamma(a)=\iota(\gamma^{-1}a\gamma),\qquad
 \lambda_\gamma=N^{-1}u_\gamma^*\lambda,\qquad
 \alpha_\gamma(x)=u_\gamma[3]^{-1}\alpha(x\gamma^{-1}).
\end{equation}
Here $\lambda_\gamma$ denotes the primitive polarization
class, which is integral and principal, and $u_\gamma[3]^{-1}$
is the inverse of the automorphism induced on $A[3]$.
The new triple satisfies the same abstract Rosati datum as the
old one.  For the induced incoming frame, the identity
\begin{equation}\label{eq:polarized-normalizer-tate}
 \alpha^{-1}u_\gamma[3]\alpha_\gamma
    =R_{\gamma^{-1}}\quad\text{on }O/3O
\end{equation}
holds on the actual three-torsion module.  Moreover,
\eqref{eq:polarized-normalizer-data} represents the right
normalizer action on the \emph{canonical} full-three PEL tower;
no independent choice of the incoming marking is involved.
\end{lemma}
\begin{proof}
Rotger's modular formula for positive Atkin--Lehner
operators \cite[\S3.1]{Rotger} simultaneously conjugates the
quaternionic action and replaces the polarization by its
normalized pullback.  We check its integrality explicitly in
the present order.  Conjugation gives
\[
 \mu_\gamma:=\gamma^{-1}\mu\gamma
           =3i+4j-2ij\in O\]
           and hence  
 $\mu_\gamma^2=-6.$
The normalized Riemann form is
\[
 N^{-1}E_0(\gamma x,\gamma y)
 =\operatorname{Trd}\bigl(\mu_\gamma^{-1}x\bar y\bigr).
\]
In the ordered basis $1,i,(1+j)/2,(i+ij)/2$,  its matrix is
\[
 \begin{pmatrix}
 0&2&-2&-1\\ -2&0&-3&-4\\ 2&3&0&2\\ 1&4&-2&0
 \end{pmatrix}
\]
with determinant one.
Thus the transported form is integral and unimodular.
Positivity follows because it is the normalized pullback of
$\lambda$ by the isogeny $u_\gamma$.  As
$u_\gamma\iota_\gamma(a)=\iota(a)u_\gamma$, the Rosati
compatibility with the original abstract involution is
preserved.  Further $\deg u_\gamma=N^2=4$, so
$u_\gamma[3]$ is invertible.  Formula
\eqref{eq:polarized-normalizer-data} is $O$-linear, and
\eqref{eq:polarized-normalizer-tate} is immediate by
composition.  The full level here records the Weil pairing
up to similitude: the displayed frame changes the pairing
by $N^{-2}$, which is $1$ modulo $3$ for $N=2$.

It remains to identify this formula with the geometric
normalizer operator rather than with an arbitrarily
constructed marking.  In the complex adelic uniformization,
write a marked QM point as $[h,g_f]$.  The rational lattice
for $g_f$ is $I_{g_f}=B\cap\widehat O g_f$, and its torsion
trivialization is $b\mapsto b g_f^{-1}$.  In particular, for
the right operator $g_f\mapsto g_f\gamma_f^{-1}$, the
trivialization is necessarily transported by right
multiplication by $\gamma_f$; on the incoming framing this
is $R_{\gamma_3^{-1}}$, with the choice of right-action
inverse fixed by the convention.  This is the explicit
compatibility (4.1) in the proof of
\cite[Theorem~4.2]{Molina}.  Since $N=2$ is a $3$-adic unit,
the lattice comparison is an isomorphism on genuine
three-torsion.  Combining it with Rotger's fixed-surface
realization gives precisely
\eqref{eq:polarized-normalizer-data}.  Hence the formula is
the canonical normalizer action, including its polarization
and integral level data.
\end{proof}

\begin{proposition}
\label{prop:pel-inertia}
Let $\mathfrak q$ be a finite prime of $K$ of residue characteristic
$p>144$ and residue cardinality $q=N\mathfrak q\equiv-1\pmod4$.
Suppose that $v_{\mathfrak q}(j)=1$.  For the marked PEL cover,
local inertia at this point is the effective order-four right class
of $\gamma=1+i+j$ (up to inverse).  If the corresponding fibre
is non-CM with $\operatorname{End}^0_{\overline K}(A)=B$,
this is the projective inertia action on the specific normalized
two-dimensional representation $r_3$, and its residual local
representation is absolutely irreducible.  If in addition
$q\equiv1\pmod3$, then the residual restriction to
$G_{K(\zeta_3)}$, where $\zeta_3$ is a primitive cube root
of unity, is absolutely irreducible.
\end{proposition}
\begin{proof}
The full-level PEL curve is taken with the explicit
similitude marking sheaf of Lemma~\ref{lem:full-three-pel}.
For an integral construction one may first work with a
full-five and full-three marking over $\Z[1/30]$.
Buzzard's smooth proper auxiliary-level model
\cite[Theorem~2.1 and Corollary~2.3]{Buzzard} supplies
this scheme: full-five level is finer than his $V_1(5)$
level, and adjoining full-three level is finite \'{e}tale
away from three.  Quotient by the full-five relabeling group
after retaining the full-three marking.  This action is free:
a geometric $O$-linear polarized automorphism preserving the
three-marking has finite order and lies in
$1+3\operatorname{End}_{\Z_3}(T_3(A))$, a torsion-free group;
thus it is the identity, including at CM points.  The
quotient is consequently the fine full-three PEL model,
smooth and proper over $\Z[1/30]$.  By the canonical Shimura
model its generic fibre is the level-$K_{\mathrm{lev},3}$
curve, and its normalizer Hecke automorphisms are defined
over $\Q$ \cite{Deligne,Rotger}.  This identifies the
integral marked fibre with the same canonical adelic level
on which the Galois and normalizer actions commute.

Consider its geometric fibre over the $j$-line.  By
Lemma~\ref{lem:Gamma} the unbranched degree is $144$ and
the deck group is the specified $\Gamma$.  The
Baba--Granath canonical quotient has an order-two elliptic
point above $j=0$ in the original norm-one Shimura curve,
and the corresponding $w_2$-involution fixes that point
\cite{BabaGranath}; its CM field is $\Q(i)$.  In the
adelic double quotient the element
$\gamma=1+i+j$ has this fixed point $h_0$ and
\[
 [h_0,1]\cdot\gamma_f^{-1}
   =[h_0,\gamma_f^{-1}]
   =[\gamma h_0,1]=[h_0,1].
\]
This is a fixed point in the actual level-three Shimura
tower, not merely in the unmarked Atkin--Lehner quotient.
Since $\gamma^2=2(i+j)$, its square represents the
nontrivial order-two elliptic stabilizer; $\gamma$ maps to
$w_2$ on the unmarked curve.  Thus it gives the entire
cyclic order-four stabilizer.  There is no larger stabilizer:
the norm-one elliptic stabilizer here has order two and the
stabilizer of this point in the four-element Atkin--Lehner
group is $\langle w_2\rangle$ (the other involutions
exchange the two points above $j=0$).  The action of
$\gamma$ has order four on the genuine three-level module
by Lemma~\ref{lem:gamma}.  Consequently the inertia
subgroup of the geometric $144$-sheeted cover above $j=0$
is exactly $\langle\gamma\rangle\simeq C_4$, up to
conjugacy and choice of generator.

At the indicated prime $\mathfrak q$, $p>144$ implies
$p\nmid30|\Gamma|$.  The smooth integral full-three model
and the distinct branch sections $j=0,-16/27,\infty$
therefore give a tame normalizer cover with good reduction.
The tame-local structure (or the generalized Abhyankar
lemma \cite[Expos\'{e}~XIII]{SGA1}) identifies its completed
local extension, after strict henselization, with
\[
        j=u w^4,\qquad u\ \text{a unit}.
\]
If $v_{\mathfrak q}(j)=1$, the corresponding Kummer
extension has ramification index four.  Hence actual
arithmetic inertia on a chosen geometric full-three fibre
maps onto the entire cyclic group generated by the deck
class of $\gamma^{\pm1}$.

We next identify this deck class with the \emph{same}
normalized Tate representation.  Let $y=(A,\iota,\lambda,
\alpha)$ denote a geometric lift of the non-CM fibre.
For $g$ in local inertia, $g(y)$ is obtained from $y$ by
the corresponding normalizer element in the canonical
PEL tower.  The identification of these two PEL points
provides an isomorphism from $({}^gA,{}^g\iota,
{}^g\lambda,{}^g\alpha)$ to the transported data of
Lemma~\ref{lem:polarized-normalizer-comparison}.
Composing with $u_\gamma$ gives a $B$-linear
quasi-isogeny $\mu_g:{}^g A\to A$.  Its action on the
three-primary torsion, with the incoming-frame convention,
therefore satisfies
\[
 \bigl[\alpha^{-1}\mu_g[3]{}^g\alpha\bigr]
      =\bigl[R_{\gamma^{\pm1}}\bigr].
\]
The equality is forced by the geometric right Hecke action:
Rotger's formula transports $\lambda$ as well as $\iota$,
while Molina's identity (4.1) transports the \emph{actual}
torsion frame.  Hence no additional projective normalizer
factor can be inserted.  If a different $B$-linear
quasi-isogeny is chosen, its quotient with $\mu_g$
commutes with $B$, and is a rational scalar because
$\operatorname{End}^0_{\overline K}(A)=B$.  The projective
action is thus independent of all these choices.

In the construction of Remark~\ref{rem:hom}, the
normalizing degree factors and the scalar cochain $\beta$
are central.  They disappear in projective space.  Splitting
$B\otimes\Q_3$ over $E_3$ and taking the two-dimensional
Morita multiplicity space therefore gives exactly
\[
 \mathbb P r_3(I_{\mathfrak q})
       =\langle[R_{\gamma^{\pm1}}]\rangle.
\]
Since $\gamma$ is a unit at three, the comparison is
integral and reduction modulo the maximal ideal is
legitimate.  Its polynomial is $X^2-2X+2$; its two roots
modulo three lie in $\F_9$ with ratio of order four.
Thus the residual inertia has two distinct eigenlines.
For $q\equiv-1\pmod4$, Frobenius acts on the tame inertia
generator by $\tau\mapsto\tau^q=\tau^{-1}$, so it
exchanges these eigenlines.  This proves absolute
irreducibility.  When also $q\equiv1\pmod3$, the local
cyclotomic character is trivial on the whole decomposition
group and that group lies in $G_{K(\zeta_3)}$; the same
argument proves the final assertion.
\end{proof}

\section{A controlled scalar cochain over \texorpdfstring{$\Q(\sqrt2)$}{Q(sqrt(2))}}\label{sec:cochain}
In this section all roots of unity are constant coefficient
modules. We use $\operatorname{Br}(K)[2]$ for the two-torsion
in the Brauer group of $K$, $\operatorname{inv}_v$ for its local
invariant at a finite prime ideal $v$, and $\mathfrak f_\chi$ for the
conductor of a finite-order Hecke character $\chi$.

\begin{lemma}[Dyadic unit interpolation]\label{lem:dyadic-interp}
Let $F_2=\Q_2(\sqrt2)$, $\pi=\sqrt2$ and $u=1+\pi$.
For a root of unity $\omega$ of two-power order and
$s\in\{\pm1\}$ there is a continuous finite-order character
$\theta:F_2^\times\to\mu_{2^\infty}$ satisfying
\[
 \theta(u)=\omega\ \ \t{and}\ \ \theta(-1)=s.
\]
If $\omega$ has order at most $2^r$, one can choose $\theta$
with order at most $2^{\max(2,r+1)}$ and conductor exponent
at most $5+2\max(r-1,0)$.
\end{lemma}
\begin{proof}
Let $U_n=1+\pi^n\OO_{F_2}$, where $\OO_{F_2}$ is the
valuation ring of $F_2$. We have
\[
 u^2=3+2\pi,\qquad u^4=17+12\pi\ \ \t{and}\ \ 
  v_\pi(u^4-1)=5.
\]
The $2$-adic logarithm gives an isomorphism
$
 \log:U_5\xrightarrow{\sim}\pi^5\OO_{F_2}$
of rank-two $\Z_2$-modules.
The vector $\log(u^4)$ is primitive since its valuation is
exactly five. Consequently, $ U_5=\overline{\langle u^4\rangle}\times V$
and $V\simeq\Z_2$.
Give $U_5$ the character taking $u^4$ to $\omega^4$ and
trivial on $V$. Since
$|\OO_{F_2}^\times/U_5|=16$, $u$ has order four
modulo $U_5$, and $-1$ does not lie in the subgroup it generates,
the group $\overline{\langle-1,u\rangle}U_5$ has index two
in $\OO_{F_2}^\times$. The above character extends to this
subgroup by the prescribed values at $-1$ and $u$.
Extend once more across the index-two inclusion, choosing a
square root of the necessary value, and set $\theta(\pi)=1$.
This gives a character of $F_2^\times$. The extension across an
index-two subgroup increases its order by at most a factor two.
If $\operatorname{ord}(\omega)\le2^r$, the character on $U_5$
annihilates $U_{5+2\max(r-2,0)}$ by the $2$-adic logarithm.
The extension does not change this restriction, so its conductor
exponent is bounded by $5+2\max(r-1,0)$, as asserted.
\end{proof}

\begin{proof}[Proof of Theorem~\ref{thm:cochain}]
It is well known that  $\OO_K=\Z[\sqrt2]$ is a principal ideal domain.
 Its unit group is
\[
 \OO_K^\times=\{\pm(1+\sqrt2)^n:n\in\Z\}.
\]

For each prime ideal of odd residue characteristic $v=\mathcal P$ with local invariant
$\operatorname{inv}_v[\varepsilon]=1/2$, choose a tame character
of $K_v^\times$ whose restriction to the cyclic residue-field
unit group has order $2^{v_2(N\mathcal P-1)}$ and sends
$-1$ to $-1$. At each odd-residue-characteristic prime ideal of invariant zero, choose
the trivial unit character. At each real embedding choose the sign
character if and only if the local invariant is nonzero.
Let $\omega$ be the inverse of the product of these local
unit characters evaluated at $1+\sqrt2$, and let $s$ be the
inverse of their product evaluated at $-1$.
Brauer reciprocity implies that $s$ is exactly the sign
required for the dyadic local invariant.
Apply Lemma~\ref{lem:dyadic-interp} to choose the character
at the unique prime ideal above $2$ with values $\omega$ and $s$.
The product of all local unit characters is now trivial on the
generators $-1,1+\sqrt2$ of $\OO_K^\times$.
Since $h_K=1$, global class field theory extends these unit
characters to a finite-order Hecke character $\chi$ on $G_K$.

The exact sequence of constant discrete modules
\[
 1\longrightarrow\{\pm1\}\longrightarrow\mu_\infty
 \xrightarrow{z\mapsto z^2}\mu_\infty\longrightarrow1
\]
induces a boundary $H^1(G_K,\mu_\infty)\to
H^2(G_K,\{\pm1\})$. By local reciprocity, the local
obstruction to a character square root is detected by its
value on $-1$. Hence the boundary of $\chi$ has the same
local Brauer invariants as $[\varepsilon]$ at every local completion (including the real embeddings).
The global-to-local map on $\operatorname{Br}(K)$ is injective,
so the two classes coincide. Therefore a continuous
pointwise square root of $\chi$ can be corrected by a sign
cochain to produce $\beta$ with
$\beta^2=\chi$ and $\delta\beta=\varepsilon$.

The odd local characters have orders dividing the product
$B_{T_0}$. The values prescribed for the dyadic character
have order dividing $2B_{T_0}$. Lemma~\ref{lem:dyadic-interp}
then permits global character order $m\le4B_{T_0}$.
It also gives a dyadic conductor exponent at most
$5+2\log_2 B_{T_0}$; at each odd prime in $T_-$
the character is tame and has conductor exponent one.
Since $N\mathcal P_2=2$, we have
\[
 N_{K/\Q}(\mathfrak f_\chi)
 \le2^{5+2\log_2 B_{T_0}}
 \prod_{\mathcal P\in T_-}N\mathcal P
 =32B_{T_0}^2\prod_{\mathcal P\in T_-}N\mathcal P.
\]

Because $\varepsilon$ is inflated from $\Gal(L_0/K)$, the
restriction $\beta|_{G_{L_0}}$ is a character, invariant under
conjugation by $G_K$. At every prime outside $T_0$ this
character is at most quadratic on inertia. Choose a principal
integral ideal whose odd prime support is precisely the finite
set of prime ideals where this inertia character is nontrivial.
The associated quadratic character has the same nontrivial
inertia restriction at each such odd-residue-characteristic prime ideal; twisting $\beta$
by it eliminates those restrictions without changing $\beta^2$
or $\delta\beta$. The resulting kernel field $L'/K$ is Galois,
unramified outside $T_0$, and
\[
 [L':K]\le 2m[L_0:K]\le8B_{T_0}[L_0:K].
\]
\end{proof}

\begin{remark}\label{rem:hom}
A full-three level structure makes the three-torsion of the
chosen fibre and all its Galois conjugates constant over a
common splitting field. The injection
\[
 \operatorname{Hom}_{\overline K}(A_1,A_2)/3
 \hookrightarrow\operatorname{Hom}(A_1[3],A_2[3])
\]
shows that the Galois action on the Hom lattice is trivial
modulo three. Since its image is finite and the principal
congruence subgroup of level three in $\mathrm{GL}_r(\Z)$
is torsion-free, all geometric homomorphisms descend to that
field. Complex homomorphisms of abelian varieties over
$\overline\Q$ already descend to $\overline\Q$, by the
Zariski density of geometric torsion points.

For each $g$, first transport the $B$-label on ${}^gA$ by the
same normalizer transition by which $g$ moves the canonical
full-three PEL point.  Choose
$\mu_g:{}^gA\to A$ to be $B$-linear for these transported
labels.  Any two such choices differ by an element of
$\operatorname{End}^0_B(A)=\Q$.  If
$\mu_g\,{}^g\mu_h=c(g,h)\mu_{gh}$ with
$c(g,h)\in\Q^\times$, then degree multiplicativity on an
abelian surface gives
$\deg(c(g,h)\mu_{gh})=|c(g,h)|^4\deg(\mu_{gh})$.
Consequently, after putting
$a_g=(\deg\mu_g)^{-1/4}>0$ and $\nu_g=a_g\mu_g$, the remaining
discrepancy is
$c(g,h)/|c(g,h)|\in\{\pm1\}$.  This is precisely the sign
cocycle $\varepsilon(g,h)$ used above.

After selecting $\beta$,
$R_\ell(g)=\beta(g)^{-1}a_g\mu_g\circ g$ is a true
$\ell$-adic representation, and on the kernel field $L'$
it equals the actual Tate action of $A$.  Let $E_\ell/\Q_\ell$
be a finite coefficient extension splitting
$B\otimes\Q_\ell$, choose
$B\otimes E_\ell\simeq M_2(E_\ell)$ and put $e=e_{11}$.
Because $R_\ell$ commutes with the transported left $B$-action,
\[
 r_\ell:=R_\ell\big|_{e(V_\ell(A)\otimes E_\ell)}
\]
is a two-dimensional representation and
$R_\ell\otimes E_\ell\simeq r_\ell^{\oplus2}$ as a
representation of the centralizer.  This is the actual
multiplicity representation used throughout the paper.  In
particular, its projectivization is the right normalizer class
identified in Proposition~\ref{prop:pel-inertia}; no further
choice of matrix units changes that projective class.
The Weil pairing gives the determinant identity
$\det r_\ell=\chi_\ell\chi^{-1}$ in compatible
normalizations, where $\chi_\ell$ denotes the cyclotomic
character.
\end{remark}

\begin{corollary}\label{cor:beta-quant}
Let $L_0/K$ be a splitting field of the full-three fibre.  It may be
chosen Galois with $[L_0:K]\le144!$: by Lemma~\ref{lem:Gamma} the
generic fibre has $144$ geometric points, so the Galois closure of
the corresponding degree-$144$ finite \'{e}tale algebra has degree at
most $144!$ (and special fibres only lower this bound).  Suppose its
ramification lies in
$T_0=S\cup\operatorname{Supp}(M_{\mathcal C,S}(z))$.
Under the hypotheses of Theorem~\ref{thm:cochain}, there is a
constant $C_S$ independent of $z$ such that
\[
 [L':K]\le C_S M_{\mathcal C,S}(z),\qquad
 N_{K/\Q}(\mathfrak f_\chi)\le
 C_S M_{\mathcal C,S}(z)^3.
\]
\end{corollary}
\begin{proof}
By $2^{v_2(N\mathcal P-1)}\le N\mathcal P$,
$B_{T_0}\le C_S M_{\mathcal C,S}(z)$.
Likewise the product of norms over the odd ramified
prime ideals in $T_-$ is at most $C_S M_{\mathcal C,S}(z)$.
Use the bounds of Theorem~\ref{thm:cochain} and the
uniform degree bound $144!$.
\end{proof}

\section{The automorphic and Jacobian height interface}\label{sec:height}
We write $h_F(A)$ for the stable Faltings height of an abelian
variety $A$ over a number field.  In this section $r_3$ denotes
the two-dimensional summand
of the normalized QM representation from
Remark~\ref{rem:hom}, after extending the coefficient field
to split $B\otimes\Q_3$. The fixed scalar character $\chi$
is the square of the cochain $\beta$; the subscript on
$\chi_3$ refers instead to the cyclotomic character.

\subsection{Fixing the same representation in cohomology}
\label{subsec:cohom-normalization}
We write $\mathrm{Fr}^{\rm ar}_v$ for arithmetic Frobenius and
$\mathrm{Fr}^{\rm geo}_v=(\mathrm{Fr}^{\rm ar}_v)^{-1}$ for its
geometric inverse.  We use arithmetic Frobenius on covariant
Tate modules.  If $X/K$ is smooth and proper and
$J=\operatorname{Pic}^0_{X/K}$, the canonical
Galois-equivariant comparison is
\begin{equation}\label{eq:etale-tate-dual}
 V_\lambda(J)\simeq
 H^1_{\mathrm{\acute et}}(X_{\overline K},E_\lambda)^\vee,
\end{equation}
where $E_\lambda$ is a coefficient field splitting the QM algebra,
$V_\lambda=T_\ell\otimes_{\Z_\ell}E_\lambda$, and the right side
is the $E_\lambda$-linear dual with the contragredient Galois action.
For a disconnected geometric curve, the statement is applied to the
product of the Jacobians of its connected components.

\begin{lemma}[Frobenius-normalized cohomological matching]
\label{lem:cohom-normalization}
Let $r:G_K\to\mathrm{GL}_2(E_\lambda)$ be a continuous,
semisimple representation, and suppose
\[
 \det r=\chi_\ell\,\vartheta
\]
for a finite-order character $\vartheta$. Suppose a smooth proper
Shimura curve $X_U/K$ has a two-dimensional semisimple constituent
$W$ of $V_\lambda(J_U)$ such that for all but finitely many
finite prime ideals $v\nmid\ell$ one has
\begin{equation}\label{eq:matching-polynomial}
 \begin{aligned}
 \det(X-r(\operatorname{Frob}_v))
 &=X^2-a_vX+(Nv)\vartheta(\operatorname{Frob}_v),\\
 \det(X-W(\operatorname{Frob}_v))
 &=X^2-a_vX+(Nv)\vartheta(\operatorname{Frob}_v).
 \end{aligned}
\end{equation}
Then $W\simeq r$ as representations of $G_K$. In particular
this identifies the \emph{same} representation, with no
additional finite-order twist or Tate twist.  The comparison is
made on the Tate side \eqref{eq:etale-tate-dual}, so that a dual
present in an $H^1$ convention has already been taken into account.
\end{lemma}
\begin{proof}
The characteristic polynomials in
\eqref{eq:matching-polynomial} imply equality of traces on a
dense set of Frobenius conjugacy classes. Chebotarev density
and continuity give equality of characters on $G_K$; the
Brauer--Nesbitt theorem and semisimplicity imply $W\simeq r$.
A Tate twist by a nonzero integer would change the determinant
by $\chi_\ell^{2n}$ and violate the displayed determinant
identity. A putative finite-order twist must pass the trace
comparison, and if it does, Brauer--Nesbitt already identifies
it with $r$. Thus no separate twist is left to be introduced.
\end{proof}

\begin{corollary}
\label{cor:cohom-hom}
Let $A/L'$ be an abelian surface whose quaternionic Tate module
becomes $r^{\oplus2}$ over $G_{L'}$ after extending coefficients.
Assume the hypotheses of Lemma~\ref{lem:cohom-normalization}
with $r=r_3$, and that $J_U$ is defined over $K$.
Then $\operatorname{Hom}_{L'}(A,J_U)\ne0$.
If $A$ is geometrically simple, it is an $L'$-isogeny factor
of $J_U$.
\end{corollary}
\begin{proof}
The inclusion $W\subset V_\lambda(J_U)$ and
$V_\lambda(A)|_{G_{L'}}\simeq r_3^{\oplus2}|_{G_{L'}}$
give a nonzero $G_{L'}$-equivariant map between the Tate modules
after extending coefficients. In particular, only one
copy of $r_3$ in $V_\lambda(J_U)$ is needed: a nonzero map
from the direct sum of two copies to that constituent
already gives a nonzero element of
$\operatorname{Hom}_{G_{L'}}(V_\lambda(A),V_\lambda(J_U))$. The invariants of the finite-dimensional
Hom space commute with this flat coefficient extension, so the
Hom space over $\Q_\ell$ is already nonzero. Faltings'
homomorphism theorem \cite{Fal} produces an actual nonzero $L'$-morphism
$A\to J_U$. When $A$ is geometrically simple its kernel is
finite, and Poincar\'e complete reducibility gives the
isogeny-factor assertion.
\end{proof}

\begin{lemma}
\label{lem:skinner-supercuspidal}
Let $\Pi$ be a cohomological cuspidal Hilbert representation
whose attached three-adic representation, in the convention fixed
above, is isomorphic to $r_3$.
Suppose $v_0\nmid3$ is a finite prime ideal at which
$\operatorname{WD}(r_3|_{G_{K_{v_0}}})$ has irreducible
two-dimensional Weil part.  Then $\Pi_{v_0}$ is
supercuspidal.  This conclusion does not assume in advance
that $\Pi$ has any discrete-series local component.
\end{lemma}
\begin{proof}
Skinner \cite[(1)]{Skinner} states the
local--global compatibility of the attached $p$-adic Galois
representation with the local Langlands parameter at every
finite $v\nmid p$, without condition (C) of
\cite[Theorem~2.1]{Saito}.  Applying this with $p=3$ and
$v=v_0$, we get
\[
 \operatorname{WD}(r_3|_{G_{K_{v_0}}})^{\mathrm{Fr\text{-}ss}}
 \simeq\iota\operatorname{Rec}_{v_0}
             (\Pi_{v_0}\otimes|\det|_{v_0}^{-1/2}).
\]
The harmless unramified normalization twist leaves the
irreducibility of the Weil representation unchanged.  By the
local Langlands classification for $\mathrm{GL}_2$ an
irreducible two-dimensional Weil parameter corresponds to a
supercuspidal representation.  Thus $\Pi_{v_0}$ is
supercuspidal.  Since Skinner's compatibility was used first,
there is no circular appeal to the Carayol--Saito realization.
\end{proof}

\begin{proposition}
\label{prop:carayol-tate}
Let $K$ be totally real of degree two, and let $D/K$ be a
quaternion algebra split at exactly one real embedding and ramified
at one fixed nonzero prime ideal $v_0$.  Let $\Pi$ be a parallel-weight-two
cuspidal representation of $\mathrm{GL}_2(\mathbb A_K)$ such that
$\Pi_{v_0}$ is discrete series.  Let $\rho^S_{\Pi,3}$ denote the
$\Pi$-isotypic two-dimensional Galois representation occurring in
$H^1$ in Saito's cohomological geometric-Frobenius convention, and put
\[
        \rho_{\Pi,3}:=(\rho^S_{\Pi,3})^\vee.
\]
Then, for a sufficiently small level $U\subset D(\mathbb A_{K,f})^\times$
and a finite coefficient extension $E_\lambda/\Q_3$, there is a
$G_K$-stable two-dimensional constituent
\begin{equation}\label{eq:carayol-tate}
 W_\Pi\subset V_3(J_U)\otimes_{\Q_3}E_\lambda,
 \qquad W_\Pi\simeq(\rho^S_{\Pi,3})^\vee=\rho_{\Pi,3},
 \qquad J_U=\operatorname{Pic}^0(X_U).
\end{equation}
The representation is attached to the \emph{same} Hecke eigensystem
$\Pi$, with no unspecified Tate or finite-order twist.
\end{proposition}
\begin{proof}
By Jacquet--Langlands \cite{JacquetLanglands} or \cite{Carayol}, the
finite discrete-series component permits transfer of $\Pi$ to the
quaternion algebra $D$.  This hypothesis is needed \emph{before}
applying the compact Shimura-curve construction when $[K:\Q]=2$;
see \cite[Theorem~2.1, condition (C)]{Saito}.  The corresponding
two-dimensional Galois representation occurs in the first cohomology
of the compact Shimura curve by
\cite[Theorem~(B'), \S2.2]{Carayol} and
\cite[Lemma~3.1(ii), equations (3.3), (3.8)]{Saito}.
For weight $(2,2;2)$, the algebraic representation $\rho(k)$ defining
Saito's coefficient sheaf is trivial \cite[\S3]{Saito}, so ordinary
$H^1$ is sufficient.

Saito's equation (3.8) places the Hecke-side representation
$\sigma_h(\Pi_v)$ in $H^1$, whereas the representation in (2.2)
is written with the contragredient $\check\sigma_h(\Pi_v)$.
We therefore denote the actual $H^1$ summand by $\rho^S_{\Pi,3}$
and its dual by $\rho_{\Pi,3}$.  The canonical Galois-equivariant
isomorphism
\[
 V_3(J_U)\otimes E_\lambda
  \simeq H^1_{\mathrm{\acute et}}(X_{U,\overline K},E_\lambda)^\vee
\]
then gives \eqref{eq:carayol-tate}.  The commuting Hecke and Galois
actions in \cite[(3.3), (3.8)]{Saito} show that the dual is taken on
the same eigensystem, rather than after an automorphic twist.
Chebotarev and semisimplicity identify this constituent with the
representation attached to that eigensystem.
\end{proof}

\begin{remark}
\label{rem:cohom-scope}
Write $\mathrm{Fr}^{\rm ar}_v$ for arithmetic Frobenius and
$\mathrm{Fr}^{\rm geo}_v=(\mathrm{Fr}^{\rm ar}_v)^{-1}$ for its
geometric inverse.  At an unramified prime suppose
\[
 \det(1-\mathrm{Fr}^{\rm geo}_v T\mid\rho^S_{\Pi,3})
  =1-A_vT+D_vT^2.
\]
Because $\rho_{\Pi,3}=(\rho^S_{\Pi,3})^\vee$, the eigenvalues of
$\mathrm{Fr}^{\rm ar}_v$ on $\rho_{\Pi,3}$ are the same as the
eigenvalues of $\mathrm{Fr}^{\rm geo}_v$ on $\rho^S_{\Pi,3}$.
Consequently
\begin{equation}\label{eq:cohom-frob-poly}
 \det(1-\mathrm{Fr}^{\rm ar}_vT\mid\rho_{\Pi,3})
   =1-A_vT+D_vT^2.
\end{equation}
Thus the geometric Hecke polynomial on the $H^1$ summand is compared
directly with the arithmetic Frobenius polynomial on the covariant
Tate module \emph{after} the cohomological dual has been taken.
There is no additional reciprocal-polynomial operation and no extra
finite or Tate twist.  In particular, if $\rho_{\Pi,3}\simeq r_3$,
then $A_v=a_v$ and
$D_v=(Nv)\chi^{-1}(\mathrm{Fr}^{\rm ar}_v)$ in the notation of
Lemma~\ref{lem:cohom-normalization}.
\end{remark}

\begin{lemma}
\label{lem:dyadic-conductor}
Let $M=M_{\mathcal C,S}(z)\ge1$ and let $r_\ell$ be a member of
the geometric compatible system obtained from the same
normalized quaternionic quasi-isogenies and the cochain
$\beta$, for an auxiliary coefficient prime $\ell$ distinct
from a given residue characteristic.  Suppose the finite
full-three field $L_0/K$ and the sign splitting have the bounds
of Corollary~\ref{cor:beta-quant}.  Then the Artin conductor
exponents at the fixed primes satisfy
\[
 a_{\mathcal P_2}(r_\ell)\ll_S 1+\log M,
 \qquad
 a_{\mathcal P}(r_\ell)\ll_S1
 \quad(\mathcal P\in S,\ \mathcal P\nmid2),
\]
where $\mathcal P_2$ is the unique prime above two.
At primes outside $S$ with positive odd contact the
representations are tame, with conductor exponent at most two;
outside $S$ and the odd-contact support they are unramified.
Consequently, for the modular Hilbert representation
$\Pi_z$ corresponding to the system, there are fixed $C,A_0>0$
such that
\begin{equation}\label{eq:uniform-conductor}
 N_{K/\Q}(\mathfrak n(\Pi_z))\le C M^{A_0}.
\end{equation}
\end{lemma}
\begin{proof}
The degree of $L_0/K$ is bounded independently of $z$.
At a local field over a fixed prime, finite extensions of this
bounded degree have uniformly bounded different exponents
and upper ramification breaks.  The QM surface has potentially
good reduction, as proved in Lemma~\ref{lem:nonordinary}.
At residue characteristic different from three, triviality of
$A[3]$ over $L_0$ and torsion-freeness of the principal
congruence subgroup $1+3M_4(\Z_3)$ force the finite inertia
image on $V_3(A)$ to be trivial.  The three-adic Tate module
is therefore unramified there over $L_0$.

For the unique dyadic prime use $\ell=3$.  On $G_{L_0}$ the
normalized representation is the Tate representation twisted
by $\beta^{-1}|_{G_{L_0}}$.  This is a genuine character with
square $\chi|_{G_{L_0}}$.  By the explicit conductor estimate in
Theorem~\ref{thm:cochain}, the local conductor exponent of
$\chi$ is $O(1+\log M)$.  After pullback to the finitely many
bounded-degree local fields, this same estimate holds with a
fixed multiplicative constant.  For a fixed local field $F$
over $\Q_2$, the squaring map sends sufficiently deep principal
units $U_F^n$ onto $U_F^{n+v_F(2)}$.  Thus a character whose
square is trivial on $U_F^n$ is trivial on the uniformly deeper
$U_F^{n+v_F(2)}$.  Applied to $\beta|_{G_{L_0}}$, this gives
an upper ramification cut $O(1+\log M)$.  Passing back to $K$
through the bounded-degree local extension and using the
conductor--different formula bounds the two-dimensional
conductor by $O(1+\log M)$.

At a fixed odd prime, the finite-order character $\chi$ is tame;
its square root $\beta|_{G_{L_0}}$ is also tame on wild inertia,
since squaring is an automorphism of each pro-odd wild inertia
quotient.  At the prime above three, an auxiliary prime
$\ell=5$ and the field obtained by adjoining $A[5]$ yield a
bounded-degree extension where the potentially good surface
has good reduction.  This again bounds the wild inertia
conductor uniformly.  Compatibility of the geometric system
at unramified prime ideals follows from the common quasi-isogenies;
Skinner's compatibility theorem away from the coefficient
prime and at primes above it \cite[equations (1)--(2),
Theorem~1]{Skinner} identifies the local Weil--Deligne
conductors of the geometric compatible system with those of
$\Pi_z$, including at primes above three.  This theorem
is available even before a nonzero prime ideal with a discrete-series local component has
been identified, so there is no circular use of Carayol's
cohomological construction.

At odd-contact primes outside $S$, the full-three cover is
tamely ramified with inertia order at most two, and the
corrected scalar character is tame.  Hence the Artin
conductor exponent of a two-dimensional representation is
at most two.  There is no ramification at the remaining
primes outside $S$.  The fixed odd primes contribute a
constant; the dyadic contribution is at most
$2^{C(1+\log M)}\ll M^{C'}$, and the nonfixed tame primes
contribute at most $M^2$.  This proves
\eqref{eq:uniform-conductor}.
\end{proof}

\begin{proposition}\label{prop:height}
Suppose that the following assertions hold uniformly for every
$z\in K\setminus\mathcal C$ satisfying \eqref{eq:tube}.
\begin{enumerate}
\item The integral full-three PEL cover is identified with the
actual $3$-torsion marking action at $j=0$, and its local
inertia on $r_3$ is projectively the right action of
$\gamma=1+i+j$.
\item The normalized $r_3$ is modular with a parallel-weight-two
Hilbert cuspidal representation $\Pi_z$ of
$\mathrm{GL}_2(\mathbb A_K)$.  Apply Skinner's
local--global compatibility \cite[equation (1)]{Skinner} at
$\mathfrak q_*$ to identify its Weil--Deligne parameter with
that of $r_3$, in the fixed local Langlands normalization.
\item The arithmetic-covariant Hecke normalization of $\Pi_z$
coincides with the already fixed representation $r_3$, meaning
that its polynomial at almost all finite prime ideals is the first polynomial
of \eqref{eq:matching-polynomial}, with
$\vartheta=\chi^{-1}$.
\end{enumerate}
If the conductor of $\Pi_z$ satisfies
$N_{K/\Q}(\mathfrak n(\Pi_z))\le C M_{\mathcal C,S}(z)^A$
for fixed $C,A\ge1$, and the chosen QM fibre $A_z$ is a
geometric simple abelian surface outside the CM locus,
then there is $H>0$ with
\[
 h(z)\le H M_{\mathcal C,S}(z)^{24A}.
\]
The CM fibres are covered by the bounded-degree finiteness statement of
Lemma~\ref{lem:cm-bounded-degree} and may be incorporated by adjusting $H$.
\end{proposition}
\begin{proof}
Let $M=M_{\mathcal C,S}(z)$. The local inertia at
$\mathfrak q_*$ is projectively of order four, and a lift
of Frobenius in the tame Weil group conjugates it to its
$367$th power. As $367\equiv-1\pmod4$, Frobenius exchanges
its two eigenlines. Therefore the two-dimensional Weil
representation is absolutely irreducible.  Apply
Skinner's unconditional local--global compatibility for Hilbert
forms at $v\nmid3$ \cite[Introduction, equation (1)]{Skinner}
\emph{before} invoking Carayol's compact Shimura-curve
construction.  This identifies the Weil--Deligne parameter of
$\Pi_{z,\mathfrak q_*}$ with that of the original $r_3$
(up to the explicitly fixed local Langlands normalization).
Lemma~\ref{lem:skinner-supercuspidal} now gives
$\Pi_{z,\mathfrak q_*}$ supercuspidal.
It admits a Jacquet--Langlands transfer to the fixed
quaternion algebra ramified at $\mathfrak q_*$ and one real
embedding.  Proposition~\ref{prop:carayol-tate} now constructs
a two-dimensional constituent $W_\Pi$ of the covariant Tate
module of the compact Shimura-curve Jacobian, with the
arithmetic Frobenius polynomial prescribed by the third
hypothesis. Lemma~\ref{lem:cohom-normalization} identifies
$W_\Pi\simeq r_3$; its role is uniqueness, while
Proposition~\ref{prop:carayol-tate} supplies existence.
Corollary~\ref{cor:cohom-hom} then identifies $A_z$ as a
factor over the original controlled field $L^\prime$,
without an additional Hom extension.

Choose a fixed sufficiently neat level $U_0$ in this
quaternion algebra and let $U\subset U_0$ be a level
realizing the given automorphic representation. By the
standard congruence subgroup index estimate, one can take
\[
 d_U=[U_0:U]\ll
 N_{K/\Q}(\mathfrak n(\Pi_z))^4\ll M^{4A}.
\]
The forgetful cover of smooth compact Shimura curves
$X_U\to X_{U_0}$ is finite etale on geometric connected
components. Fix a Belyi map on each component of the fixed
base curve. If the degrees of the connected components of
$X_U\to X_{U_0}$ are $d_1,\ldots,d_t$, then
$\sum_i d_i\ll d_U$, and composition with the fixed Belyi
maps gives
\[
 \deg_B(X_{U,i})\ll d_i\ \ \t{and}\ \  g(X_{U,i})\ll d_i.
\]
Javanpeykar \cite[Theorem~1.1.1]{JB} gives,
for every smooth projective connected curve of genus at least one,
\[
 h_F(\operatorname{Jac}X)
 \le 13\cdot10^6 \,g(X)\deg_B(X)^5.
\]
Thus $h_F(J_{U,i})\ll d_i^6$. The genus-zero components
have trivial Jacobian and contribute zero. Additivity of stable Faltings
height and $\sum_i d_i^6\le(\sum_i d_i)^6$ give
\[
 h_F(J_U)\ll d_U^6\ll M^{24A}.
\]
Here $h_F$ denotes stable Faltings height, and $J_U$ is the
identity component of the relative Picard variety of the
entire proper Shimura curve. It is defined over the base
field even when its geometric components are permuted.

By Proposition~\ref{prop:carayol-tate}, the matching lemma,
and Faltings' isogeny theorem \cite{Fal}, the geometrically simple fibre $A_z$
is an isogeny factor of $J_U$. The field of comparison can
be chosen with degree polynomially bounded in $M$:
Corollary~\ref{cor:beta-quant} controls the field of the
original Tate realization and the scalar cochain.
We use the exact two-step quantitative bound in
\cite[Corollary~1.5 and Theorem~1.4]{GR2014}.  This avoids
choosing an uncontrolled complement or an unbounded-degree
isogeny.  Let $B_z\subset J_{U,L'}$ be the image of the
nonzero $L'$-homomorphism $A_z\to J_U$ from
Corollary~\ref{cor:cohom-hom}.  Since $A_z$ is geometrically
simple, $B_z$ is a two-dimensional abelian subvariety defined
over $L'$, and $A_z$ and $B_z$ are $L'$-isogenous.

Put $g_U=\dim J_U\ll M^{4A}$ and $d'=[L':\Q]\ll M$.
We use the stable Faltings-height version of
R\'emond's inequality \cite[Corollary~1.2, and the observation
following Corollary~1.3]{RemondFaltings}:
for an abelian subvariety $B$ of an abelian variety $J$,
both defined over the same number field,
\begin{equation}\label{eq:remond-subvariety}
 h_F(B)+h_F(J/B)\le h_F(J).
\end{equation}
The uniform lower bound for stable Faltings heights
$h_F(C)\ge-\kappa_0\dim C$ (in the original normalization)
therefore implies
\[
 h_F(B_z)\le h_F(J_U)+\kappa_0(g_U-2)
           \ll M^{24A}.
\]
Here \eqref{eq:remond-subvariety} is applied over $L'$:
$B_z$ is an $L'$-subvariety by
Corollary~\ref{cor:cohom-hom}, and no descent of $B_z$ to
$K$ is assumed.  This avoids a large-dimension application
of the quantitative isogeny theorem to $J_U$.

For the remaining isogeny $B_z\to A_z$ the dimension is
fixed at two.  Recall the Gaudron--R\'emond quantity
\cite[p.~2058]{GR2014}:
\[
 \begin{gathered}
 \alpha(g)=2^{10}g^3=1024g^3,\\
 \kappa(C;k)=\Bigl((14g)^{64g^2}[k:\Q]
          \max\{h_F(C),\log[k:\Q],1\}^{2}\Bigr)^{\alpha(g)},
 \qquad g=\dim C.
 \end{gathered}
\]
Applying
\cite[Theorem~1.4]{GR2014} to the $L'$-isogenous
abelian surfaces $B_z$ and $A_z$, with $B_z$ as the
known-height source, gives an $L'$-isogeny of degree at
most $\kappa(B_z;L')$.  Since $\dim B_z=2$ and
$[L':\Q]\ll M$, we have
\[
 \begin{split}
 \log\kappa(B_z;L')
 &\ll 1+\log[L':\Q]+
          \log\max\{h_F(B_z),1\}
 \ll 1+\log M.
 \end{split}
\]
The isogeny inequality for stable Faltings heights then yields
\begin{equation}\label{eq:explicit-height-final}
 h_F(A_z)\le h_F(B_z)+\frac12\log\kappa(B_z;L')
        \ll M^{24A}.
\end{equation}
The implicit constants depend on the fixed base field,
Shimura datum and choice of a fixed Belyi map.

To compare the height of $z$, pass to the \emph{fixed}
compact fine PEL curve $X$ with its level-three and auxiliary
prime-to-three labels, rather than trying to express Igusa
invariants as functions on the Atkin--Lehner quotient.
Let $\lambda_X$ denote the Hodge determinant line bundle
of its universal principally polarized abelian surface.
The map $X\to\mathcal A_2^{\rm Sat}$ to the Satake
compactification is nonconstant: on the complex
uniformization the defining quaternionic upper-half-plane
period map is nonconstant.  The ample Hodge line on
$\mathcal A_2^{\rm Sat}$ therefore has positive-degree
pullback to each relevant compact component of $X$.
Consequently, the height machine on a projective curve gives
\[
  h(j(x))\ll 1+h_{\lambda_X}(x),
\]
where $j$ is the fixed coarse-quotient function and $x$ is
a fine lift representing the polarized fibre underlying
$A_z$.  Changing this lift by one of the finitely many
Atkin--Lehner transitions changes the underlying abelian
surface by an isogeny of degree bounded solely in terms
of the fixed discriminant-six datum.  The isogeny inequality
therefore changes $h_F$ by at most a fixed constant.
The Hodge line is metrized by the Faltings metric.  To make the
comparison uniform, fix once and for all an even theta level on the
fine PEL curve.  The associated theta embedding gives a projective
height $h_\Theta$ whose pullback to $X$ is a Weil height for a positive
multiple of $\lambda_X$, up to $O(1)$.  Pazuki's
\cite[Theorem~1.1, Remark~1.2 and Corollary~1.3]{Pazuki} comparison,
with the dimension and theta level fixed, gives
\[
 \left|h_\Theta-\frac12h_F(A_x)\right|
 \le C\log(\max\{h_\Theta,1\}+2)
\]
for a constant $C$ depending only on the fixed level and dimension.
The logarithmic term is absorbed for large $h_\Theta$ (and the bounded
range is absorbed into the constant), so
\[
  h_{\lambda_X}(x)\ll1+\max\{0,h_F(A_x)\}.
\]
All heights here are absolute normalized heights, hence the constants
are independent of the field of definition of $x$.  This argument uses
the polarized fine cover, where the moduli map exists; it does
\emph{not} assert that Igusa invariants descend as rational functions
of the unpolarized Atkin--Lehner coordinate $j$.
Finally $j(x(z))$ is a nonconstant rational function of $z$
over the fixed number field, so the elementary height formula
for a rational map of $\mathbb P^1$ implies
$h(z)\ll1+h(j(x(z)))$.  Combining the three inequalities
with \eqref{eq:explicit-height-final} yields the claimed
$h(z)\ll M^{24A}$.
\end{proof}

\begin{remark}
\label{rem:exactmatch}
The matching polynomial of Lemma~\ref{lem:cohom-normalization}
is imposed for the original $r_3$ with determinant
$\chi_3\chi^{-1}$.  The source of the polynomial is
Proposition~\ref{prop:carayol-tate}, not the matching lemma
itself.  Saito's (2.2) uses geometric Frobenius and the local
representation $\check\sigma_h(\Pi_v)$, whereas the $H^1$
constituent in (3.8) uses $\sigma_h(\Pi_v)$.  The latter is dualized into the \emph{covariant Tate module};
formula \eqref{eq:cohom-frob-poly} shows that its geometric
characteristic polynomial is then exactly the arithmetic Tate
polynomial in \eqref{eq:matching-polynomial}.  The finite character $\chi$
was defined by the sign-cocycle construction before any
modularity assertion.
\end{remark}
\section{The residual dihedral image and modularity}\label{sec:residual}
Fix a parameter $z\in K\setminus\mathcal C$ satisfying
\eqref{eq:tube} and $\sigma(P_{\mathcal C}(z))>0$ at both real
embeddings. We first treat the non-CM fibres, whose geometric
endomorphism algebra is the quaternion division algebra $B$.
The CM fibres are treated separately in Theorem~\ref{thm:height-final}.
The quaternion order is $O=\Z[i,(1+j)/2]$, and $r_3$
is the same normalized two-dimensional representation defined in
Section~\ref{sec:cochain}. In this section $D_8$ denotes the
dihedral group of order eight.

\begin{lemma}\label{lem:Gamma}
Write $\mathbb A_f$ for the finite adele ring of $\Q$.
Let $K_0=\widehat O^{\times}$ and
$K_{\mathrm{lev},3}=\ker(K_0\to(O/3O)^\times)$. If $N$ denotes the
finite-adelic normalizer of $\widehat O$, the effective
normalizer group
$
\Gamma=N/(\mathbb A_f^\times K_{\mathrm{lev},3})$
acts on the generic full-three-marked PEL fibre. It has order
$144$ and a normal subgroup $U\simeq C_3^2$ with
\begin{equation}\label{eq:Gamma}
\Gamma/U\simeq D_8\times C_2\ \ \t{and}\ \ 
(\ker\delta_3)/U\simeq C_4\times C_2,
\end{equation}
where $\delta_3$ is the Atkin--Lehner valuation parity at three.
The full generic unbranched fibre is a torsor for $\Gamma$.
\end{lemma}
\begin{proof}
For a non-CM geometric QM surface the centralizer of the left
$O$-action in $\operatorname{End}^0(A)=B$ is $\Q$; hence its
polarized $O$-linear automorphism group is precisely $\{\pm1\}$.
The fine full-three functor exists over $\Z[1/30]$ as the
quotient of the auxiliary \emph{full-five}-and-three marking
scheme by the free full-five relabeling group.  Using full-five
level makes the relabeling group action genuinely defined;
$V_1(5)$ alone is not a normal level.  The underlying prime-to-six
fine integral construction is
\cite[Theorem~2.1 and Corollary~2.3]{Buzzard}, while
Lemma~\ref{lem:full-three-pel} proves that the added finite
\'{e}tale sheaf
$\operatorname{Isom}_O(O/3O,A[3])$ is exactly the complete PEL
similitude level at three.  In particular, there is no omitted
pairing or determinant component.  Its multiplier is the reduced
norm modulo three, and $\operatorname{Nrd}(\gamma)=2$ shows that
both multiplier components occur.  Thus the generic forgetful
fibre over a single unmarked $O$-surface consists of
\[
 |(O/3O)^\times|/|\{\pm1\}|=72/2=36
\]
actual marked points.

The polarized forgetful map for discriminant six is the quotient
by the entire four-element positive Atkin--Lehner group
$W=\langle w_2,w_3\rangle$, not by a proper subgroup
\cite[Theorem 2.1]{BabaGranath}.  The four unmarked points
above a generic $j$ are distinct; the marked generic fibre
therefore contains exactly $4\cdot36=144$ points.  This
establishes its size directly from the moduli functor, without
assuming the arithmetic torsor property.

At split finite primes the normalizer of a maximal order is the
centre times its units. At the division primes two and three,
the quotient by centre and units is the valuation-parity group
$C_2$.  Hence $N/(\mathbb A_f^\times K_0)\simeq C_2^2$.
At three, with $\varepsilon$ the class of a local division
uniformizer, the residue ring is
\[
 O/3O\simeq\F_9\oplus\F_9\varepsilon\ \ \t{with}\ 
 \varepsilon^2=0\ \t{and}\ 
 \varepsilon a=a^3\varepsilon.
\]
Its units have order $8\cdot9=72$ and its radical-unit subgroup
$U=1+\F_9\varepsilon$ is normal, of order nine.  Modulo the
central $\{\pm1\}$ the level kernel has order 36.  Central
finite ideles act trivially: they differ from rational scalars
by integral units, whose residue at three is $\pm1$, already
trivial in the effective marking group.  This proves
$|\Gamma|=4\cdot36=144$.

Modulo $U$ let $t,s,r$ be represented by $1+i,i,j$, respectively.
Their norms have independent valuation parities at two and
three for $s,r$, while $t$ is a unit of norm $-1$ and has
residue order four modulo $\F_3^\times$.  Quaternion
conjugation gives
\[
 t^4=s^2=r^2=1,\qquad st=ts,\quad sr=rs,\quad
 rtr^{-1}=t^{-1}.
\]
These generators give all 16 elements and prove
$\Gamma/U\simeq D_8\times C_2$; the kernel of the three-adic
valuation parity is generated modulo $U$ by $t,s$, giving
$C_4\times C_2$.

The level-unit subgroup acts by literal relabeling of the
universal three-torsion.  For the nonunit Atkin--Lehner
component at three, that description would be incorrect:
the action is the normalizer Hecke action on the canonical
Shimura tower, with the fractional-lattice transport of
Lemma~\ref{lem:PEL-lattice-transport}.  The assertion that $K_{\mathrm{lev},3}$ is normalized can be checked
locally.  Outside three it is the maximal-order unit group;
at three it is $1+3O_3$, and the division maximal order
$O_3$ is unique, so every $n_3\in B_3^\times$ preserves
$3O_3$ by conjugation.  Thus each adelic normalizer
conjugates $K_{\mathrm{lev},3}$ into itself, and its Hecke correspondence
is a one-valued automorphism of the canonical level-$K_{\mathrm{lev},3}$
curve.  The normalizer correspondence is defined over the reflex
field $\Q$ by the adelic Shimura tower: its generic-fibre
moduli interpretation is the right action on the Tate
trivialization from \cite[\S\S4.11--4.15]{Deligne} and
\cite[Theorem~4.2]{Molina}, while the fixed-polarization
presentation is \cite[\S3.1]{Rotger}.  In the present
discriminant-six case the underlying canonical quotient
and its Atkin--Lehner operators are explicitly defined over
$\Q$ by \cite[Theorem~3.10 and Corollary~3.11]{BabaGranath}.
Consequently, $\Gamma$ acts on the actual
$\Q$-defined generic marked fibre.  Its level subgroup acts
transitively on the 36 markings of any unmarked point, and its
image $W$ acts transitively on the four unmarked points.
Thus $\Gamma$ acts transitively on all 144 points.  Since it
also has order 144, the action is free and transitive: the
actual generic fibre is a $\Gamma$-torsor. The action is
$\Q$-defined, so its commuting arithmetic Galois action
embeds in the right regular centralizer $\Gamma^{\rm op}$.
\end{proof}

\begin{lemma}
\label{lem:balanced-morita}
Put $D=B\otimes_{\Q}\Q_3$ and let $\mathcal O_D$ be its
maximal order, with maximal two-sided ideal $\mathfrak P$.
There is a finite coefficient extension
$F=\Q_3(i,\pi)$ (with $i^2=2$ and $\pi^2=-3$) with residue field
$\F_9$, and a splitting $D\otimes_{\Q_3}F\simeq M_2(F)$
such that the projective normalizer $D^\times/\Q_3^\times$
preserves the homothety class of the standard rank-two
$\mathcal O_F$-lattice.  In this normalization the reduction
of the projective normalizer has image $D_8$ and kernel
$\Q_3^\times(1+\mathfrak P)$.
\end{lemma}
\begin{proof}
The extension $E=\Q_3(i)$ is unramified quadratic, since
$X^2-2$ is irreducible modulo three, and $j$ is a division
uniformizer satisfying $j^2=-3$ and $ja=\bar a j$
for $a\in E$.  Consequently
$\mathcal O_D=\mathbb Z_3[i]\oplus\mathbb Z_3[i]j$.
On $F^2$ take the matrices
\[
 i\longmapsto
 \begin{pmatrix}i&0\\0&-i\end{pmatrix},\qquad
 j\longmapsto
 \begin{pmatrix}0&-3\\1&0\end{pmatrix},
 \qquad H=\operatorname{diag}(\pi,1).
\]
Conjugating by $H$ gives, for $a,b\in\mathbb Z_3[i]$, the formulas
\begin{equation}\label{eq:balanced-three-integral}
 H^{-1}(a+bj)H
 =\begin{pmatrix}a&\pi b\\\pi\bar b&\bar a\end{pmatrix}
 \ \ \t{and}\ \ 
 H^{-1}(j/\pi)H=\begin{pmatrix}0&1\\1&0\end{pmatrix}.
\end{equation}
All matrices in the first display with $a+bj\in\mathcal
O_D^\times$ are integral, with unit determinant; reduction
modulo $\pi$ is
\[
 \begin{pmatrix}x&0\\0&x^3\end{pmatrix}\ \ \t{with}\ \ 
  x\in\F_9^\times.
\]
The second matrix is integral and invertible.  As
$D^\times=\Q_3^\times\mathcal O_D^\times\sqcup
 \Q_3^\times\mathcal O_D^\times j$,
these matrices represent the whole projective local
normalizer on one balanced lattice.
After removing scalar matrices, the diagonal group is
$\F_9^\times/\F_3^\times\simeq C_4$ and the swap matrix
acts on it by inversion.  The kernel in $\mathcal
O_D^\times$ consists of units whose residue in $\F_9$
belongs to $\F_3^\times$, namely
$\mathbb Z_3^\times(1+\mathfrak P)$.
Adding central powers of three proves the claimed kernel.
In particular, the radical units $1+\mathfrak P$ reduce to
the identity on the balanced lattice, even though they need
not reduce trivially on a different lattice in the same
split rational representation.
\end{proof}

\begin{remark}
\label{rem:balanced-morita-r3}
For the continuous two-dimensional representation $r_3$
constructed in Remark~\ref{rem:hom}, the global canonical
PEL comparison of Lemma~\ref{lem:full-three-pel} and
Proposition~\ref{prop:pel-inertia}, together with the
QM-isogeny-class comparison of \cite[Theorem~4.2]{Molina},
identifies $\mathbb P r_3(G_K)$ with a subgroup of the
above projective normalizer.  After a finite extension of
the coefficient field, $r_3$ preserves the balanced lattice;
its residual projective image is the reduction of the
\emph{actual Galois subgroup} inside $D_8$.

Fix a geometric full-three
PEL point $y$ above the $K$-rational generic parameter.
The canonical $\Q$-defined $\Gamma$-torsor structure of
Lemma~\ref{lem:Gamma} gives, for each $g\in G_K$, a
unique $\delta_g\in\Gamma$ with ${}^g y=y\cdot\delta_g$.
Choose its adelic normalizer representative $n_g$.
Because the normalizer acts on the \emph{canonical}
marked moduli functor, this equality compares genuine
$O$-linear three-torsion frames and the transported
polarization.  It determines the local Tate action only
\emph{modulo the principal three-level kernel}, not as an
exact equality of $3$-adic matrices.  More precisely,
Rotger's polarized transport in
Lemma~\ref{lem:polarized-normalizer-comparison} and Molina's
Theorem~4.2, including the explicit torsion identity
(4.1) in its proof, give a $k_g\in1+3\mathcal O_D$
(depending on $g$) with
\[
 [r_3(g)]=[R_{n_{g,3}^{-1}k_g}]
 \qquad\text{in }\operatorname{PGL}_2(\overline\Q_3),
\]
up to the consistent opposite-action convention and
conjugation by a fixed splitting basis.  For a nonunit
$n_{g,3}$ the right side uses the fractional Tate-lattice
transport of Lemma~\ref{lem:PEL-lattice-transport}.
The ambiguity from the chosen $B$-linear quasi-isogeny is
central because $\operatorname{End}_B^0(A)=\Q$.
Crucially, $1+3\mathcal O_D\subset1+\mathfrak P$, and
the preceding matrix calculation shows directly that
$1+\mathfrak P$ reduces to the identity on the
balanced lattice.  Thus the principal-level ambiguity
vanishes in the \emph{residual projectivization}.
This proves the global projective inclusion in the full
local normalizer and identifies its finite residual quotient.
It does not assert that $r_3(G_K)$ itself is finite: over
the full-three splitting field it may retain an infinite
principal-congruence image.

Choose a finite coefficient field $F'$ containing both
$F$ and the splitting field of the original $r_3$.
Each $r_3(g)$ therefore differs by an $F'^\times$-scalar
from a matrix in $\mathrm{GL}_2(\mathcal O_{F'})$
representing the right normalizer (up to opposite action,
which replaces the matrix by its inverse).  The continuous determinant
$\det r_3:G_K\to F'^\times$ has compact image, so its
valuation is identically zero.  The normalizer matrix also
has unit determinant.  Thus the scalar has valuation zero,
so $r_3(g)\in\mathrm{GL}_2(\mathcal O_{F'})$ for every $g$.
This proves that the same lattice is $G_K$-stable, and its
reduction is exactly the explicit projective reduction above.
\end{remark}

\begin{lemma}\label{lem:projective-normalizer}
Let $D=B\otimes_\Q\Q_3$, let $\mathcal O_D$ be its maximal order and
let $\mathfrak P$ be the maximal two-sided ideal of $\mathcal O_D$.
There is a canonical quotient
\begin{equation}\label{eq:local-D8-quotient}
 D^\times/\bigl(\Q_3^\times(1+\mathfrak P)\bigr)
 \simeq
 (\F_9^\times/\F_3^\times)\rtimes
 \Gal(\F_9/\F_3)\simeq D_8.
\end{equation}
Under the notation $t,s,r$ of Lemma~\ref{lem:Gamma}, the images of
$t=[1+i]$, $s=[i]$ and $r=[j]$ are respectively a rotation of order
four, its square, and a reflection.  In particular the kernel of the
map
\[
 \Gamma/U\longrightarrow D_8
\]
is the order-two subgroup generated by $st^{-2}$; it is contained in
$\ker\delta_3$.  Hence $\delta_3$ descends to the quotient $D_8$, where
its kernel is exactly the rotation subgroup $C_4$.
\end{lemma}
\begin{proof}
The residue division field of $\mathcal O_D$ is $\F_9$.
Reduction gives
$\mathcal O_D^\times/(1+\mathfrak P)\simeq\F_9^\times$.
Modulo central units this becomes
$\F_9^\times/\F_3^\times\simeq C_4$.  A division
uniformizer has odd $D$-valuation; modulo the centre its valuation
has order two, and conjugation by it induces the nontrivial
Frobenius $a\mapsto a^3$ of $\F_9/\F_3$.  On
$\F_9^\times/\F_3^\times$ this is inversion, proving
\eqref{eq:local-D8-quotient}.

Modulo $\mathfrak P$, $i^2=2=-1$ and
$(1+i)^2=2i$, so projectively $[i]=[1+i]^2$.  Thus $t$ maps to
a generator of the rotation group and $s$ to its square.  The class
of $j$ has odd division valuation and maps to the Frobenius coset,
hence to a reflection.  It follows from the presentation in
Lemma~\ref{lem:Gamma} that the only additional element of
$\Gamma/U$ killed by this projective quotient is $st^{-2}$.
Both $s$ and $t$ have even three-adic valuation parity, so
$st^{-2}\in\ker\delta_3$.  The last assertion follows.

Finally, Remark~\ref{rem:balanced-morita-r3} supplies a single
stable rank-two integral coefficient lattice on which this quotient
is the genuine projective reduction; Lemma~\ref{lem:PEL-lattice-transport}
identifies the corresponding integral Tate transports.  Applying
\cite[Theorem~4.2]{Molina} to the QM $K$-surface produced
from the rational Atkin--Lehner quotient, as in
Proposition~\ref{prop:pel-inertia}, identifies it with the
actual residual projectivization of the same $r_3$.
The coefficient-field extension only conjugates the
representation; it does not change the residual projective
image or its absolute irreducibility.
\end{proof}

\begin{lemma}\label{lem:global-D8}
For the residual semisimplification of the original $r_3$, we have
\begin{equation}\label{eq:global-D8}
\mathbb P\bar r_3(G_K)\simeq D_8.
\end{equation}
Its rotation quadratic character is the restriction of
$\delta_3$ to $G_K$.
\end{lemma}
\begin{proof}
Let $L_0/K$ be the Galois splitting field of a generic fibre
of the full-three cover. The Galois permutations commute with the
$\Q$-defined free transitive normalizer action, hence
$\Gal(L_0/K)$ embeds in $\Gamma^{\rm op}$.
Over $G_{L_0}$ the three-torsion marking is constant, and
$R_3(g)=\beta(g)^{-1}V_3(A)(g)$.
As $R_3\simeq r_3^{\oplus2}$ over a splitting coefficient field,
all eigenvalues of the residual matrix $\bar r_3(g)$ coincide.
Its projective image is unipotent and therefore a finite
three-group. It is normal in the total projective image.
A nontrivial finite three-subgroup of
$\operatorname{PGL}_2(\overline\F_3)$ has a unique fixed point
in $\mathbb P^1$. Normality would make that point fixed by the
whole image, contradicting the residual absolute irreducibility
from Proposition~\ref{prop:pel-inertia}. Thus
$\mathbb P\bar r_3(G_{L_0})=1$.
The intersection of $\Gal(L_0/K)$ with the normal three-group
$U$ in Lemma~\ref{lem:Gamma} similarly has trivial projective
image. Hence the total projective image is a subquotient of
$D_8\times C_2$, a group of exponent at most four and order
at most sixteen. Finite two-subgroups of
$\operatorname{PGL}_2(\overline\F_3)$ are cyclic or dihedral.
A cyclic or dihedral group of order sixteen has an element of
order eight. Hence the image has order at most eight.
At $\mathfrak q_*$, fourfold inertia and the Frobenius inversion
proved in Proposition~\ref{prop:pel-inertia} generate a dihedral
group of order eight. Thus \eqref{eq:global-D8} follows.
The element $\gamma=1+i+j$ satisfies
$[\gamma]=[t]^3[s]$ in $\Gamma/U$, so
$\delta_3(\gamma)=1$.  By Lemma~\ref{lem:projective-normalizer},
$\delta_3$ descends through the actual projective normalizer quotient
and its kernel there is exactly the rotation subgroup $C_4$.
Since the projective image has just been proved to be $D_8$, this
identifies the rotation quadratic character with the restriction of
$\delta_3$ to $G_K$.
\end{proof}

\begin{lemma}\label{lem:rotation-field}
The quadratic field corresponding to the rotation subgroup
of \eqref{eq:global-D8} is
\begin{equation}\label{eq:Mz}
M_z=K\bigl(\sqrt{-P_{\mathcal C}(z)}\bigr).
\end{equation}
It is a CM extension of $K$.
\end{lemma}
\begin{proof}
The kernel of the three-adic parity character maps under
forgetting the full-three marking to $\langle w_2\rangle$.
In view of Baba--Granath \cite[Theorem 3.10, Corollary 3.11 and Theorem 3.17]{BabaGranath}, the corresponding quadratic intermediate
cover is
$
\Q(j)\bigl(\sqrt{-27j-16}\bigr).
$

The square-class identity \eqref{eq:squareclass} identifies its
pullback with \eqref{eq:Mz}. Because
$\sigma(P_{\mathcal C}(z))>0$ at both real embeddings,
$M_z/K$ is totally imaginary quadratic.
\end{proof}

\begin{lemma}\label{lem:nonordinary}
The representation $r_3$ is potentially Barsotti--Tate but
not potentially ordinary at the unique prime of $K$ above three.
It is odd at each real embedding of $K$.
\end{lemma}
\begin{proof}
After a finite extension the QM abelian surface acquires
semistable reduction. Its torus part cannot admit the action
of the quaternion division algebra $B$, since its dimension
is at most two; hence the surface has potentially good reduction.
The rational Tate module of the slope-zero part of its height-four
$3$-divisible group is stable under $B\otimes\Q_3$.  Any nonzero
module over this quaternion division algebra has $\Q_3$-dimension at
least four.  Hence a nonzero slope-zero part would already have
multiplicity four and exhaust the whole height.  Polarization symmetry
would then force a slope-one part of the same positive multiplicity,
which is impossible.  Thus no Newton slope is zero, and by the same
symmetry there is no slope one.  By the Dieudonn'e--Manin slope-denominator
condition, a rational slope $a/b$ in lowest terms has multiplicity
divisible by $b$.  A slope strictly between zero and one other
than $1/2$ has denominator $b\ge3$; its distinct complementary
slope $1-a/b$ also has multiplicity at least $b$, exceeding the
total height four.  The only remaining possibility is four slopes
all equal to $1/2$.
This excludes potential ordinarity of the two-dimensional factor.
The factor is potentially crystalline with Hodge--Tate weights
$\{0,1\}$ at each embedding. At each real embedding, complex
conjugation acts outside the rotation subgroup because
$M_z/K$ is CM. Its residual representation has trace zero;
as $r_3(c)^2=1$ in characteristic zero, its trace is also zero
and its determinant is $-1$.
\end{proof}

\begin{theorem}
\label{thm:r3-modular}
For each non-CM $z$ in the positive tube, the normalized $r_3$ is
associated with a cuspidal parallel-weight-two Hilbert eigenform
over $K$.
\end{theorem}
\begin{proof}
By Lemmas~\ref{lem:global-D8} and \ref{lem:rotation-field},
the irreducible residual representation is induced from a
finite-order character of the CM quadratic field $M_z$.
Choose, after allowing a sufficiently divisible finite conductor,
an algebraic Hecke character of $M_z$ of infinity type $(1,0)$ at
every CM pair; the standard unit-compatibility obstruction is killed by
this fixed finite conductor.  The ratio of its residual
three-adic Galois character to the inducing character has a
finite-order Teichmuller lift, corresponding by class field theory
to a finite-order Hecke character. After twisting, automorphic
induction yields a cuspidal weight-two Hilbert representation
with the prescribed residual representation
\cite{ArthurClozel}. Cuspidality follows from residual
absolute irreducibility. The original $r_3$ is potentially
Barsotti--Tate, nonordinary, odd, unramified outside finitely
many prime ideals, and has determinant $\chi_3\chi^{-1}$ with
$\chi$ of finite order. Proposition~\ref{prop:pel-inertia}
provides absolute irreducibility on $G_{K(\zeta_3)}$.
We apply Gee's corrected modularity-lifting theorem
\cite[Theorem~1.1, conditions~(1)--(3)]{GeeErratum},
with $p=3$, $F=K$ and $\rho=r_3$ (after an enlargement of
its coefficient field and a stable lattice).  For precision we
check the hypotheses in exactly the order of the published
erratum.
\begin{enumerate}
\item The representation $r_3$ is continuous, unramified
outside finitely many prime ideals and has determinant
$\chi_3\chi^{-1}$; it is potentially Barsotti--Tate at
the unique prime ideal $\mathcal P_3$ of $K$ above three,
by Lemma~\ref{lem:nonordinary}.
\item The preceding automorphic induction gives a
parallel-weight-two Hilbert modular form $f$ with
$\bar\rho_f\simeq\bar r_3$.  The corrected hypothesis
requires $\rho_f$ to be potentially ordinary at
$\mathcal P_3$ \emph{only if} $r_3$ is potentially ordinary
there.  Lemma~\ref{lem:nonordinary} proves that $r_3$ is
not potentially ordinary, so this implication is vacuous.
The erratum does not impose an additional equality of
local deformation components in this situation.
\item Proposition~\ref{prop:pel-inertia} proves
$\bar r_3|_{G_{K(\zeta_3)}}$ absolutely irreducible.
The exceptional extra condition in the erratum concerns only
$p=5$ and is irrelevant here.
\end{enumerate}
Thus the three hypotheses of the \emph{corrected}
Gee theorem hold, and its conclusion proves the modularity
of the original $r_3$.
\end{proof}

\begin{lemma}\label{lem:cm-bounded-degree}
On the fixed discriminant-six Shimura curve with the fixed finite
PEL level used above, there are only finitely many CM points of
bounded degree over $\Q$.  Consequently, for each fixed quadratic
base field used in this paper, the CM fibres that occur in a
bounded-degree lift of a $K$-rational coarse parameter form a finite
set and have bounded height.
\end{lemma}
\begin{proof}
A CM point on this Shimura curve is a Heegner point attached to an
optimal embedding of an imaginary quadratic order.  Shimura
reciprocity identifies its Galois orbit, up to the uniformly bounded
stabilizer contributed by the fixed level and automorphisms, with the
corresponding ring class group; see the discussion of CM points and
Artin reciprocity in \cite[\S7, especially Proposition~7.1]{Molina}.
Thus bounded degree bounds the ring class number.  The classical
ring-class number formula, together with the Heilbronn finiteness
theorem for imaginary quadratic class numbers, implies that only
finitely many imaginary quadratic orders have class number below a
fixed bound.  For each order there are only finitely many optimal
embedding classes at the fixed level.  This proves the assertion.
\end{proof}

\begin{theorem}[Restricted seven-contact height theorem]
\label{thm:height-final}
There are positive integers $c$ and $H$ for which the estimate
\eqref{eq:height-interface} holds.
\end{theorem}
\begin{proof}
Proposition~\ref{prop:pel-inertia} gives the local order-four
projective inertia.  The relation
$\mathrm{Fr}\,\tau\,\mathrm{Fr}^{-1}=\tau^{367}$ and
$367\equiv-1\pmod4$ make the associated two-dimensional
Weil--Deligne parameter irreducible.  The local--global
compatibility at $\mathfrak q_*$ comes from
\cite[Introduction, equation (1)]{Skinner}, whose statement
is independent of any finite discrete-series hypothesis.
By Lemma~\ref{lem:skinner-supercuspidal}, the Hilbert form
corresponding to Theorem~\ref{thm:r3-modular} is
supercuspidal at $\mathfrak q_*$.  This establishes the
needed discrete-series prime ideal \emph{before} Carayol--Saito
is invoked.
It has a Jacquet--Langlands transfer to the
quaternion algebra over $K$ ramified at $\mathfrak q_*$ and one
real embedding \cite{JacquetLanglands}. For the corresponding compact Shimura curve,
Proposition~\ref{prop:carayol-tate} supplies the actual
$H^1$ constituent and its dual in the covariant Tate module.
The duality identity \eqref{eq:cohom-frob-poly} handles the
geometric-versus-arithmetic Frobenius convention.  The
matching Lemma~\ref{lem:cohom-normalization} then identifies
the Tate constituent with the original $r_3$; no additional
finite character or Tate twist enters \cite{Carayol,Saito}.
Lemma~\ref{lem:dyadic-conductor} yields fixed $C,A>0$ with
\[
 N_{K/\Q}(\mathfrak n(\Pi_z))
 \le C M_{\mathcal C,S}(z)^A.
\]
In particular, the dyadic conductor exponent is allowed to
grow as $O(1+\log M)$; it was incorrect to claim a uniform
constant exponent at that fixed prime.  The matching
constituent $W_\Pi$ is supplied by
Proposition~\ref{prop:carayol-tate}, and the comparison with
the original $r_3$ is Lemma~\ref{lem:cohom-normalization}.
All hypotheses of Proposition~\ref{prop:height} now hold.
For geometrically simple fibres it gives
$h(z)\ll M_{\mathcal C,S}(z)^{24A}$.
The CM fibres belong to a finite set by
Lemma~\ref{lem:cm-bounded-degree}; their heights are absorbed into the
constant.
This proves the stated inequality with $c=24A$.
\end{proof}

\section{Shared local certificates}
\label{sec:shared-certificates}
The norm--trace selector below is the relative $L/K$ form
of the reciprocal-trace selection method in
\cite[Section~14]{Sun38}.  We include its local proof to
make the valuation normalizations explicit.
\begin{lemma}
\label{lem:shared-normtrace}
Let $L/K$ be the fixed quadratic extension and let
$\mathcal P$ be a nonramified finite prime of $K$ of odd
residue characteristic. Let $R\in K^\times$ have odd
valuation $\alpha=v_{\mathcal P}(R)$, and let
$D_1,\ldots,D_J\in L$, $k_1,\ldots,k_J\ge1$.
On the nonzero locus, put
\[
 W=\prod_{i=1}^J\left(1+R\left(\frac{D_i}{R^{k_i}}\right)^2\right),
 \quad \Theta=\operatorname{Tr}_{L/K}(W^{-1}),
 \quad B=R/\Theta^2.
\]
Assume $W\Theta\ne0$ and that a local pole predicate on $B$
excludes negative odd valuations at $\mathcal P$.
Then $\alpha>0$, and there is at least one
$\mathcal Q\mid\mathcal P$ of $L$ for which
\[
 v_{\mathcal Q}(D_i)\ge k_i\qquad(1\le i\le J).
\]
Here valuations over $\mathcal P$ agree with the integer
normalization on $K$ because the extension is unramified.
\end{lemma}
\begin{proof}
Write $\lambda_{\mathcal Q}=v_{\mathcal Q}(W)$.
For every nonzero $D_i$, the nonconstant summand in the
corresponding factor has valuation
$\alpha+2v_{\mathcal Q}(D_i)-2k_i\alpha$, which is odd.
Thus it cannot cancel with $1$ and
\[
 \lambda_{\mathcal Q}
 =\sum_{D_i\ne0}\min\{0,\alpha+2v_{\mathcal Q}(D_i)
                           -2k_i\alpha\}\le0.
\]
Let $\lambda=\max_{\mathcal Q\mid\mathcal P}
\lambda_{\mathcal Q}$.  The integral local trace gives
$v_{\mathcal P}(\Theta)\ge-\lambda$.
As $v_{\mathcal P}(B)=\alpha-2v_{\mathcal P}(\Theta)$
is odd, the pole exclusion forces it to be positive. Hence
$\alpha\ge1$ and $\lambda\ge(1-\alpha)/2$.
Choose $\mathcal Q$ attaining $\lambda$.  Each negative
summand in its displayed sum is at least $\lambda$, whence
\[
 2v_{\mathcal Q}(D_i)
 \ge(2k_i-1)\alpha+\frac{1-\alpha}{2}
 \ge2k_i-1.
\]
Integral valuations imply $v_{\mathcal Q}(D_i)\ge k_i$;
zero inputs satisfy this automatically.
\end{proof}

The two-anchor comparison combines the formal-logarithm estimate
\cite[Lemma~3.2]{Sun38} with the subtraction argument
in \cite[Section~6]{Sun38}.
\begin{lemma}
\label{lem:shared-two-anchor}
Fix $m\ge1$ and an integer $c_*\ne0$. At a good unramified
prime of residue characteristic $p>m+1$ with $p\nmid c_*$,
let $P_i=[n_i]P$, $H=[h]P$, and $U_i=[u_i]P$. Set
$V_0=H-[c_*]P$ and $V_i=U_i-[c_*]P_i$.
Suppose the points on the primary and conjugate curves are
in the respective formal groups, and the two truncated-logarithm
ratios $L_i,L'_i$ are defined and obey
$v_{\mathcal Q}(L_i-L'_i)\ge m$.
Then
\[
 v_{\mathcal Q}(L_i-n_i)\ge m
 \qquad\text{for every }i.
\]
Consequently, any integer polynomial $f$ satisfying
$v_{\mathcal Q}(f(L_1,\ldots,L_r))\ge m$ also satisfies
$p^m\mid f(n_1,\ldots,n_r)$.
\end{lemma}
\begin{proof}
The exact formal logarithm gives the two index ratios
$u_i/h$ and $(u_i-c_*n_i)/(h-c_*)$.
As $p\nmid c_*$, at least one of $h,h-c_*$ is a unit.
For that ratio, the standard truncated-logarithm isometry
(\cite[Lemma~3.2]{Sun38}) gives the index quotient modulo
$\mathcal Q^m$. The difference condition makes the other
truncated ratio integral; the isometry then makes its index
quotient integral and gives the same congruence. Therefore
\[
 c_* L_i\equiv hL_i-(h-c_*)L'_i
 \equiv u_i-(u_i-c_*n_i)=c_*n_i\pmod{\mathcal Q^m}.
\]
Divide by the unit $c_*$. Polynomial substitution preserves
these congruences; a rational integer has $\mathcal Q$-valuation
$v_p$ at the selected unramified prime.
\end{proof}

\section{Aggregate trace heights and finite algebraic choices}
\subsection{Aggregate trace height for integer indices}
For a nonempty subset $I\subseteq\{1,\dots,10\}$ and integers
$n_i\ne0$ ($i\in I$), define
\[
 T_I=\frac12\Tr_{F/\Q}\biggl(\sum_{i\in I}x([n_i]P)\biggr)\in\Q.
\]
The nonzero points $[n_i]P$ are doubles of points
$Q_i\in E(L)$ by \eqref{eq:carrier}.  The duplication formula
for $y^2=x^3-a_E^2x$ gives, whenever $[2]Q_i\ne O$,
\[
 x([2]Q_i)=\left(\frac{x(Q_i)^2+a_E^2}{2y(Q_i)}\right)^2
 \in L^{\times2}.
\]
The left side belongs to $F$; since $L$ is totally real, this
identity makes it strictly positive at both real embeddings
of $F$.  The assertion that it must be a square \emph{in $F$}
would require an extra argument and is not used.  Hence
$T_I>0$.
\begin{lemma}\label{lem:aggregate-trace}
There are constants $C_1>0$ and $C_2\ge0$, independent of $I$ and of
the indices, such that
\begin{equation}\label{eq:aggregate-trace}
 \max_{i\in I}|n_i|^2\le C_1h(T_I)+C_2.
\end{equation}
\end{lemma}
\begin{proof}
We use the aggregate-trace proof of
\cite[Lemma~3.1]{Sun38}, with the present curve.
Over an algebraic extension, $E$ and its conjugate $E^\sigma$
become $E_0:y^2=x^3-x$ after adjoining the square roots of
$a_E=4+\sqrt3$ and $a_E^\sigma=4-\sqrt3$. The classes of these
elements are independent in $F(i)^\times/F(i)^{\times2}$:
they are totally positive, neither is an $F$-square, and their
product $13$ is not an $F$-square. A totally positive
$F$-element which becomes a square in $F(i)$ was already a
square in $F$.  The Galois involution changing the sign of just
one square root proves that the two resulting nontorsion points
on $E_0$ are independent over $\Z[i]$, modulo torsion.
As in \cite[proof of Lemma~3.1]{Sun38}, the sum of the abscissas
on $E_0^{2|I|}$ has a nonconstant restriction to every positive-
dimensional abelian translate through the selected multiple
point.  Within one block the polar leading terms at a common
pole have positive square ratios and cannot cancel; between
blocks, a common pole contradicts the independence just shown.
Von Buhren's height theorem \cite[Theorem~1.1]{VB} therefore gives
$\sum_{i\in I}n_i^2\widehat h_E(P)\le C_Ih(T_I)+C'_I$.
Since $\widehat h_E(P)>0$ and there are only finitely many
nonempty sets $I$, the constants may be chosen uniformly.
\end{proof}

The same aggregate-trace lemma holds for the carrier $E_3,P_3$ below, by the conjugate-squareclass independence of Lemma~\ref{lem:sqrt3-carrier}. In that case, use $F_3$ in place of $F$ and $A_3$ in place of $a_E$. The proof is unchanged with these substitutions.

\subsection{Reduction orders for the two carriers}
\begin{lemma}\label{lem:common-CM-reduction-orders}
Let $F_0/\Q$ be a real quadratic field, let $c\in F_0^\times$,
and suppose
\[
 N_{F_0/\Q}(c)=m^2\ell
\]
for a nonzero rational number $m$ and an odd rational prime
$\ell$. Let $p$ be a rational prime of good reduction for
$E_c:y^2=x^3-c^2x$, splitting in $F_0$, and suppose
$p\equiv1\pmod{16}$ and $\bigl(\frac\ell p\bigr)=-1$.
Choose $a,b\in\Z$ with $p=a^2+b^2$ and
$a+bi\equiv9\pmod{16\Z[i]}$.  Then the reduction orders
at the two primes of $F_0$ over $p$ are, in some order,
\[
       d_+(p)=(a+1)^2+b^2\ \ \t{and}\ \ 
       d_-(p)=(a-1)^2+b^2.
\]
\end{lemma}
\begin{proof}
The curve $E_c$ is the quadratic twist of
$E_1:y^2=x^3-x$ by $c$.  The classical $j=1728$
complex-multiplication formula gives the Frobenius trace
$2a$ for $E_1$ with the indicated Gaussian normalization
(the unordered pair of orders is independent of changing
the sign of $a$).  If $\mathfrak p,\mathfrak p'$ are the two
primes above $p$, the product of the two quadratic twist
characters equals
\[
 \left(\frac{c\bmod\mathfrak p}{p}\right)
 \left(\frac{c\bmod\mathfrak p'}{p}\right)
 =\left(\frac{N_{F_0/\Q}(c)}p\right)
 =\left(\frac\ell p\right)=-1.
\]
Thus the traces are opposite and the orders are
$p+1\pm2a=(a\pm1)^2+b^2$, as claimed.
\end{proof}

\subsection{A single finite-algebra witness}
The following compiler is the specialization-safe finite-algebra
construction of \cite[Section~19]{Sun38}, applied
with the present base field and finitely many chart branches.

\begin{lemma}\label{lem:finite-algebra-compiler}
Over the real quadratic field $K$, all the finite choices
of elliptic ordinates, square roots, chart branches, and
nonvanishing guards are encoded with one additional
$K$-unknown, after the freely chosen abscissa components
and structural parameters have been fixed.
\end{lemma}
\begin{proof}
We spell out the specialization-safe step because the variable count
uses it.  First stratify by the finitely many possible actual degrees,
squarefree parts, gcd degrees and nonzero pivots of the fixed-degree
polynomials occurring in the charts.  On each stratum the remaining
finite choices form a finite reduced algebra $A$ over the coefficient field of the stratum. Here only fixed-degree monic equations for ordinates and square roots, group-law charts, and finite branches occur. Stratifying by the ranks and pivots of their multiplication tables and taking squarefree parts gives the specialized reduced finite algebra uniformly.

Let $Z_1,\ldots,Z_r\in A$ be the coordinate functions of the
finite algebraic choices, including all of the roots and chart
coordinates under consideration.  By the definition of the affine
coordinate algebra, these elements generate $A$ over the coefficient
field.  Write $D=\dim A$, and let $\Omega$ be the set of its $D$
distinct geometric points.  Since $A$ is finite and reduced in
characteristic zero, its geometric fibre is reduced and $|\Omega|=D$.

We construct a primitive element without requiring that an
\emph{intermediate} pair of coordinates distinguish all of $\Omega$.
Put $T_1=Z_1$.  Suppose inductively that $T_{j-1}$
distinguishes exactly those pairs of geometric points that are
distinguished by at least one of $Z_1,\ldots,Z_{j-1}$.  For an integer
$b$ put $T_j=T_{j-1}+bZ_j$.  If a pair $x,y\in\Omega$
is distinguished by $Z_1,\ldots,Z_j$, then it is distinguished
by the ordered pair $(T_{j-1},Z_j)$.  The equality
\[
 T_{j-1}(x)+bZ_j(x)=T_{j-1}(y)+bZ_j(y)
\]
excludes at most one integer $b$: if the values of $Z_j$ agree,
then the values of $T_{j-1}$ differ, and the equality is
impossible; otherwise the equality determines $b$ uniquely.
Therefore, among the fixed integers
\[
       1\le b\le 1+\binom D2,
\]
there is at least one choice for which every pair distinguished
by the first $j$ coordinates remains distinguished by $T_j$.
Pairs on which the first $j$ coordinates agree may of course
remain indistinguishable at this stage; this is harmless.

Iterating over all coordinates gives $T_r\in A$ separating the
whole of $\Omega$, because the coordinates generate the reduced
algebra $A$.  Consequently, the characteristic polynomial
$P(T)=\det(T-M_{T_r})$
has degree $D$ and it is squarefree, and the map from the finite
scheme to the affine line defined by $T_r$ is an isomorphism
onto the zero scheme of $P$.  All choices of the integers
$b$ lie in a fixed finite list determined solely by $D$.
For an arbitrary specialization, one branches over the finitely
many resulting tuples of coefficients and retains exactly the
branches on which the \emph{final} discriminant of $P$ is
nonzero.  Each valid reduced specialization has at least one
such branch by the preceding argument; a zero intermediate
discriminant is never used as a reason to discard it.  The
construction is Galois equivariant.  Hence a root of the final
$P$ in $K$ corresponds exactly to a $K$-rational point of the
finite algebra.  This supplies the primitive-element step of
\cite[Sections~10.3--10.5]{Sun38} uniformly over the three
real quadratic fields used here.
It remains to impose equations and nonvanishing guards that may
\emph{depend on the finite algebraic choice}.  This point must be
settled before clearing denominators.  On a fixed reduced branch let
$P(T)$ be its squarefree primitive-element polynomial.  Write the
additional equations and guards as $E_i(T)=0$ and $G_j(T)\ne0$,
respectively.  Every denominator depending on $T$ is itself
included among the $G_j$, so no forbidden pole becomes a root.
After clearing denominators and stratifying by the degrees and
nonzero subresultant pivots, these polynomials and the following
monic gcds specialize with the stated degrees.  Put
\[
 Q(T)=\gcd\bigl(P(T),E_1(T),\ldots,E_m(T)\bigr),
 \qquad
 R(T)=\frac{Q(T)}{\gcd\bigl(Q(T),\prod_jG_j(T)\bigr)}.
\]
Because $P$ is squarefree, the geometric roots of $R$ are
\emph{exactly} the roots representing choices satisfying every
equation and every nonvanishing guard.  This filtering is important:
a guard that vanishes at one root must not eliminate a
different valid root of the same finite algebra.  On each fixed-degree stratum,
the subresultant formulas express the coefficients of $R$ as rational
functions in the freely chosen parameters.  All denominator,
subresultant-pivot, and branch guards that remain at this stage
therefore concern the \emph{parameters alone}.  Empty branches,
for which $R=1$, are discarded.

Write $d=\deg R\ge1$.  Include the $d+1$ fixed shifts
$R(T-b)$, for $0\le b\le d$, as further finite branches.
At any specialization $R$ has at most $d$ roots; hence some
shift makes all its roots nonzero.  Conversely, a root of any
shift corresponds to a root of $R$, and so to a valid finite
choice.  On one shifted branch write the monic polynomial as
\[
 R_b(T)=T^d+\sum_{j=1}^d\frac{N_j}{B}T^{d-j},
\]
where $B$ clears the coefficient denominators, and let $H$ be
the product of the remaining \emph{parameter-only} nonzero
chart and subresultant guards.  Define
\[
 \Psi(U)=1+\sum_{j=1}^d N_j B^{2j-1}H^jU^j.
\]
If $BH=0$, the polynomial $\Psi$ is identically $1$,
so an invalid specialization contributes no solution.  For
$BH\ne0$, set $S=B^2H$.  The identity
\[
 \Psi(U)=(SU)^dR_b\left(\frac1{SU}\right)
 \qquad(U\ne0)
\]
shows that a root of $\Psi$ is equivalent to a root of $R_b$;
moreover $\Psi(0)=1$.  The polynomial equalities defining
this parameter stratum are combined with $\Psi=0$ by taking
the sum of their squares, which is valid because $K$ is
formally real.  A finite disjunction over all strata, chart
choices and fixed shifts is their product, with the \emph{same}
unknown $U$.  This gives one additional field unknown and
retains equivalence after every valid specialization.
The construction uses the primitive-element and guard-absorption
methods of \cite[Section~19]{Sun38}.
\end{proof}

\section{Carrier and height data over the field \texorpdfstring{$\Q(\sqrt3)$}{Q(sqrt(3))}}\label{sec:sqrt3}

In this section, we put
\[
 K_3=\Q(\sqrt3),\ \ F_3=\Q(\sqrt2),\ \ \t{and}\ \ 
 L_3=K_3F_3=\Q(\sqrt2,\sqrt3).
\]
The nontrivial element of $\Gal(L_3/K_3)$ restricts to the
nontrivial automorphism of $F_3/\Q$.  Thus the relative trace and
conjugate tests used in the preceding construction are again
available, now in the basis $1,\sqrt2$ of $L_3/K_3$.

Put
\[
 e_3=5+2\sqrt2\ \ \t{and}\ \ A_3=3e_3=15+6\sqrt2,
\]
and consider
\begin{equation}\label{eq:sqrt3-curve}
 E_3:\ y^2=x^3-A_3^2x
\end{equation}
over $F_3$.

\begin{lemma}\label{lem:sqrt3-carrier}
For the curve \eqref{eq:sqrt3-curve}, we have
\[
\begin{gathered}
 \operatorname{rank}E_3(F_3)=\operatorname{rank}E_3(L_3)=1,\\
 E_3(L_3)_{\rm tors}=E_3[2]\ \ \t{and}\ \
 2E_3(L_3)\subset E_3(F_3).
\end{gathered}
\]
Moreover, for a suitable nontorsion point $P_3\in E_3(F_3)$, we have
\begin{equation}\label{eq:sqrt3-carrier}
             2E_3(L_3)=\Z P_3.
\end{equation}
The square classes of $A_3$ and its nontrivial $F_3/\Q$
conjugate are independent in $F_3(i)^\times/F_3(i)^{\times2}$.
\end{lemma}
\begin{proof}
For $c\in F_3^\times$ write
$E_c:y^2=x^3-c^2x$ and $E'_c:y^2=x^3+4c^2x$.
The standard two-isogeny descent gives
\[
 2^{\operatorname{rank}E_c(F_3)}
 =\frac{|\operatorname{Im}\alpha_c|\cdot
        |\operatorname{Im}\alpha'_c|}{4};
\]
for a square class $d$ the two covering equations are
\begin{equation}\label{eq:sqrt3-coverings}
 W^2=dU^4-\frac{c^2}{d}V^4\ \ \t{and}\ \ 
 W^2=dU^4+\frac{4c^2}{d}V^4.
\end{equation}
In applying these equations to a descent class, $(U,V)$ must
be nonzero.  At a finite completion we use the primitive
normalization of Lemma~\ref{lem:primitive-quartic-local};
in particular, the identically zero solution cannot establish
membership in a descent image.
The ring $\Z[\sqrt2]$ has class number one and unit group
$\langle-1,1+\sqrt2\rangle$.  The element $e_3$ has norm $17$,
and $3$ is inert in $F_3$.  The usual valuation argument outside
$2c$, together with the two real embeddings for the dual curve,
reduces the possible descent classes to
\[
\begin{array}{ll}
 \operatorname{Im}\alpha_{e_3}
   \subset\langle-1,1+\sqrt2,e_3\rangle,
 &\operatorname{Im}\alpha'_{e_3}\subset\langle e_3\rangle,\\
 \operatorname{Im}\alpha_{A_3}
   \subset\langle-1,1+\sqrt2,3,e_3\rangle,
 &\operatorname{Im}\alpha'_{A_3}\subset\langle3,e_3\rangle.
\end{array}
\]
At the prime $(e_3)$ one has $\sqrt2\equiv6\pmod{17}$.
Hence $1+\sqrt2\equiv7$ and $3\equiv3\pmod{17}$, both
quadratic nonresidues, while $-1$ is a square.  The first
covering in \eqref{eq:sqrt3-coverings} therefore excludes
$1+\sqrt2$ from $\operatorname{Im}\alpha_{e_3}$ and excludes
$1+\sqrt2$ and $3$ from $\operatorname{Im}\alpha_{A_3}$.
The two-torsion classes give
\[
 \operatorname{Im}\alpha_{e_3}=\langle-1,e_3\rangle
 \ \ \t{and}\ \ 
 \operatorname{Im}\alpha_{A_3}
       \subset\langle-1,A_3,3(1+\sqrt2)\rangle.
\]
The remaining class $3(1+\sqrt2)$ is excluded at the dyadic
prime.  Indeed its covering is
\begin{equation}\label{eq:sqrt3-dyadic-cover}
 W^2=(3+3\sqrt2)U^4-(21+39\sqrt2)V^4.
\end{equation}
Modulo $16\Z[\sqrt2]$, a dyadic unit fourth power has the form
$1+4k\sqrt2$.  If exactly one of $U,V$ is a unit, the
$\sqrt2$-coefficient in \eqref{eq:sqrt3-dyadic-cover} is odd,
whereas that of a square is even.  If both are units, the right
side is congruent to
\[
 14+8(k+h)+(12+12(k+h))\sqrt2\pmod{16}.
\]
For $k+h$ even its constant term is $14$; for $k+h$ odd it is
$6$, and a square with constant term $6$ modulo $16$ has
$\sqrt2$-coefficient $4$ or $12$ modulo $16$, not $0$ or $8$.
Thus this class has no primitive dyadic solution, and
$|\operatorname{Im}\alpha_{A_3}|=4$.  The finite-ring
check can be made exhaustive: modulo $16$ the unit fourth
powers are exactly $1+4k\sqrt2$ ($k\bmod4$),
whereas nonunit fourth powers are $0$ and $4$.  These are
all the possibilities under the primitive normalization.

For the dual descent of $E_{e_3}$, the only possible nontrivial
class is $e_3$.  Its covering is
\[
 W^2=e_3(U^4+4V^4).
\]
We exclude it modulo $64\Z[\sqrt2]$.  In
$R=(\Z/64\Z)[X]/(X^2-2)$, a unit fourth power has the form
\[
 (1+16a)+4bX\qquad(a\bmod4,\ b\bmod16),
\]
whereas a nonunit fourth power is one of
$0,16$, or $4+16bX$ with $b\bmod4$.  If $U$ is a unit, the
constant coordinate of
$(5+2X)(U^4+4V^4)$ is one of
\[
 5,9,21,25,37,41,53,57,
\]
and its $X$-coordinate is $2\pmod4$.  Constants
$5,21,37,53$ do not occur for squares in $R$, while a square
with constant coordinate $9,25,41$ or $57$ has $X$-coordinate
divisible by four.  If $U$ is not a unit, primitivity forces $V$
to be a unit; the possible constant coordinates are $20,36,40$.
The values $20$ and $40$ do not occur for squares, and a square
with constant coordinate $36$ has $X$-coordinate divisible by
$16$, whereas the right side has $X$-coordinate $8\pmod{16}$.
Thus there is no primitive solution modulo $64$, and hence
$\operatorname{Im}\alpha'_{e_3}=\{1\}$ and
$\operatorname{rank}E_{e_3}(F_3)=0$.  More explicitly, the
nonunit fourth powers in this residue ring are precisely
$0$, $16$, and $4+16b\sqrt2$ ($b\bmod4$), and
the unit fourth powers are precisely
$(1+16a)+4b\sqrt2$ ($a\bmod4$, $b\bmod16$).
The two cases displayed above exhaust these finite lists.
For $E'_{A_3}$, the class $3$ is excluded at $(e_3)$, and the
class $e_3$ would require a fourth power congruent to
$-36\equiv15\pmod{17}$; the fourth powers modulo $17$ are
$1,4,13,16$.  Thus
$|\operatorname{Im}\alpha'_{A_3}|\le2$ and
$\operatorname{rank}E_3(F_3)\le1$.

Now put $B_3=9+10\sqrt2$.  Since $A_3^2-B_3^2=16$, we have
\[
 P_0=\left(\frac{A_3^2}{4},\frac{A_3^2B_3}{8}\right)
      \in E_3(F_3).
\]
The prime $73$ splits completely in $L_3$.  At the two primes
of $F_3$ with $\sqrt2\equiv\pm32\pmod{73}$, direct reduction
gives
\[
     \#E_3(\F_{73})=80\ \ \t{and}\ \ 
     \#E_3^\sigma(\F_{73})=68.
\]
Good reduction, together with the standard formal-group
injectivity for the $73$-primary torsion, yields that
$|E_3(L_3)_{\rm tors}|$ divides $\gcd(80,68)=4$.
Since $E_3[2]\subset E_3(F_3)$, equality holds.  The ordinate
of $P_0$ is nonzero, so $P_0$ is nontorsion; hence
$\operatorname{rank}E_3(F_3)=1$.

The quadratic twist of $E_3$ by $3$ is $E_{3A_3}=E_{9e_3}$,
which is $F_3$-isomorphic to $E_{e_3}$.  The quadratic-twist
rank decomposition therefore gives
$\operatorname{rank}E_3(L_3)=1$.  The involution of $L_3/F_3$
fixes a nontorsion point, hence acts trivially on the rank-one
free quotient.  It follows that $2E_3(L_3)\subset E_3(F_3)$,
and Lemma~\ref{lem:common-biquadratic-carrier}, applied with
$k=K_3$, $F_0=F_3$, and $L_0=L_3$, now gives
\eqref{eq:sqrt3-carrier} and the two-coordinate representation
of the half-points.

Finally, $e_3$ and $e_3^\sigma$ generate the two distinct
primes of $F_3$ above $17$.  The valuations of $A_3$ and
$A_3^\sigma$ at these two primes have opposite parity.  Since
$17\equiv1\pmod4$, this parity remains visible after adjoining
$i$, proving the stated square-class independence.
\end{proof}

We next construct the local geometry.  Put
\[
 q_3=347\ \ \t{and}\ \ s_3=359.
\]
Choose an integer $v_*>3$ with
\[
 v_*\equiv1+q_3\pmod{q_3^2}\ \ \t{and}\ \ 
 v_*\equiv2\pmod{s_3^2},
\]
and set
\[
 \delta_* =\frac{(v_*^2-1)(v_*^2-9)}{16v_*^2},\qquad
 u_1=\frac{v_*^2+3}{4v_*},\qquad
 u_4=\frac{v_*^2-3}{2v_*}.
\]
Thus $u_1^2=1+\delta_*$ and $u_4^2=1+4\delta_*$.  Choose a
positive rational number $\lambda<u_4^{-1}$ satisfying
\[
 \lambda\equiv1+q_3\pmod{q_3^2}\ \ \t{and}\ \ 
 \lambda\equiv8\cdot7^{-1}+s_3\pmod{s_3^2}.
\]
 Rational Chinese
remainder approximation gives such a $\lambda$.  Define
\begin{equation}\label{eq:sqrt3-local-data}
\begin{gathered}
 t_*=\frac{\lambda^{-2}-1}{\delta_*},\qquad
 b_*=\frac{\lambda^2\delta_*}{4q_3},\\
 x_3(z)=t_*+\frac{s_3^2z+s_3^4b_*}{q_3z^2},\qquad
 j_3(z)=\frac{16x_3(z)^2}{(x_3(z)-4)^3}.
\end{gathered}
\end{equation}
For $c\in\{0,1,4\}$ put
\[
 F_c(z)=q_3(t_*-c)z^2+s_3^2z+s_3^4b_*.
\]

Let
\begin{equation}\label{eq:sqrt3-contacts}
 \mathcal C_3=\{0\}\cup
 \{z\in\Q:F_0(z)F_1(z)F_4(z)=0\}\ \t{and}\ 
 P_{3,\mathcal C}(z)=\frac{F_1(z)F_4(z)}{4q_3^2}.
\end{equation}
In choosing $\lambda$ above, require in addition
\[
       F_0(1-q_3)F_1(1-q_3)F_4(1-q_3)\ne0.
\]
This excludes only finitely many possible values of $\lambda$.
Indeed, for fixed $v_*$ and each $c\in\{0,1,4\}$, multiplying
$F_c(1-q_3)$ by $\delta_*\lambda^2$ gives a quadratic
polynomial in $\lambda^2$ with nonzero leading coefficient
$s_3^4\delta_*^2/(4q_3)$. The specified congruence classes
and positive real interval still contain infinitely many rational
choices.  Thus $1-q_3\notin\mathcal C_3$, so the artificial
contact $-1$ used below is distinct from every member of
$\Theta_3=\{q_3/(c-1):c\in\mathcal C_3\}$.
Let $S_3$ contain all prime ideals of $K_3$ over rational
primes at most $324$, all primes over $q_3s_3$, and all prime
divisors of the denominators and pairwise differences of the
contacts.  Define $M_{\mathcal C_3,S_3}$ by the same positive
odd-contact rule as before.

\begin{lemma}\label{lem:sqrt3-seven}
The set $\mathcal C_3$ consists of seven distinct rational
points, all in $(-\infty,0]$, and all have residual orbifold
order two.  For $z>0$ one has $x_3(z)>4$ and $P_{3,\mathcal C}(z)>0$.
Moreover,
\begin{equation}\label{eq:sqrt3-squareclass}
 -(27j_3(z)+16)=-P_{3,\mathcal C}(z)
 \left(\frac{8(x_3(z)+8)}{z^2(x_3(z)-4)^2}\right)^2.
\end{equation}
If $\mathfrak q_3\mid347$ and $\mathfrak s_3\mid359$ are
fixed primes of $K_3$, then
\begin{equation}\label{eq:sqrt3-local-conditions}
\begin{aligned}
 z\equiv1\pmod{\mathfrak q_3}
   &\Longrightarrow v_{\mathfrak q_3}(j_3(z))=1,\\
 z\in\OO_{K_{3,\mathfrak s_3}}^\times
   &\Longrightarrow
   v_{\mathfrak s_3}\left(j_3(z)+\frac{16}{27}\right)=1.
\end{aligned}
\end{equation}
\end{lemma}
\begin{proof}
The discriminants of $F_0,F_1,F_4$ are respectively
\[
 s_3^4\lambda^2,
 \qquad s_3^4\lambda^2u_1^2,
 \qquad s_3^4\lambda^2u_4^2.
\]
They are nonzero rational squares.  The three quadratics have
no common root, and their coefficients are positive, so their six
roots are distinct negative rational numbers.  At $z=0$ the map
$x_3$ has local degree two above the order-four point at infinity;
this changes its residual order to two.  The six roots above
$x=0,1,4$ have residual order two by the same factorization
\[
 j(x)+\frac{16}{27}
 =\frac{16(x-1)(x+8)^2}{27(x-4)^3}
\]
used in Proposition~\ref{prop:seven}.  Since
$0<\lambda<u_4^{-1}$, one has $t_*>4$; all other terms in
$x_3(z)$ are positive for $z>0$, proving $x_3(z)>4$.
The identity
$P_{3,\mathcal C}(z)=z^4(x_3(z)-1)(x_3(z)-4)/4$ gives the positivity and,
together with the displayed formula for $j+16/27$, proves
\eqref{eq:sqrt3-squareclass}.

Modulo $q_3^2$ the defining congruences give
\[
 \delta_*\equiv-q_3,\ \ 
 \lambda^{-2}-1\equiv-2q_3,
 \ \ t_*\equiv2\pmod{q_3},
 \ \ b_*\equiv-1/4\pmod{q_3}.
\]
Since $s_3\equiv12\pmod{q_3}$,
$s_3^2+s_3^4b_*\not\equiv0\pmod{q_3}$.  Hence
$z\equiv1\pmod{\mathfrak q_3}$ gives
$v_{\mathfrak q_3}(x_3(z))=-1$, and therefore
$v_{\mathfrak q_3}(j_3(z))=1$.
In view of the congruences $v_*\equiv2\pmod{s_3^2}$ and
$\lambda\equiv8/7+s_3\pmod{s_3^2}$, we have
\[
 v_{s_3}(t_*-1)=1.
\]
If $z$ is an $\mathfrak s_3$-unit, the two remaining terms in
$x_3(z)-t_*$ have $s_3$-adic valuation at least two.  Thus
$v_{\mathfrak s_3}(x_3(z)-1)=1$, while $x_3(z)-4$ and
$x_3(z)+8$ are units.  The factorization above gives the second
formula in \eqref{eq:sqrt3-local-conditions}.
\end{proof}

Let
\[
 \Theta_3=\left\{\frac{q_3}{c-1}:\ c\in\mathcal C_3\right\}
 \ \ \t{and}\ \ \xi=1.
\]
This is a fixed set of seven distinct negative rational numbers.

\begin{lemma}\label{lem:sqrt3-dyadic}
Let $K_{3,2}=\Q_2(\sqrt3)$, $\pi=1+\sqrt3$ and
$\varepsilon=2+\sqrt3$.  Given a root of unity $\omega$ of
$2$-power order and $s\in\{\pm1\}$, there is a finite-order
character
\[
 \theta:K_{3,2}^\times\longrightarrow\mu_{2^\infty},
 \qquad \theta(\varepsilon)=\omega,
 \qquad \theta(-1)=s.
\]
If $\operatorname{ord}(\omega)\le2^m$, it can be chosen with
order at most $2^{m+3}$ and conductor exponent at most $6+2m$.
\end{lemma}
\begin{proof} Clearly. $\pi^2=2\varepsilon$ and
$U_6=1+\pi^6\OO_{K_{3,2}}=1+8\OO_{K_{3,2}}$.  Moreover,
\[
 \varepsilon^4-1=8(12+7\sqrt3),
\]
and the last factor is a $\pi$-unit.  The $2$-adic logarithm
identifies $U_6$ with the free rank-two $\Z_2$-module
$\pi^6\OO_{K_{3,2}}$, with $\log(\varepsilon^4)$ primitive.
Prescribe the value $\varepsilon^4\mapsto\omega^4$ on one direct
summand and take the character trivial on a complement.  It is
trivial on $U_{6+2\max(m-2,0)}$.  The class of $\varepsilon$ has
order four modulo $U_6$, while $-1$ does not lie in its cyclic
subgroup.  Hence the character extends first to
$\langle U_6,\varepsilon,-1\rangle$ with the prescribed sign and
then across the remaining index-four quotient.  This increases
the order by at most four.  Taking $\theta(\pi)=1$ completes the
extension and proves the conductor bound.
\end{proof}

The class number of $K_3$ is one and
$\OO_{K_3}^\times=\{\pm(2+\sqrt3)^n:n\in\Z\}$.
Repeating the class-field-theoretic proof of
Theorem~\ref{thm:cochain}, using Lemma~\ref{lem:sqrt3-dyadic}
at the unique dyadic prime, gives
\begin{equation}\label{eq:sqrt3-cochain-bounds}
\begin{gathered}
 m\le8B_{T_0},\ 
 N_{K_3/\Q}(\mathfrak f_\chi)
 \le64B_{T_0}^2\prod_{\mathcal P\in T_-}N\mathcal P,\ 
 [L':K_3]\le16B_{T_0}[L_0:K_3].
\end{gathered}
\end{equation}
Indeed, the local characters away from two may be chosen exactly
as in that proof.  Their product on the totally positive generator
$2+\sqrt3$ and on $-1$ is cancelled by
Lemma~\ref{lem:sqrt3-dyadic}; class number one then globalizes the
local character.  The proof of the Bockstein identity and the final
quadratic twist is unchanged.  In particular the degree and
conductor are polynomially bounded by the odd-contact product.

\begin{lemma}\label{lem:sqrt3-height}
There are positive constants $c_3$ and $H_3$ such that
\begin{equation}\label{eq:sqrt3-height}
 h(z)\le H_3M_{\mathcal C_3,S_3}(z)^{c_3}
\end{equation}
for any $z\in K_3\setminus\mathcal C_3$ with $z\in\OO_{K_{3,\mathfrak s_3}}^\times$,
$ z\equiv1\pmod{\mathfrak q_3},$ and $ \sigma(P_{3,\mathcal C}(z))>0$ for all 
$\sigma:K_3\hookrightarrow\R.$
\end{lemma}
\begin{proof}
The discriminant-six quaternion algebra, the maximal order, the
principal PEL datum and the full-three cover of
Lemmas~\ref{lem:source-pel}--\ref{lem:PEL-lattice-transport}
are defined over $\Q$, and hence base change to $K_3$.  By the
first formula in \eqref{eq:sqrt3-local-conditions}, the first two conclusions of
Proposition~\ref{prop:pel-inertia} give an order-four projective
inertia element represented by $\gamma=1+i+j$ at
$\mathfrak q_3$.  Since $347\equiv-1\pmod4$, local Frobenius
inverts this rotation.  The local decomposition group therefore
contains a dihedral group of order eight and the residual
representation is already absolutely irreducible over $G_{K_3}$.
The proof of Lemma~\ref{lem:global-D8}, with $347$ in place of
$367$, now shows
\[
             \mathbb P\bar r_3(G_{K_3})\simeq D_8,
\]
and again identifies the rotation character with the
three-adic Atkin--Lehner parity $\delta_3$.

The quadratic intermediate cover of the full Atkin--Lehner
quotient is $\Q(j)(\sqrt{-27j-16})$
by \cite[Theorem~3.10, Corollary~3.11 and Theorem~3.17]{BabaGranath}.
The second formula in \eqref{eq:sqrt3-local-conditions} therefore
makes $\delta_3$ nontrivial on $I_{\mathfrak s_3}$.  By
Lemma~\ref{lem:projective-normalizer}, $\delta_3$ descends to the
actual projective residual $D_8$-quotient with kernel precisely the
rotation subgroup.  Hence the image of $I_{\mathfrak s_3}$ contains
a reflection.
On the other hand, $
 K_3(\zeta_3)=K_3(i).$
The three quadratic subfields of
$\Q(\sqrt3,i)$ have discriminants $12,-4,-3$, and the
biquadratic discriminant formula gives
$|\operatorname{disc}\Q(\sqrt3,i)|=144
 =|\operatorname{disc}K_3|^2$.  Thus $K_3(i)/K_3$ is
unramified at every finite prime.  Both finite inertia groups above
therefore lie in $G_{K_3(\zeta_3)}$; their images contain a
fourfold rotation and a reflection.  Consequently
\begin{equation}\label{eq:sqrt3-cyclo-D8}
 \mathbb P\bar r_3(G_{K_3(\zeta_3)})\simeq D_8,
\end{equation}
so the restriction is absolutely irreducible.

By \eqref{eq:sqrt3-squareclass} the rotation quadratic field is
$ K_3(\sqrt{-P_{3,\mathcal C}(z)}),$
which is CM under the stated positivity condition.  The proof of
Lemma~\ref{lem:nonordinary} is local in the quaternion division
algebra at three and does not require the base local field to be
unramified over $\Q_3$.  It therefore applies to
$K_{3,\mathfrak p_3}=\Q_3(\sqrt3)$ and shows that the normalized
$r_3$ is potentially Barsotti--Tate and not potentially ordinary.
The CM automorphic-induction argument in the proof of
Theorem~\ref{thm:r3-modular}, followed by Gee's corrected theorem
\cite[Theorem~1.1]{GeeErratum}, applies verbatim.  The erratum
allows ramification of the prime above three, and the
potentially-ordinary matching condition is vacuous here because
$r_3$ is not potentially ordinary.  Thus $r_3$ is modular of
parallel weight two.

At $\mathfrak q_3$ the fourfold inertia and Frobenius inversion
make the two-dimensional Weil--Deligne parameter irreducible.
Lemma~\ref{lem:skinner-supercuspidal} gives a supercuspidal local
component before Jacquet--Langlands is invoked.  The corrected
Carayol--Saito realization in Proposition~\ref{prop:carayol-tate}
and Lemma~\ref{lem:cohom-normalization} identify the corresponding
Jacobian Tate constituent with the same $r_3$.
The bounds \eqref{eq:sqrt3-cochain-bounds}, the fixed-degree
full-three splitting field, and the proof of
Lemma~\ref{lem:dyadic-conductor} give
\[
 N_{K_3/\Q}(\mathfrak n(\Pi_z))
 \ll M_{\mathcal C_3,S_3}(z)^A
\]
for a fixed $A$.  Proposition~\ref{prop:height} now yields
\eqref{eq:sqrt3-height} outside the CM locus; the CM fibres are
absorbed by Lemma~\ref{lem:cm-bounded-degree}.
\end{proof}

For the prime-pattern step, if $\pi=a+bi$ is a Gaussian prime of
prime norm $p=a^2+b^2$, put

\section{Carrier and height data over the field \texorpdfstring{$\Q(\sqrt6)$}{Q(sqrt(6))}}\label{sec:sqrt6}

In this section, we put
\[
 K_6=\Q(\sqrt6),\ \ F=\Q(\sqrt3),\ \ \t{and}\ \ 
 L_6=K_6F=\Q(\sqrt2,\sqrt3).
\]
We use the same curve $E/F$ and the same generator $P$ as in
Theorem~\ref{thm:carrier}.  The nontrivial element of
$\Gal(L_6/K_6)$ restricts to the nontrivial automorphism of
$F/\Q$; all traces and conjugate tests below are taken with
respect to this relative extension.

\begin{lemma}\label{lem:sqrt6-carrier}
For the curve $E(F)$ as in \eqref{eq:E}, we have
\[
\begin{gathered}
 \operatorname{rank}E(F)=\operatorname{rank}E(L_6)=1,\\
 E(L_6)_{\rm tors}=E[2]\ \ \t{and}\ \ 
 2E(L_6)=\Z P\subset E(F).
\end{gathered}
\]
Every point in $2E(L_6)$ can be represented by the abscissa
of a half-point in $E(L_6)$, using two $K_6$-coordinates.
\end{lemma}
\begin{proof}
Since $L_6=F(\sqrt2)$, this is the same extension of $F$
used in Sections~\ref{sec:descent} and~\ref{sec:stability}.
In particular, the rank calculation, the torsion computation
and the inclusion \eqref{eq:2EL} apply without change.
The common biquadratic carrier Lemma~\ref{lem:common-biquadratic-carrier}
then gives the cyclic group and coordinate assertions in the
basis $1,\sqrt3$ of $L_6/K_6$.
\end{proof}

We next replace the distinguished rational prime $367$ by $331$.
Let
\begin{equation}\label{eq:six-local-data}
\begin{gathered}
 q_6=331,\qquad \mathfrak q_6=(25+7\sqrt6)\subset\OO_{K_6},
 \qquad D_6=q_6^2+3,\\
 a_6=\frac{(q_6-3)(q_6+1)}{D_6},\qquad
 b_6=\frac{(q_6-1)(q_6+3)}{D_6^2},\\
 x_6(z)=\frac{a_6z+b_6}{q_6z^2},\qquad
 j_6(z)=\frac{16x_6(z)^2}{(x_6(z)-4)^3}.
\end{gathered}
\end{equation}
Define the seven contacts
\begin{equation}\label{eq:six-contacts}
\begin{split}
 \mathcal C_6=\bigg\{\infty,0,
 -\frac{(q_6-1)(q_6+3)}{(q_6-3)(q_6+1)D_6},
 \frac{q_6-1}{D_6},-\frac{q_6+3}{q_6D_6},
\frac{(q_6-1)(q_6+3)}{4q_6D_6},-\frac1{D_6}\bigg\}
\end{split}
\end{equation}
and put
\begin{equation}\label{eq:six-PC}
 P_6(z)=\frac{(a_6z+b_6-q_6z^2)
                   (a_6z+b_6-4q_6z^2)}{4q_6^2}.
\end{equation}
Let $S_6$ contain all prime ideals of $K_6$ over rational
primes at most $324$, all prime ideals over $331$, and all
prime divisors of the denominators and pairwise differences
of the finite contacts.  Define $M_{\mathcal C_6,S_6}(z)$
by the same positive odd-contact rule as before.  The
corresponding local condition is
$z\equiv1\pmod{\mathfrak q_6}$.

\begin{lemma}\label{lem:sqrt6-seven}
The contacts \eqref{eq:six-contacts} are rational, distinct,
and have residual orbifold order two.  If
$z\equiv1\pmod{\mathfrak q_6}$, then
$v_{\mathfrak q_6}(j_6(z))=1.$
Moreover, with $\theta_0=0$ and
$\theta_i=q_6/(c_i-1)$ for the six finite contacts
$c_i$ of \eqref{eq:six-contacts},
\begin{equation}\label{eq:six-range}
 \Theta_6=\{\theta_0,\ldots,\theta_6\}
 \subset(-332,0]\cap\Q,
\end{equation}
and the denominators of these seven numbers are units at
$\mathfrak q_6$.  For $1/2\le z\le3/2$ one has
$0<x_6(z)<1$ and $P_6(z)>0$.
\end{lemma}
\begin{proof}
The rational identities
\[
 a_6^2+4q_6b_6=1\ \ \t{and}\ \ 
 a_6^2+16q_6b_6
 =\left(\frac{q_6^2+6q_6-3}{q_6^2+3}\right)^2
\]
are the identities used in the proof of
Proposition~\ref{prop:seven}.  Together with the
factorization \eqref{eq:factor-j}, they give precisely the
seven contacts in \eqref{eq:six-contacts}, with the same
local degrees and residual orders.

Since
\[
 N_{K_6/\Q}(25+7\sqrt6)=25^2-6\cdot7^2=331,
\]
$\mathfrak q_6$ has residue field $\F_{331}$, in which
$\sqrt6\equiv91$.  Reduction of \eqref{eq:six-local-data}
gives
\[
 a_6\equiv-1\pmod{\mathfrak q_6}\ \ \t{and}\ \ b_6\equiv-\frac13\equiv110
 \pmod{\mathfrak q_6}.
\]
Thus $a_6z+b_6\equiv109\pmod{\mathfrak q_6}$ on the
specified tube, and $v_{\mathfrak q_6}(x_6(z))=-1$.
As $v_{\mathfrak q_6}(x_6(z)-4)=-1$, the defining formula
for $j_6$ gives $v_{\mathfrak q_6}(j_6(z))=1$.
Direct substitution of the six rational contacts into
$\theta_i=q_6/(c_i-1)$ gives \eqref{eq:six-range};
all denominators are prime to $331$.  Finally,
$0<a_6<1$ and $0<b_6<1$, and the formula for $x_6$
shows $0<x_6(z)<1$ on $[1/2,3/2]$.
The identity
\[
 P_6(z)=\frac{z^4}4(x_6(z)-1)(x_6(z)-4)
\]
then proves the last assertion.  We shall also use the corresponding
square-class identity
\begin{equation}\label{eq:six-squareclass}
 -(27j_6(z)+16)
 =-P_6(z)
 \left(\frac{8(x_6(z)+8)}{z^2(x_6(z)-4)^2}\right)^2,
\end{equation}
which follows by substituting the preceding expression for $P_6$ in
$27j+16=16(x-1)(x+8)^2/(x-4)^3$.
\end{proof}

The next local calculation replaces
Lemma~\ref{lem:dyadic-interp}.  It also gives the required
version of Theorem~\ref{thm:cochain} over $K_6$.

\begin{lemma}\label{lem:sqrt6-dyadic}
Let $F_{6,2}=\Q_2(\sqrt6)$, $\pi=\sqrt6$ and
$\eta_6=5+2\sqrt6$.  For every root of unity $\omega$
of two-power order and every $s\in\{\pm1\}$ there is
a finite-order character
$\theta:F_{6,2}^{\times}\to\mu_{2^\infty}$ such that
\[
 \theta(\eta_6)=\omega\ \t{and}\ \ \theta(-1)=s.
\]
If $\omega$ has order at most $2^r$, the character can be
chosen of order at most $2^{\max(2,r+1)}$ and conductor
exponent at most $3+2r$.
\end{lemma}
\begin{proof}
We normalize $v_\pi(\pi)=1$.  Then
\[
 v_\pi(2)=2\ \ \t{and}\ \ v_\pi(\eta_6-1)=3.
\]
For $U_n=1+\pi^n\OO_{F_{6,2}}$, the logarithm gives
$\log:U_3\simeq\pi^3\OO_{F_{6,2}}$.
Furthermore, $\log\eta_6$ is primitive in this free
rank-two $\Z_2$-module.  Choose a complementary direct
summand and define a character of $U_3$ taking $\eta_6$
to $\omega$ and trivial on the complement.  This character
is trivial on $U_{3+2r}$.

Now $-1\in U_2\setminus U_3$ and
$[\OO_{F_{6,2}}^\times:U_3]=4$.  Thus the character
extends to $\langle-1\rangle U_3$ with the independently
prescribed value $s$.  That subgroup has index two in the
full unit group.  Extension across this index-two subgroup
increases the character order by at most a factor two;
when $r=0$ its order is at most four.  Put
$\theta(\pi)=1$ to extend to $F_{6,2}^\times$.
The restriction to $U_{3+2r}$ is unchanged, proving the
conductor bound.
\end{proof}

The class number of $K_6$ is one. Indeed its discriminant is
$24$, and the real quadratic Minkowski bound is
$\sqrt{24}/2<5/2$; the prime ideal above two is principal,
generated by $2+\sqrt6$ of norm $-2$.  Its unit group is
\[
 \OO_{K_6}^{\times}=\{\pm(5+2\sqrt6)^n:n\in\Z\}.
\]
Repeating the class-field-theoretic proof of
Theorem~\ref{thm:cochain}, with
Lemma~\ref{lem:sqrt6-dyadic} in place of
Lemma~\ref{lem:dyadic-interp}, gives for the corresponding
sign cocycle the same quantitative bounds
\begin{equation}\label{eq:six-cochain-bounds}
\begin{gathered}
 m\le4B_{T_0},\ 
 N_{K_6/\Q}(\mathfrak f_\chi)
 \le32B_{T_0}^{2}\prod_{\mathcal P\in T_-}N\mathcal P,\ 
 [L':K_6]\le8B_{T_0}[L_0:K_6].
\end{gathered}
\end{equation}
Here $T_0$ contains the unique dyadic prime of $K_6$, and
$T_-$ and $B_{T_0}$ are defined as in
Theorem~\ref{thm:cochain}.  In fact, the prescribed value of
$\eta_6$ has order at most $2B_{T_0}$, so
Lemma~\ref{lem:sqrt6-dyadic} gives a conductor exponent
at most $5+2\log_2B_{T_0}$, as required for
\eqref{eq:six-cochain-bounds}.

\begin{lemma}\label{lem:sqrt6-height}
There are positive constants $c_6$ and $H_6$ such that
\begin{equation}\label{eq:six-height}
 h(z)\le H_6 M_{\mathcal C_6,S_6}(z)^{c_6}
\end{equation}
for every $z\in K_6\setminus\mathcal C_6$ with
$z\equiv1\pmod{\mathfrak q_6}$ and
$\sigma(P_6(z))>0$ for each real embedding $\sigma$ of $K_6$.
\end{lemma}
\begin{proof}
We verify that the arguments of Sections~\ref{sec:seven}
through~\ref{sec:residual} apply to the new base field.
The discriminant-six quaternion algebra, the maximal order
$O$, its principal polarization and the full-three PEL
cover are defined over $\Q$.  The change-of-marking group,
the $144$-element normalizer action and the quadratic
intermediate cover in Lemma~\ref{lem:Gamma} are therefore
unchanged after base change to $K_6$.  The fine PEL model is smooth
over $\Z[1/30]$, hence has good reduction at $331$.  Moreover the
full-three generic fibre has $144$ points by Lemma~\ref{lem:Gamma}, so
its splitting extension has Galois degree at most $144!$; in particular
$331>144$ excludes wild inertia in that fixed finite permutation
extension.  The three branch values remain distinct modulo $331$.  Lemma~\ref{lem:sqrt6-seven} and the
proof of Proposition~\ref{prop:pel-inertia} identify local
projective inertia on the same normalized two-dimensional
$r_3$ with the right class of $\gamma=1+i+j$.
Since $331\equiv1\pmod3$, the final, cyclotomic conclusion
of Proposition~\ref{prop:pel-inertia} applies at this prime.
The element $\gamma$ has reduced norm two and is invertible
on the three-primary Tate module. Its projective characteristic
polynomial $T^2-2T+2$ has eigenvalue ratio of exact order
four modulo three.  Since
\[
 331\equiv-1\pmod4\ \ \t{and}\ \ 331\equiv1\pmod3,
\]
Frobenius inverts the inertia class and
$\zeta_3\in(K_6)_{\mathfrak q_6}$.  Hence
$\bar r_3|_{G_{K_6(\zeta_3)}}$ is absolutely irreducible.
The group-theoretic argument of
Lemma~\ref{lem:global-D8} gives
$\mathbb P\bar r_3(G_{K_6})\simeq D_8$.
Moreover, the square-class identity \eqref{eq:six-squareclass}
identifies the rotation field as
$K_6(\sqrt{-P_6(z)})$, 
which is CM under the stated positivity condition.

At the prime above three, the argument of
Lemma~\ref{lem:nonordinary} depends on the division
algebra $B\otimes\Q_3$, the four-dimensional rational
Tate module, and polarization symmetry, and does not
require the local field to be unramified over $\Q_3$.
In the present case there is a unique prime above three,
and
\[
 (K_6)_{\mathcal P_3}
 =\Q_3(\sqrt6)=\Q_3(\sqrt{-3})=\Q_3(\zeta_3),
\]
since $-2$ is a square in $\Q_3$.  The normalized $r_3$
is potentially Barsotti--Tate and not potentially ordinary,
and is odd at both real embeddings.  The CM automorphic
induction of the residual representation therefore supplies
its residual modularity.  Gee's corrected modularity-lifting
theorem \cite[Theorem~1.1]{GeeErratum} applies at $p=3$:
it allows the prime above three to ramify in the totally
real base field, and its potentially-ordinary matching
condition is vacuous for this $r_3$.  Thus the original
normalized $r_3$ is modular of parallel weight two.

At $\mathfrak q_6$, its Weil--Deligne parameter is
irreducible. Skinner's local--global compatibility
\cite{Skinner} therefore yields a supercuspidal local
component.  As in the proof of
Theorem~\ref{thm:height-final}, this precedes the
Jacquet--Langlands transfer and the Carayol--Saito
realization in the Jacobian of a compact Shimura curve.
The cohomological normalization of
Section~\ref{sec:height} identifies the resulting
Tate constituent with the same $r_3$.

Outside $S_6$, the pulled-back PEL cover is ramified
only at positive odd contacts, with inertia order at
most two. The character bounds \eqref{eq:six-cochain-bounds}
and the fixed-degree full-three splitting field give the
same polynomial conductor bound as in
Lemma~\ref{lem:dyadic-conductor}:
\[
 N_{K_6/\Q}(\mathfrak n(\Pi_z))
 \ll M_{\mathcal C_6,S_6}(z)^A
\]
for some fixed $A$.  At the finitely many primes of
$S_6$, any changes in residue degree or ramification
index affect only fixed constants.  The arguments of
Proposition~\ref{prop:height} (including the
Gaudron--R\'emond and Javanpeykar estimates) now give
\eqref{eq:six-height} outside the CM locus.  The CM fibres are
absorbed by Lemma~\ref{lem:cm-bounded-degree}.
\end{proof}

\section{Common arithmetic data for the two reductions}\label{sec:new-data}
In the remainder of the paper $K$ is any of the three fields in the main theorems, $L=\Q(\sqrt2,\sqrt3)$, and $F$ is the complementary quadratic field carrying $E$ and $P$: take $F=\Q(\sqrt2)$ for $K=\Q(\sqrt3)$, and $F=\Q(\sqrt3)$ otherwise. Write $\tau$ for the nonidentity automorphism of $L/K$. The cyclic carrier is $2E(L)=\Z P$. We use the seven-contact height estimates above, and write $c>0$ for an exponent valid for the chosen field. Write $S$ for the fixed height-exception set in $K$, enlarged as specified below by all prime ideals above finitely many fixed rational primes. This changes the height constants but not $c$.

Let $q=367,347,331$, respectively, for $K=\Q(\sqrt2),\Q(\sqrt3),\Q(\sqrt6)$. Write $\mathcal C_K$ for the contact set $\mathcal C$, $\mathcal C_3$ or $\mathcal C_6$ of \eqref{eq:contacts}, \eqref{eq:sqrt3-contacts} or \eqref{eq:six-contacts}, respectively, in the corresponding height theorem and $\Theta_K$ for its image under $z\mapsto q/(z-1)$ (infinity maps to zero). The coordinate change $t\mapsto1+q/t$ preserves normalized contact depths outside fixed primes and changes the height by $O(1)$. We also use a fixed polynomial $r_K\in\Z[T]$ of degree fourteen, as follows. Put $Q=q$ for the first and third fields and $Q=347\cdot359$ for the second. Set $D(T)=1+Q^2T^2$ and $\psi(T)=b+Q^2T/D(T)$, where $b=10q+1$ for the first and third fields and $b=1+Q$ for the second. Take
\[r_K(T)=C_K D(T)^{e_K}\prod_{a\in\mathcal C_K\setminus\{\infty\}}\bigl((b-a)D(T)+Q^2T\bigr)\]
and \[e_K=\begin{cases}1&\t{if}\ \infty\in\mathcal C_K,\\0&\text{otherwise},\end{cases}
\]
with $C_K$ a positive integer clearing denominators. Every factor is positive on $\R$, by $|Q^2T/D(T)|\le Q/2$ and the displayed contact sets. At the distinguished finite places $-1$ is a nonsquare in the residue field, so $QT/(1+Q^2T^2)$ is integral and $\psi(T)\equiv1$ there. The real inequalities in the height theorems hold for $\psi(T)$, and $h(\psi(T))=2h(T)+O(1)$. Away from the fixed coefficient, resultant and leading-coefficient primes, for integral $T$ at most one homogeneous pullback factor can have positive valuation; its depth is precisely the contact depth. For $v_{\mathcal P}(T)<0$, every factor has valuation $2v_{\mathcal P}(T)$, so the total valuation of $r_K(T)$ is even and there is no positive odd normalized contact. Thus the scalar height is bounded by
\begin{equation}\label{eq:new-scalarheight}
 h(T)\ll 1+\left(\prod_{\substack{\mathcal P\notin S\\v_{\mathcal P}(r_K(T))>0\text{ odd}}}N\mathcal P\right)^c.
\end{equation}
Here and below primes at which the displayed coefficients are not units are removed once and for all. In particular, $r_K(0)r_K(1)\ne0$.

For any nonempty collection of nonidentity root points let $T_{\rm ag}$ be half the rational trace of the sum of their abscissas in $F$. All the abscissas are positive at both real places, by the duplication identity in Lemma~\ref{lem:aggregate-trace}; hence $T_{\rm ag}>0$. For at most eight indices, with identity branches omitted from the sum, that lemma (and the conjugate-class independence in Lemma~\ref{lem:sqrt3-carrier} for the second field) gives
\begin{equation}\label{eq:new-aggregate}
 B^2\ll 1+h(T_{\rm ag})\ \ \t{with}\ \ B=\max(1,\max_i|n_i|).
\end{equation}

\subsection{The two sets of source equations} Take a nonrecursive recursively enumerable set $\mathcal A\se \N$, and fix the integer precision $m$, larger than all the degree-height bounds specified below. The first reduction uses the split integer source $F_{1,a}$ and $F_{2,a}$ of \cite[Section 2, Lemmas 2.2--2.3]{Sun38}, with the eight root indices
\[
 (f,g,h,k,l,x,y_0,m_1)\ \ (fkx\ne0),
\]
and a field parameter $w_*\ne0,1$, where $a$ ranges over $\N=\{0,1,2,\ldots\}$  The first residual is an integer polynomial in the seven entries $(f,g,h,k,x,y_0,m_1)$. The second uses $(f,g,h,k,l)$, and it is a monic polynomial in $w_*$. Namely,
\begin{equation}\label{eq:new-F2}
 F_{2,a}=(w_*-O_0(D_{\rm size}+1))^2-4O_0^2D_{\rm size},
\end{equation}
where $O_0$ and $D_{\rm size}$ are given in \cite[Section 2]{Sun38}, with its $\varepsilon_+=1$. In particular,
\[
 \begin{gathered} D_{\rm size}=(T_0^2-4)M_0^2A_0^2+4\ \t{with}\ T_0=2M_0+2,\\ 2\mid A_0
  \ \ \t{and}\ \ D_{\rm size}/4\equiv1\pmod4,\end{gathered}
\]
where $M_0$ and $A_0$ are as in \cite[Section 2]{Sun38}. The cited source lemma asserts equivalence to $a\in\mathcal A$ when both equations hold with integer indices, $w_*>1$ integral, and the indicated nonzero guards. In the affirmative direction the five coding indices and $w_*$ can be chosen first, and $x>0$ subsequently chosen arbitrarily large and divisible by any prescribed integer, with the other two indices satisfying the first equation. On the guarded domain with $F_{1,a}=0$, its growth lemma gives $|x|\ll_a1+\log B$. Both degrees in the indices admit a uniform upper bound.

For the second reduction, take the input $X\in K$ and the finite sign choice $\delta\in\{1,-1\}$, with $\beta=3+\delta(X^3+X)$. For integral $X$ choose $\beta\ge3$. Use the three indices $n,t,z$, with $z\in M\Z\setminus\{0\}$ for a large fixed positive integer $M$. The two residuals are
\begin{equation}\label{eq:new-pell}
 \Delta=(t-nz+1)^2-\mathfrak D(n)z^2-1\ \t{and}\  n-\beta\ \t{with}\  \mathfrak D(n)=(1+r_K(n)^2)^2-1.
\end{equation}
The construction of the shifted Pell source and the index divisibility argument are adapted from \cite{LiSun19}. The polynomial $\Delta$ has degree at most $58$. The guard on $z$ is imposed by taking the root point of $z$ to be $[M]$ of a freely represented cyclic point, with nonidentity output. Thus it costs no extra point. For every integral $X$, the choice $n=\beta$ admits arbitrarily large positive $z,t$ with $\Delta=0$ and $z$ divisible by any prescribed integer: use powers of the norm-one unit $1+r_K(n)^2+\sqrt{\mathfrak D(n)}$. Taking its powers divisible by the orders of its reductions modulo prescribed integer moduli also makes $t-nz+1\equiv1$ to arbitrarily high local precision.

In both reductions, denote the first and second ordinary groups by $G=1,2$. Use seven common comparison points $U_j$ in the first reduction, three in the second, and one common anchor $H$. If $P_{G,j}$ is the root point for slot $j$ of group $G$, put $V_{G,j}=U_j-[c_*]P_{G,j}$ and $V_0=H-[c_*]P$, where $c_*$ is a fixed positive multiple of $64$. The slots are the lists given above (for the second reduction $(n,t,z)$ and $(n)$). Set $\rho_E=1/x_E$ and $\mathfrak z=-x_E/y_E$, with their regular zero values at the origin. Use the degree-$m$ truncated formal logarithm $\ell_m$ of $E$ at its origin and its coefficient conjugate $\ell_m^\tau$, and put
\[
 L_{G,j}=\frac{\ell_m(\mathfrak z(U_j))}{\ell_m(\mathfrak z(H))}\ \ \t{and}\ \ 
 L'_{G,j}=\frac{\ell_m^\tau(\mathfrak z(V_{G,j}^\tau))}{\ell_m^\tau(\mathfrak z(V_0^\tau))}.
\] The ordinary list for group $G$ consists of $\rho_E(H),\rho_E(V_0^\tau)$ and $\rho_E(U_j),\rho_E(V_{G,j}^\tau)$ with tag $1$, $L_{G,j}-L'_{G,j}$ with tag $m$, and the source residual evaluated at the $L_{G,j}$ with tag $m$. All denominators are guarded. At a successful unramified place these tests recover the integer indices to precision $m$. For the monic residual \eqref{eq:new-F2}, if $w_*$ has negative valuation its leading term is uniquely of least valuation, since the other coefficients evaluated at the local indices are integral. Hence successful testing implies $w_*$ is integral there and the actual residual has valuation at least $m$. The same comparison proves this bound for $n-\beta$ in the second reduction, without any preliminary integrality hypothesis on $\beta$.

The scalar list in the first reduction comprises the two coefficients in the basis of $L/K$ of each of $\rho_E(H),\rho_E(P_x),\mathfrak z(P_x)/\mathfrak z(H)$, with tags $1,1,m$ respectively. In the second reduction take the two coefficients of each of
\[
 \rho_E(H),\quad\rho_E(P_t),\quad\rho_E(P_z),\quad
 \mathfrak z(P_z)/\mathfrak z(H),\quad \mathfrak z(P_t)/\mathfrak z(P_z)-\beta
\]
with tags $1,1,1,m,1$. We use the bases $1,\sqrt2$ or $1,\sqrt3$, as appropriate. The coefficient functions lie in $K$ even for nonrational input; in the affirmative integer constructions they are rational. Outside the fixed ramified and bad primes, successful scalar tests give $p^m\mid x$ in the first reduction and $p^m\mid z$ in the second. Indeed the first tests enter the formal group, and the isometric formal logarithm computes the valuation of the parameter quotient as the difference of the integer-index valuations, as in \cite[Lemma 3.2]{Sun38}. If $\beta$ is integral at the chosen prime ideal, the last test in the second reduction also gives $t/z\equiv\beta$ there; the raw quotient and the logarithm quotient differ by a positive-valuation quantity when the quotient is integral. All the scalar coefficient inequalities hold in $K$, so the conclusion is valid for the same chosen prime ideal, whether split or inert in $L/K$.

\subsection{Fixed ordinary positions and a square scale}\label{subsec:new-positions}
For each scale $j=0,1$ and each group $G=1,2$ choose four translates if $q=367,331$, or five if $q=347$. They are of the form $\kappa_l=lD_\sharp/Q+d_{jG}$, where $D_\sharp$ is coprime to $Q$, the spacing exceeds $4q+\operatorname{diam}\Theta_K+2$, and the integer offsets $d_{jG}$ are chosen large enough to make every number $\kappa_l-\theta$ positive. Choose offsets separated enough that all these numbers within a scale are distinct. The arguments are $t_{jGl}=A_j\eta+\kappa_l$, with $A_0=1,A_1=A$ and a single nonzero slope $\eta\in K$. Their seven contacts are the factors $A_j\eta+\kappa_l-\theta$. At each distinguished prime at most one translate fails the tube, since multiplication by the rational prime gives distinct leading residues in the only possible cancellation case. At each real place at most one translate fails $|t|>2q$ (or $t\notin[-q,0]$ in the second field). For the first and third fields, $1/2<1+q/t<3/2$ gives $0<x(z)<1$ by the displayed pullback formula, so the required sign polynomial is positive. For the second field the surviving argument gives $z>0$ and Lemma~\ref{lem:sqrt3-seven} supplies its sign condition. Four or five translates therefore suffice, and the height estimate applies to a surviving argument in each group at each scale. This holds for arbitrary nonzero $A\in K$.

Write $I_g,J_g$ for the genuine contact positions in the first and second scales respectively, subdivided according to $G$. Their sizes are the same even integer $N=2\cdot7\cdot(4\text{ or }5)$. Enlarge each scale by one parity-adjustment position and one leaf for every genuine or parity-adjustment position. Enlarge the second scale by one more position. Write $I,J$ for the resulting two sets; thus $|I|=2(N+1)$ and $|J|=2(N+1)+1$. Their coefficients will be positive rational numbers $b_i$ and $c_j$, respectively. Choose distinct positive rational normalizer constants $e_0,e_1$, different from the contact coefficients at the corresponding scale. For clarity the actual slope positions are $\lambda_i=b_i$ ($i\in I$) and $\lambda_j=c_j/A$ ($j\in J$). The normalizer positions are $e_0,e_1/A$.

We require the following fixed private-prime property. For each unordered pair of contact positions with at least one auxiliary position, choose a separate sufficiently large odd prime $\ell_e$. At a provisional positive integral square $a_*$ larger than the genuine cross-scale ratios, choose these primes avoiding all genuine-position and normalizer differences at $a_*$. By CRT choose the auxiliary coefficients (positive integers suffice) such that, at $A=a_*$, the assigned pair of slope positions agrees modulo $\ell_e$ but not modulo $\ell_e^2$, and every other pair of positions, including the normalizers, is distinct modulo $\ell_e$. This can be done one prime at a time: if one endpoint is genuine its residue is prescribed, otherwise choose a fresh residue, and give the auxiliary endpoint a different lift modulo $\ell_e^2$. The other residues can all be distinct because the prime exceeds the total number of positions. Multiplication by $a_*$ at the second scale is invertible. The finite CRT prescriptions allow all coefficients positive and distinct within their scale. Now choose a positive integer $C_A$ congruent to $\sqrt{a_*}$ modulo every $\ell_e^2$, large enough that
\begin{equation}\label{eq:new-scale-bound}
 C_A^2>\max_{i,j}c_j/b_i,
\end{equation}
and larger than the cross ratios involving normalizers as well. The private-prime property persists whenever $A\equiv C_A^2\pmod{\ell_e^2}$. This is the reason for adding fixed auxiliary positions rather than trying to prescribe all genuine--genuine characters.

Choose a fixed finite core $S_0$ of rational primes containing all bad primes of the preceding data, all private primes, $2,3,5,7,13,17$, all primes dividing the coefficients or the nonzero differences of the slope and normalizer positions at $A=C_A^2$, the primes dividing $r_K(0)r_K(1)$, and all primes up to a sufficiently large fixed threshold depending on the number of positions and $m$. The threshold is chosen so that the elementary quadratic-character prescription proved below works for every prime outside it. Further fixed enlargements are made before choosing the coefficient primes. Put $\Omega=\prod_{p\in S_0}p$.

Choose three distinct fixed primes $q_I,q'_I,q_J$ outside $S_0$, avoiding $(1+\Omega)(9+\Omega)$, with the following additive Legendre profiles (write $[a/p]=0$ or $1$ according as the nonzero residue of $a$ is a square or a nonsquare):
\[
 \begin{gathered} f(p)=\bigl([{-1}/p],[3/p],[2/p],[5/p]\bigr),\\
 f(q_I)=f(q'_I)=(0,0,0,1)\ \ \t{and}\ \  f(q_J)=(0,0,0,0).\end{gathered}
\]
Their mutual Legendre bits are zero. Existence follows from Dirichlet's theorem applied to the compatible residue and character conditions (the three primes are chosen successively). All the fixed exceptional primes for height estimates include these three primes as well as the finite multiplier primes chosen in the next section. They need not be in $S_0$.

For HTP set $Q_*=1+L_0w_*$, where $L_0\ge2$ is fixed, guard $Q_*\ne0$, and use $\beta=Q_*^2$. For integer definition keep its $\beta$ above. Take
\begin{equation}\label{eq:new-scalar}
 \begin{array}{ll}
 r=\mu\Gamma^2Q_*r_K(\beta),&\rho=\nu\Gamma_1^2r(Q_*-1)\quad\text{(HTP)},\\
 r=\mu\Gamma^2r_K(\beta),&\rho=\nu\Gamma_1^2r\quad\text{(integer definition)}.
 \end{array}
\end{equation}
Here $\mu,\nu$ belong to fixed finite lists of positive squarefree prime products specified below. The fixed integers $\Gamma,\Gamma_1$ provide arbitrarily high positive valuations at the core, and $\Gamma_1$ also at the three coefficient primes. Guard $r\rho\ne0$. Put $L_\rho=1+\rho^2$,
\begin{equation}\label{eq:new-A}
 Y=\frac{L_\rho(T_{\rm ag}^2+\rho)}{T_{\rm ag}^2+L_\rho\rho}\ \ \t{and}\ \  A=C_A^2Y^2.
\end{equation}
Guard all denominators and $A\ne0$. If $\rho>0$ at a real embedding then $A>C_A^2$ there. At any rational prime where $v_p(\rho)>0$ is odd and $T_{\rm ag}$ is rational, we have
\[
 Y-1=\frac{\rho^2T_{\rm ag}^2}{T_{\rm ag}^2+L_\rho\rho}\ \ \t{and}\ \ v_p(Y-1)\ge2v_p(\rho).
\]
This follows from the unequal parities of $v_p(T_{\rm ag}^2)$ and $v_p(L_\rho\rho)$. In particular $A$ has the required private-prime congruences. Finally $T_{\rm ag}^2=L_\rho\rho(Y-1)/(L_\rho-Y)$, with $L_\rho-Y\ne0$ because $\rho\ne0$. The elementary height inequalities give
\begin{equation}\label{eq:new-scale-height}
 h(T_{\rm ag})\ll1+h(A)+h(\rho)\ \ \t{and}\ \  h(\rho)\ll1+h(\beta).
\end{equation}
The latter uses $h(w_*)\le\tfrac12h(\beta)+O(1)$ in the HTP case. These inequalities hold on all guarded branches, without assuming positivity in soundness.

Let $E_0=\eta+e_0$ and $E_1=A\eta+e_1$. Define
\begin{equation}\label{eq:new-raw}
 R_0=q_Iq'_I\frac{\prod_{i\in I}(\eta+b_i)}{E_0^{|I|+2}}\ \ \t{and}\ \ 
 R_1=2q_J\frac{\prod_{j\in J}(A\eta+c_j)}{E_1^{|J|+1}}.
\end{equation}
All factors and normalizers are guarded nonzero. For source $j$ at scale $j$ and group $G$, divide its raw tests $D$ of tag $k$ by $E_j^{2k+1}$. Use the group gate
\begin{equation}\label{eq:new-gates}
 g_{jG}=q_Iq'_Iq_J\frac{\prod_{i\in C_j\setminus C_{jG}}(A_j\eta+d_i)}{E_j^{|C_j\setminus C_{jG}|}},
\end{equation}
where $C_0=I,C_1=J$, $d_i=b_i$ or $c_i$, and $C_{jG}$ consists only of the genuine contacts of the designated group. Thus the gate is a unit at a genuine own-group contact outside fixed primes, vanishes to positive depth at every other selected contact of its scale, and has valuation one at the coefficient primes when all slope and normalizer factors are units. At polar slopes its quotient tends to one. At a genuine odd contact away from the fixed exceptional set, all other factors of the same scale and the normalizer are units; $R_j$ has exactly that odd valuation. This last assertion does not assume that the other scale is a unit.

\section{A five-vertex Selmer certificate and its uniform expansion}\label{sec:new-selmer}
In this section all graph operations are over $\F_2$. For an odd squarefree positive $d\equiv17\pmod{24}$, with prime divisors $q_i$ including $5$, use the profile $f=(s,\lambda,t,\phi)$ just defined and the convention $f(5)=(0,1,1,0)$. Let $W_{ij}=[q_j/q_i]$ for $i\ne j$ and $W_{ii}=0$; quadratic reciprocity gives $W_{ij}+W_{ji}=s_is_j$. The directed Laplacian $\mathcal K$ has off-diagonal entries $W_{ij}$ and diagonal entries the row sums of $W$. Write $D_v$ for the diagonal matrix of a vector $v$.

We use the local equations of \cite[Lemma 12.1]{Sun38} for the curve $E_d:y^2=x(x-d)(x+3d)$. Its positive Kummer pairs in the coordinates $([x+3d],[x-d])$ have the form
\[
 (2^\epsilon3^\delta\prod q_i^{u_i},\;2^\epsilon\prod q_i^{v_i})
\]
with
\begin{equation}\label{eq:new-local-equations}
 \begin{gathered}
 \mathcal Ku+D_\lambda v=\delta\lambda+\epsilon t,\quad
 (\mathcal K+D_s)v=\epsilon t,\\
 \epsilon+(s+\lambda)^Tv=0,\quad \epsilon=\delta+s^T(u+v),\quad
 (s+t)^Tu+t^Tv=0.
 \end{gathered}
\end{equation}
For completeness the rank criterion in \cite[Lemma 12.2]{Sun38} follows as follows. Reciprocity gives $\mathcal K\mathbf1=0$, $\mathcal K^T\mathbf1=s$ and $\mathcal K+D_s=\mathcal K^T(I+\mathbf1s^T)$, with $(I+\mathbf1s^T)^2=I$. If $\operatorname{rank}\mathcal K=n-1$ and $\lambda\notin\im\mathcal K$, its left null vector $l$ satisfies $\lambda^Tl=1$. In the homogeneous system with $\varepsilon=0$, the second vector equation forces $v$ onto the line spanned by $(I+\mathbf1s^T)l$; the first scalar equation kills that line because $\mathbf1^T\lambda=1$ and $\mathbf1^Ts=0$. Then $\delta=0$ and $u\in\langle\mathbf1\rangle$. Existence of one solution with $\varepsilon=1$ gives positive subgroup dimension two, and the sign-changing torsion class gives full Selmer dimension three.

\begin{lemma}[Five-vertex certificate]\label{lem:new-five}
Use vertices $(5,S_1,S_2,D_1,D_2)$ with profiles and target columns
\[
\begin{array}{c|rrrr|rr}
 &s&\lambda&t&\phi&u&v\\\hline
 5&0&1&1&0&1&0\\ S_1&0&0&0&0&1&1\\S_2&1&0&1&1&1&1\\
 D_1&1&1&0&1&0&0\\D_2&0&1&0&1&0&0
\end{array}
\]
and Laplacian
\[
 \mathcal K_*=\begin{pmatrix}1&0&1&1&1\\0&1&1&0&0\\1&1&1&1&0\\1&0&0&1&0\\1&0&0&0&1\end{pmatrix}.
\]
Then \eqref{eq:new-local-equations} holds with $\varepsilon=\delta=1$, $\operatorname{rank}\mathcal K_*=4$ and $\operatorname{rank}[\mathcal K_*\mid\lambda]=5$.
\end{lemma}
\begin{proof}
The off-diagonal entries satisfy reciprocity and the row-sum rule. Substitution gives the two vector right sides $(0,0,1,1,1)^T$ and $(1,0,1,0,0)^T$ and the three scalar equalities. The kernel is the all-ones line. The left kernel is generated by $(0,1,1,1,0)^T$, whose pairing with $\lambda=(1,0,0,1,1)^T$ is one. This proves both ranks.
\end{proof}

\begin{lemma}\label{lem:new-contactgraph}
Any prescribed undirected graph on the genuine positions of $I,J$ admits the following completion. The cross block between $I,J$ has even row and column sums, and the matrix $M_C=I_{|I|+|J|}+L_C$, where $L_C$ is its internal undirected Laplacian, is invertible. Attach every $I$ vertex to $q_I$ and every $J$ vertex to $q_J$, with all other incident edges to the fixed coefficient vertices zero. Eliminating the full contact block leaves the graph on the retained vertices unchanged, provided the column sums to each scalar group have become zero after contraction of that group.
\end{lemma}
\begin{proof}
For the genuine cross block add one row and column, given by the original column and row parities, and put the total parity in the corner. All row and column sums are zero. Within each group retain the given genuine edges, choose the extra intra-group edges zero, then attach one internal leaf to each vertex of the enlarged block. Finally add the internal isolated $J$ vertex. All these additional edges are incident to at least one auxiliary position. The matrix after adding the leaves has blocks $\left(\begin{smallmatrix}B&I\\I&0\end{smallmatrix}\right)$ with inverse $\left(\begin{smallmatrix}0&I\\I&B\end{smallmatrix}\right)$; the extra isolated vertex contributes a diagonal $1$. Cross row and column parity gives $M_C\mathbf1_I=\mathbf1_I$ and $M_C\mathbf1_J=\mathbf1_J$. Hence the diagonal terms of the Schur complement at $q_I,q_J$ cancel their incident degree parities, and the off-diagonal term is zero. The augmented $\lambda$ column is unchanged since every contact profile is zero. The rank rises by precisely the number of contact vertices, for either ordinary or augmented matrix.
\end{proof}

\subsection{Scalar normalization and expansion}
Use the extended profile $\tilde f(p)=(f(p),[p/q_I],[p/q'_I],[p/q_J])$ for primes outside the core and the three coefficient primes. Select a fixed prime for each of the $2^7$ profiles and one further profile-zero reserve prime, all distinct and outside those sets. They exist by CRT and Dirichlet, since the four characters modulo $120$ and the three separate prime characters are independent. Include them in the fixed height-exception set. The list for $\mu$ comprises products with exponents zero or one of primes in the core, the three coefficient primes and these profile and reserve primes. In an affirmative instance the unnormalized scalar in \eqref{eq:new-scalar} is a positive integer. By toggling parities arrange
\begin{equation}\label{eq:new-scalar-class}
 [r]=[2\prod_{p\in\mathcal B}p],
\end{equation}
where $\mathcal B$ is disjoint from the core and coefficient primes. Force the reserve prime into $\mathcal B$, then toggle the other representative primes so that, apart from the reserve, the count of every extended profile is even except $(1,0,1,1,0,0,0)$, which is odd. Pair equal profiles after removing one of that last profile. Put the reserve and the pairs into $\mathcal B_1$ and the removed prime into $\mathcal B_2$. Both groups are nonempty, their four-bit sums are the profiles of $S_1,S_2$, and their three extra bit sums vanish. The argument works with a scalar support of arbitrary finite size.

Choose $\Gamma,\Gamma_1$ sufficiently divisible at the specified fixed primes, and let $\nu$ range over products of the same fixed primes. In affirmative instances make $v_p(\rho)>0$ odd at every core prime, coefficient prime and $p\in\mathcal B$. At fixed primes this is accomplished by the toggle $\nu$, with sufficient initial square factors. At a nonfixed $p\in\mathcal B$ the valuation of $r$ is already positive and odd. In the HTP case $Q_*-1$ is a unit at such a prime, since otherwise $Q_*\equiv1$, $\beta\equiv1$ and $r_K(1)$ is a unit. Thus no input-dependent multiplier is needed. This proves the scale congruence at every prime where it is needed, in particular $A\equiv C_A^2$ modulo every scalar-support and coefficient prime. All constants and multiplier lists are fixed before the input.

\begin{lemma}\label{lem:new-scalargraph}
For the groups $\mathcal B_1,\mathcal B_2$ just constructed there exist finitely many distinct new denominator primes, with product $D$, and prescriptions of the contact characters at the scalar, coefficient and denominator primes such that the entire graph contracts to the five-vertex certificate. At each contact its profile is zero and its target coordinates are $(1,0)$ in $I$, $(0,0)$ in $J$. The desired positive Kummer pair is
\begin{equation}\label{eq:new-kummer}
 ([15rR_0],[r])
\end{equation}
whenever $[R_0]=[q_Iq'_I\prod_{i\in I}p_i]$ and $[R_1]=[2q_JD\prod_{j\in J}p_j]$. All denominator primes can avoid an arbitrary finite set fixed before their selection.
\end{lemma}
\begin{proof}
The target coordinates are $(1,1)$ at every scalar-support prime, $(1,0)$ at $5,q_I,q'_I$, and $(0,0)$ at $q_J$ and all denominator primes. First add $q_I,q'_I$ as leaves at $5$ and $q_J$ as a leaf at $D_1$. These leaves have zero $s,\lambda,t$ and their target agrees with their neighbour. The bits to the scalar groups vanish after contraction by the extended profiles. Their other prescribed edges are zero, except the stated leaf edges.

For each scalar group choose a representative. Introduce one profile-zero denominator vertex joined to the representative and to each other vertex of that group, separately. The corresponding degree-two elimination identifies the two coordinates, adds their incident bits and adds their profiles. It removes two rows and columns and decreases both ranks by two: its row imposes equality of the two coordinates, and one of the two rows solves for the auxiliary coordinate. This applies also to the augmented column, which is zero at an auxiliary vertex. Thus each group contracts to its designated $S_i$. The total scalar--scalar edge between groups may differ from the desired bit one. If it does, introduce one profile-zero denominator vertex joined to both representatives and $D_1$; its degree-three elimination toggles the three edges between these vertices. Indeed the Schur complement is $\mathcal K+\operatorname{diag}(b)+bb^T$ with $b$ their indicator vector. Prescribe the physical $D_1$ edges to the groups so that this toggle, when used, gives exactly the $S_1,D_1$ and $S_2,D_1$ bits of Lemma~\ref{lem:new-five}. Prescribe zero group sums to $D_2$, and the remaining edges between denominator vertices and fixed vertices as required by the five-vertex matrix and the three leaves. The two main denominator primes have the profiles of $D_1,D_2$; all other denominator primes have profile zero. Reciprocity fixes the reverse bits. All these primes, including the degree-two and degree-three auxiliaries, may be chosen in sequence by Dirichlet with the indicated Legendre symbols to all previously chosen primes, and with the stated finite exclusions. The group contraction leaves exactly the five-vertex matrix with its three attached coefficient leaves.

The local equations at a physical scalar vertex $p$ require its incident bits to $I$ and $J$ to have sums, respectively,
\[
 s_p+\phi_p+a_p\ \ \t{and}\ \  t_p+d_p,
\]
where $a_p$ sums its edges to $q_I,q'_I$ and $d_p$ its edges to all denominator vertices and $q_J$ (not to $5$). This follows by inserting $u_p=v_p=1$ into the two vector equations \eqref{eq:new-local-equations}. The sums over each scalar group are zero: the $I$ assertion follows from the profile and extra-bit sums, and the $J$ assertion follows from the group sums of the physical denominator edges after accounting for the degree-three correction (the correction's extra edge cancels the change of the physical $D_1$ edge). In group one $t=0$ and the net denominator bit is zero; in group two $t=1$ and that bit is one. Degree-two auxiliaries contribute twice within a group. Consequently choose the rows with the required sums and with zero column sum within each scalar group, for instance placing all required bits of a group in one designated column. Give all denominator vertices zero edges to contacts. At $q_I$ prescribe the $I$ edges one and the $J$ edges zero, at $q'_I$ all zero, and at $q_J$ the $J$ edges one and the $I$ edges zero.

At each contact the scalar-group column sums are zero, so its two equations reduce to the cross-block parity equations of Lemma~\ref{lem:new-contactgraph}. The two equations at the coefficient vertices follow from their matching target coordinates on the added neighbours. For the denominator vertices they follow from the prescribed group sums, with the auxiliary vertices incident to equal pairs of scalar targets (and possibly a zero target). The degree-three correction leaves the equations at the retained vertices valid because its toggled edges join two equal scalar targets to a zero target. The three global scalar equations depend only on the summed profiles of each group, so are exactly those of Lemma~\ref{lem:new-five}. This verifies \eqref{eq:new-local-equations} before contraction. The contact expansion has zero aggregated scalar edges, so it may be eliminated after the group contractions, and then the three coefficient leaves may be eliminated. Their profiles have $\lambda=0$. Thus the two ranks increase from $4,5$ by precisely the number of added vertices. The total $s,\lambda,t$ sums are $0,1,0$, respectively; the corresponding odd squarefree $d$ is $17\pmod{24}$. Finally the two coordinate products, including their common $2$ and first-coordinate $3$, are precisely \eqref{eq:new-kummer}. This proves the assertion.
\end{proof}

\section{Realizing the entire graph by prime factors}\label{sec:new-realization}
The following realization uses the position-dependent Gaussian ray prime-pattern proposition of \cite[Proposition 10.2]{Sun38}, including its arbitrary fixed length and order-obstruction tail. In particular this is an arithmetic input from that paper, not a consequence of Dirichlet's theorem for the simultaneous contact primes.

Fix affirmative scalar data, root indices and hence $A$ as above. Choose all denominator vertices by Lemma~\ref{lem:new-scalargraph}, avoiding the supports of $r,A$, all slope and normalizer coefficients and differences, $(1+\Omega)(9+\Omega)$, and the fixed exceptional primes. Choose an integer $\Lambda>0$ clearing all rational position coefficients so that $D\Lambda^2\lambda_i$ and the two similarly scaled normalizer positions are integers. At primes in $\mathcal B$ and the coefficient primes these denominators are units, by the scale congruence and the exclusions, so $\Lambda$ can be chosen a unit there. Set
\begin{equation}\label{eq:new-primes}
 \eta=\frac{u}{D\Lambda^2w^2}\ \ \t{and}\ \  p_i=u+D\Lambda^2\lambda_i w^2.
\end{equation}
The denominator contains $D$, not $D^2$. Its squarefree part is essential. The positions are distinct by \eqref{eq:new-scale-bound} and $A>C_A^2$. Require all $p_i$ distinct primes avoiding all fixed data and coefficients; then $w$ is a unit at each selected prime, since otherwise that prime divides every $p_i$. The exact factorization of each numerator and normalizer gives
\begin{equation}\label{eq:new-factorizations}
 [R_0]=[q_Iq'_I\prod_{i\in I}p_i]\ \ \t{and}\ \ 
 [R_1]=[2q_JD\prod_{j\in J}p_j].
\end{equation}
Indeed the contribution of the common denominator $D\Lambda^2w^2$ and square scale $A$ to a product with an even number of numerator factors is a square, while for the odd number it contributes the class of $D$. Both normalizing powers are even. This proves the squareclass identities at every prime, not just the selected primes.

Here are the CRT conditions realizing the graph. At a prime in $\mathcal B$ or among the three coefficient primes, $A\equiv C_A^2$ and all contact and normalizer positions are distinct modulo that prime, with unit denominators. Take $w\equiv1$ there and prescribe the characters of $p_i$ from Lemma~\ref{lem:new-scalargraph}. These prescriptions are possible for every sufficiently large prime: for $k$ distinct positions the number of $u$ realizing any specified nonzero signs is at least
\[
 2^{-k}\bigl(p-((k-2)2^{k-1}+1)\sqrt p\bigr)-k-3.
\]
This is obtained by expanding the product of the $k$ half-indicators $(1+\varepsilon_i\chi(u+a_i))/2$ and applying the squarefree-polynomial Weil character-sum bound, as in \cite[Lemma 10.1]{Sun38}; the three extra exclusions allow both normalizers and $u=0$. We took the fixed threshold large enough for this lower bound to be positive. At every denominator prime require $w$ arbitrarily divisible and $u\equiv1$; thus the characters from the contact vertices to all denominator primes are zero. At every core prime require $w$ arbitrarily divisible. The Gaussian residues are $9\pmod{16}$, $i\pmod3$, $i\pmod5$, $1+3i\pmod7$, and $1+2i\pmod{\ell_K}$, where $\ell_K=13$ for the first and third fields and $17$ for the second. In particular all contact profiles are zero and the carrier orders have the required dyadic valuations $2,6$ and correct norm-$\ell_K$ twist signs. At all other fixed odd primes one may select a nonzero norm with either prescribed character, except when the norm has already been prescribed. At odd primes at least seven, every nonzero Gaussian norm fibre has $p-(-1/p)$ elements, and at most four violate $\Nm(1\pm\pi)\ne0$. Thus all required norms have Gaussian lifts with both orders units. The small primes have the displayed explicit residues. Higher odd norm congruences lift by surjectivity of the norm on principal units.

It remains to realize the contact--contact graph, rather than just its prescribed external rows. Since $p_i\equiv1\pmod4$ and $p_i\nmid w$, we have
\[
 (p_j/p_i)=\bigl(D\Lambda^2(\lambda_j-\lambda_i)/p_i\bigr).
\]
At an odd coefficient prime the contribution is the parity of its valuation in the rational difference times the character of $p_i$ modulo that prime, by reciprocity. Signs and powers of two contribute zero. At private prime $\ell_e$ the actual difference of its assigned pair has valuation one and every other difference valuation zero, by the fixed construction and $A\equiv C_A^2\pmod{\ell_e^2}$. Since $w$ is divisible there, all $p_i$ have the same nonzero residue; its quadratic character independently toggles precisely that edge. At every other prime dividing the coefficients or differences take $w$ divisible and common norm one, unless a previous prescription fixes the contribution. No such prime is a scalar or coefficient prime at which variable characters are required, by the distinctness at $A=C_A^2$; denominator primes were also excluded from the differences. Thus the genuine--genuine graph is fixed by the nonprivate prescriptions. Apply Lemma~\ref{lem:new-contactgraph} to that actual graph and then choose the private-prime characters to realize all its added edges. Every added edge has a private prime. The order of choices is legitimate: $D,\Lambda,A$ and every difference prime are fixed before the character prescriptions, and private primes never alter a genuine--genuine edge.

All additional primes required for the later finite local index moduli and order avoidance can be included in the Gaussian ray modulus with compatible norm conditions, choosing $w$ divisible and norm one at a new prime. All conditions are on a fixed finite modulus; CRT supplies a compatible pair $(u,w)$ modulo it, with each prime position having the specified norm. The pure Gaussian ray case of the proposition has abelian compatibility exactly this norm congruence. Apply the proposition with the position-dependent rays and the fixed modulus, obtaining arbitrarily large prime patterns with $u/w^2\to\infty$, the odd parts of all distinct-position reduction orders pairwise coprime, and avoidance at all prescribed odd primes. The reduction-order identification uses the common CM norm calculation of the original carrier: $N_{F/\Q}(a_E)=13$ or $N_{F/\Q}(A_3)=9\cdot17$, and the two conjugate orders are $d_\pm(p)=\Nm(1\pm\pi)$. In particular each selected $p$ splits in $K$ and $F$. At each selected prime its own numerator factor has valuation one and all other numerator and normalizer factors are units; this follows from the excluded coefficient differences. At denominator primes sufficiently deep polar $\eta$ gives $v_p(R_0)=-2v_p(\eta)$ and $v_p(R_1)=-v_p(\eta)$, with $v_p(\eta)$ negative and odd; $A$ and the coefficients are units there. At the coefficient primes exactly their designated source has odd valuation one. At scalar-support primes both ordinary sources are units. At core primes the two ordinary sources have arbitrarily large positive valuations as $w$ becomes sufficiently divisible. Thus every exterior odd-source pattern in the affirmative construction is a singleton of positive odd depth. These assertions follow also for all other primes from \eqref{eq:new-factorizations} and \eqref{eq:new-scalar-class}.

For $d=5rR_0R_1>0$,  consider the common covering
\begin{equation}\label{eq:new-cover}
 u_1^2=-5R_0 Z,\ \ u_0^2=-r(Z-1),\ \ u_s^2=R_1(Z+3),\ \ u_1u_0u_s\ne0.
\end{equation}
It maps to $E_d$ by $x=dZ$ and $y=du_1u_0u_s$. Its Kummer pair is $([5rR_0],[-5R_0R_1])$; multiplying by the torsion class $(3,-d)$ gives \eqref{eq:new-kummer}. The squarefree representative of $d$ is the odd integer $17\pmod{24}$ determined above. The full graph is the one of Lemmas~\ref{lem:new-five}--\ref{lem:new-scalargraph}, so its Selmer dimension is three. By the pointwise $2$-converse/finiteness implication \cite[Theorem 1.1]{OAIPointwise}, as used in \cite[Corollary 12.7]{Sun38}, $\Sha(E_d/\Q)$ is finite: the $2^\infty$-Selmer corank is at most one by the dimension bound and rational two-torsion. The alternating nondegenerate pairing on finite $\Sha$ makes its $2$-torsion dimension even. Hence rank is one and $\Sha[2]=0$, so the covering class has a rational representative of infinite order (choose one from its coset modulo $2E_d(\Q)$). This use of the pointwise input is essential.

At each exterior singleton prime, the covering class has the valuation vector of the corresponding torsion point with $Z=0,1,-3$, respectively. The local quotient by twice the elliptic group at the odd prime has order four, and those torsion classes distinguish its nonzero elements. Thus the difference from the torsion point is in $2E_d(\Q_p)$, and odd multiples of sufficiently high $p$-divisibility approach it. Here the assertion is the local movement argument of \cite[Lemma 15.1]{Sun38}, with the source labels permuted; the same torsion limits and same local component group apply. At the real place the second Kummer sign is negative, so the point is on the bounded component. The appropriate odd progression is dense on this component: its common difference is a nonzero even multiple of a nontorsion point on the real identity circle. Hence we may simultaneously approach all these finitely many torsion limits and any specified nonempty real open arc with $-3<Z<0$. The three square roots are recovered from the squareclasses, with $u_1,u_s>0$ on that arc and the sign of $u_0$ dictated by $y$. A finite set of nonidentically zero rational-function guards can be avoided. This proves the required rational covering interface.

\section{A common three-source square-denominator certificate}\label{sec:new-merger}
Here and throughout the soundness argument the covering equations \eqref{eq:new-cover} are over $K$, with no primality or positivity conditions imposed on $\eta,r,R_0,R_1$. Write $a_1=-5R_0$, $a_0=-r$, $a_s=R_1$ and $c_1=0,c_0=1,c_s=-3$. Define
\[
 F_Z=Z(Z-1)(Z+3),\quad e_i=F_Z/(Z-c_i),\quad Z_i=u_i e_i,\quad D_Z=Z^2+\Omega.
\]
Use the corrected valuation-normalizing expressions of \cite[Section 14]{Sun38}:
\begin{equation}\label{eq:new-LcQc}
 L_c=\frac{Z_1^5+\Omega Z_0^5+\Omega^2 Z_s^5}{Z_1^4+\Omega Z_0^4+\Omega^2 Z_s^4}\ \ \t{and}\ \ 
 Q_c=\frac{D_Z^4+\Omega F_Z^4}{D_Z(D_Z^3+\Omega F_Z^3)}.
\end{equation}
Put
\begin{equation}\label{eq:new-J}
 J_\Phi=a_1a_0a_s\left(\frac{L_c}{u_1u_0u_s Q_c}\right)^2\ \ \t{and}\ \ \gamma_i=e_i/D_Z.
\end{equation}
The guard $u_1u_0u_s\ne0$ in particular excludes the three branch values. Fix an integer $C\ge6$. For real preparation with a possibly small $R_1$, replace the constant coefficient in the original merger by the rational function
\begin{equation}\label{eq:new-adaptive}
  U=\frac{1-C_*^2J_\Phi}{1+C_*^2J_\Phi}\ \ \t{with}\ \ C_*=C\prod_{i\in\{1,0,s\}}(1+a_i^{-2}).
\end{equation}
The scalar list belongs to source $0$. The two ordinary groups at scale zero belong to source $1$, and those at scale one to source $s$. Write $D$ for a test argument, after the normalization following \eqref{eq:new-raw} in the ordinary case, and $g$ for its group gate in \eqref{eq:new-gates} or for $1$ in the scalar case. Its precision tag is $k$. Set
\begin{equation}\label{eq:new-W}
 \begin{gathered}
 C_{J_\Phi}(D)=\frac{D^2(1+J_\Phi)}{1+J_\Phi D^2},\qquad T_{i,D}=U\gamma_i^{2k}g^{2k+1}C_{J_\Phi}(D),\\
 W=\prod_{i,D}(1+J_\Phi^{1-4k}T_{i,D}^2),\qquad
 \Theta=\Tr_{L/K}(W^{-1}),\qquad B_\Phi=J_\Phi/\Theta^2.
 \end{gathered}
\end{equation}
Guard all denominators and $J_\Phi UW\Theta\ne0$. The terminal relation is
\begin{equation}\label{eq:new-BB}
v^2=(1-b)(s_b^4-b^2)(t_b^4-b^2)\ \ \t{with}\ b=112B_\Phi.
\end{equation}
It has two free parameters $s_b,t_b$. We use the rational sufficiency statement of the square-denominator lemma in \cite[Lemma 2.2]{Bur}: a positive rational $b<1$ with square denominator, $v_2(b)\ge3$ and positive squarefree representative congruent to $7$ modulo $8$ admits such a rational ordinate. Its necessary local pole exclusion over $K$ at any unramified odd prime outside the coefficient primes follows directly: if $v_{\mathcal P}(b)$ is negative and odd, the two terms of each difference $s_b^4-b^2$, $t_b^4-b^2$ have unequal valuations, and each difference has even valuation. The factor $1-b$ has odd valuation. This contradicts the equation.

\begin{proposition}\label{prop:new-detection}
Outside the fixed height-exception set, every odd source whose group gate and normalizer are units has its original tests at their specified precisions at a common place above the given prime ideal of $K$. This holds for the scalar source without a group gate or normalizer.
\end{proposition}
\begin{proof}
Fix a prime ideal $\mathcal P$ outside that set. The three square equations imply that a nonempty odd-valuation pattern of the $a_i$ is a singleton or all three. A singleton $i$ forces $Z$ to be integral and $v_{\mathcal P}(Z-c_i)>0$ odd; the other two differences are units. An all-three pattern forces a polar $Z$ of odd valuation. The product $a_1a_0a_s$ has odd valuation in either case. By \eqref{eq:new-J} the valuation $d=v_{\mathcal P}(J_\Phi)$ is odd. At each place above $\mathcal P$ every nonzero correction $J_\Phi^{1-4k}T_{i,D}^2$ has odd valuation, and hence its inverse factor in $W^{-1}$ is integral. The trace is integral. If $d<0$, $B_\Phi$ has a negative odd valuation, contradicting \eqref{eq:new-BB}.

For $d>0$, if every place above $\mathcal P$ fails some inequality $v(T_{i,D})\ge2k$, its corresponding inverse factor has valuation at least $(4k-1)d-2(2k-1)\ge d$. The remaining factors are integral. Thus $v_{\mathcal P}(\Theta)\ge d$, again making $B_\Phi$ of negative odd valuation. There is a single successful place for all tests. At that place the valuation of $C_*^2J_\Phi$ is odd, so $U$ is a unit even though $C_*$ need not be integral. The formula
\[
 v(C_{J_\Phi}(D))=\max(2v(D),-d)
\]
holds because $1+J_\Phi$ is a unit and the denominator has no valuation tie (zero arguments are interpreted separately). At a singleton active source $e_i$ is a unit and $D_Z$ integral, so $v(\gamma_i)\le0$. At an all-three pole $v(e_i)=v(D_Z)=2v(Z)$, so $\gamma_i$ is a unit. Thus, if its group gate is a unit, $v(T_{i,D})\ge2k$ implies $v(D)\ge k$. A unit normalizer gives the same for the raw argument. This proves the proposition.
\end{proof}

\begin{lemma}\label{lem:new-singleton}
In an affirmative rational construction, let $p\notin S_0$ be a prime at which only source $i$ has an odd valuation, of positive depth $n_i=v_p(a_i)$. Set $b_i'=v_p(c_i^2+\Omega)\ge0$. All group gates are supposed to have nonnegative valuation. It suffices at one place to prepare each active test of tag $k$ to depth $k(n_i+b_i')$, or, for a suppressed test, to have $v(g)\ge n_i+b_i'$. Then a sufficiently close approach of the rational cover point to $Z=c_i$ makes $v_p(B_\Phi)=n_i$.
\end{lemma}
\begin{proof}
At the torsion limit $Z_i$ dominates the two weighted sums in \eqref{eq:new-LcQc} and $L_c/Z_i\to1$, $Q_c\to1$. The identity $a_1a_0a_s(Z_i/(u_1u_0u_s))^2=a_ie_i$ gives $v_p(J_\Phi)=n_i$, independent of the approach depth. The active $\gamma_i$ has valuation $-b_i'$; all the other $\gamma$'s become arbitrarily deep. The adaptive factor $U$ is a unit, by oddness. The clip formula of the preceding proof gives a lower bound $-n_i$ for all arguments and twice the raw valuation for positive arguments. Hence for an active prepared test $v(T_{i,D})\ge2kn_i$; for a suppressed active test the bound is
\[
 -2kb_i'+(2k+1)(n_i+b_i')-n_i\ge2kn_i.
\]
For an inactive source the arbitrarily deep $\gamma$ gives the same bound, using the nonnegative group gates. Every correction at the chosen place has valuation at least $n_i$. At any other place the inverse factors have residues zero or one by the odd-valuation argument. In the affirmative rational case all tests are in $F$ (scalar coefficients in $\Q$), so the trace to $K$ is the rational trace from $F$. Its residue is $1$ or $2$ if $p$ splits, and $2$ if it is inert. All these are nonzero for $p>2$. Thus $\Theta$ is a unit and the asserted valuation follows.
\end{proof}

\begin{lemma}\label{lem:new-core}
Let $q\in S_0$. Suppose $n_i=v_q(a_i)$ are all sufficiently large positive integers (with a fixed lower bound depending only on $q,C$ and the three branch values). At all places over $q$, normalize the valuations to restrict to $v_q$ on $\Q$. If each source-$i$ normalized argument of tag $k$ has valuation at least $k(n_i+1)$, and all its group gates are nonnegative, then every guarded rational covering point has
\[
 d_q=v_q(J_\Phi)\ge\min_i n_i,\qquad v(T_{i,D})\ge2kd_q,\qquad v_q(B_\Phi)=d_q-2v_q(2)>0.
\]
\end{lemma}
\begin{proof}
Put $h_i=v_q(Z-c_i)$, $H_Z=\sum h_i$, $E_i=H_Z-h_i$ and $h_Z=v_q(Z)$. The valuation of $Z_i$ is $(n_i+h_i)/2+E_i$, an integer for rational roots. The weight valuations $0,1,2$ in the sums of fourth and fifth powers force unique least terms, at the same index of least $v_q(Z_i)$ (ties decided by weight). Thus the valuation before the $Q_c$ correction in $J_\Phi$ is $d_0=\min_i(n_i+E_i)$. Also $b'=v_q(D_Z)=\min(2h_Z,1)$; the weights in $Q_c$ give $v_q(Q_c)=\min(0,H_Z-b')$. At a polar $Z$, $E_i=b'=2h_Z$ and $d_q=\min n_i$. At an integral $Z$, $H_Z\ge b'$ and $d_q=d_0\ge\min n_i$. In either case
\[
 d_q\le \max_i n_i+2\max_{i\ne j}v_q(c_i-c_j).
\]
For the integral case choose an index with maximal $h_i$ and use the bounds on the other two $h$'s furnished by the differences. For the polar case the inequality is immediate. Consequently, $v_q(C_*^2J_\Phi)=2v_q(C)-4\sum_i n_i+d_q<0$ once the fixed lower bounds on $n_i$ hold. Thus $U$ is a unit even at the core. The clips have valuations $2v(D)$ since the arguments and $J_\Phi$ have positive depth. Using $v(\gamma_i)=E_i-b'$ gives
\[
 v(T_{i,D})\ge2k(n_i+1+E_i-b')\ge2kd_q.
\]
The last inequality follows from $b'\le1$ in the integral case and $E_i=b'$ in the polar case. All corrections are at least $d_q$-deep, and so is $W^{-1}-1$. The rational trace gives $v_q(\Theta-2)\ge d_q$ even at ramified primes, by the stated normalization of valuations. Taking $d_q>2v_q(2)$ proves the claim.
\end{proof}

\begin{lemma}\label{lem:new-real}
In the affirmative construction, assume the preparation of Lemmas~\ref{lem:new-singleton} and \ref{lem:new-core}. A rational covering point and rational terminal parameters satisfying \eqref{eq:new-BB} exist.
\end{lemma}
\begin{proof}
Here $a_1,a_0<0$ and $a_s=R_1>0$. On the real cover arc $-3<Z<0$ with $u_1,u_s>0$, the expressions in \eqref{eq:new-LcQc} have positive denominators. As $Z\to-3$, $J_\Phi\to12R_1$, while the numerator of $L_c$ changes sign between the two endpoints (the $Z_s$ and $Z_1$ terms respectively dominate), and $J_\Phi$ vanishes at an interior point. Hence $J_\Phi$ takes all positive values below $12R_1$. Our adaptive choice satisfies $C_*^{-2}<12R_1$: for $R_1\ge1$ use $C_*\ge C$, and for $0<R_1<1$ use $C_*\ge CR_1^{-2}$. At a point with $J_\Phi=C_*^{-2}$ the clips are regular at both real embeddings, $U=0$, $W=1$, $\Theta=2$, and $B_\Phi=1/(4C_*^2)$. Since $C\ge6$, on a neighbouring nonempty open arc with $J_\Phi>0$ and $U\ne0$ one has $0<112B_\Phi<1$. All rational functions required nonzero can be taken nonzero there; in particular the inverse factors are positive at both real places. The rational covering interface of Section~\ref{sec:new-realization} gives a rational point simultaneously in this arc and sufficiently close to the torsion limits for all exterior singleton primes. For that point $B_\Phi$ has positive valuation at those primes and at all core primes. At all other primes its valuation is even, since its squareclass is $[a_1a_0a_s]=[d]$. Thus $b=112B_\Phi$ has square denominator and $v_2(b)\ge3$. Its positive squarefree representative is that of $7d$, and is congruent to seven modulo eight because the odd representative of $d$ is $17$ modulo $24$. Removing a common square factor of $7d$ does not change this congruence. Now apply the rational sufficiency of \cite[Lemma 2.2]{Bur} to \eqref{eq:new-BB}.
\end{proof}

\section{Simultaneous completeness for the two sources}\label{sec:new-complete}
We give the order of choices, including the local index moduli needed in the prime pattern. All fixed constants, including the contact configurations, the finite multiplier lists, and $C,\Gamma,\Gamma_1$, are independent of the input. The fixed factors in $\Gamma$ ensure that the scalar $a_0=-r$ has the lower depths needed in Lemma~\ref{lem:new-core} at every core prime on affirmative branches.

First fix an affirmative input. In the HTP reduction choose the five coding indices and $w_*$ as in the split-source lemma before the unbounded indices. In the integer definition choose $n=\beta\ge3$. Make the finite scalar choices of Section~\ref{sec:new-selmer}, fixing $r,\rho$ and their scalar support. At every core prime $q$ prescribe an anchor index in a local congruence class of sufficiently high fixed $q$-adic divisibility (a class with exact finite valuation), so that $H$ belongs to a sufficiently deep formal subgroup at both embeddings of $F$ and its formal parameter has a fixed finite valuation. At every scalar-support prime prescribe such a class as well. For these finite requirements, one can first take a multiple of all the relevant local reduction orders or of an index entering a sufficiently small formal subgroup at bad primes. Refine that multiple by imposing exact additional $q$-adic depths, and require the anchor to be divisible by $64$. The formal parameter and each of its two coefficient functions are continuous on these local classes. In particular, the valuation of the nonzero parameter can be kept fixed. At all places over the specified primes we can make the tag-one scalar arguments involving $H$ as deep as prescribed. Use the following scalar precision requirements: $k(n_0+1)$ at the core and $k(n_0+b'_0)$ at the scalar-support primes, with $n_0=v_p(r)$ and $b'_0=v_p(1+\Omega)$.

In the HTP construction choose $x$ a sufficiently large positive multiple of an integer depending on these local classes and the prescribed depths, and choose the other two indices from the split-source lemma. The small-formal-subgroup logarithm, an analytic isometry after entering a sufficiently deep subgroup, makes both coefficients of $\rho_E(P_x)$ and $\mathfrak z(P_x)/\mathfrak z(H)$ satisfy these scalar bounds throughout the chosen anchor classes. In the Pell construction choose a sufficiently high power of the norm-one unit from \eqref{eq:new-pell}, divisible in its exponent by the orders of the unit in every required finite quotient. This ensures $z$ is divisible by $M$ and to arbitrary depth at the specified primes, $u=t-nz+1$ is congruent to $1$ to arbitrary depth, and $t=nz+u-1$. The identity
\[
 \frac{t}{z}=n+\frac{\mathfrak D(n)z}{u+1}
\]
then makes $t/z$ arbitrarily close to $\beta=n$ locally (also at two, where $v_2(u+1)=1$ for sufficiently deep $u-1$). For sufficiently deep $z,t$ the raw parameter quotient differs from the logarithm quotient by an arbitrarily deep amount. This gives the prescribed depths of all ten scalar coefficients in the chosen anchor classes. Passing between the two coefficients and the two embedded values costs at most a fixed additive depth at each ramified prime, which is included in the chosen margins. This argument uses the analytic logarithm on a small subgroup at the finitely many bad primes, not the good-prime isometry at those primes. It fixes the root points, and hence $T_{\rm ag}$ and $A$, before any ordinary primes are chosen. The root indices can have arbitrarily large maximum in the affirmative HTP case.

Choose the denominator primes as in Section~\ref{sec:new-realization}. At the core and denominator primes we will make the ordinary slope polar. First take a reference anchor $H$ in the scalar index classes, and reference comparison indices all divisible by $64$, avoiding all undefined raw expressions. This is possible: each point ranges over an infinite cyclic progression, the group-law chart functions and the nonzero truncated logarithm polynomial have only finitely many excluded points, and the anchor is chosen different from $O,[c_*]P$. Use the continuity of the raw tests in local index classes at the finitely many core and denominator primes. Refining the index classes around these reference points gives fixed lower bounds on the valuations of all raw ordinary arguments, independently of subsequent choices of indices. These local congruence classes are congruences on the integer indices: continuity follows by first passing to a subgroup on which the $p$-adic formal parameter is analytic in the integer multiplier. They are compatible with the earlier scalar classes and the divisibility by $64$.

Now prescribe polar depths of $\eta$ by requiring high divisibility of $w$ at the core and denominator primes and unit $u$ there. At a core prime the depth of $A_j\eta$ can be made arbitrarily large and negative. If $e=v_q(A_j\eta)<0$ is sufficiently negative, all linear factors have valuation $e$, the normalizer has valuation $e$, and the raw ordinary source depth is $-2e$ at scale zero, $-e$ at scale one, plus its fixed coefficient valuation. Dividing a raw argument by $E_j^{2k+1}$ gains $(2k+1)|e|$, so it gives the bound $k(n_i+1)$ even with the fixed lower bound on the raw argument. This proves the ordinary preparation of Lemma~\ref{lem:new-core}. At a denominator prime only scale one is odd and the same bound proves the active singleton preparation, since the selected denominator primes avoid $(1+\Omega)(9+\Omega)$. All ordinary group gates have nonnegative valuation at these polar slopes. At a coefficient prime all its slope factors are units, and the gates have valuation one. As the active source has depth one and $b_i'=0$ at these primes, the suppressed-test part of Lemma~\ref{lem:new-singleton} suffices without any extra comparison-index conditions.

The finitely many primes dividing the index-class moduli are now prescribed as avoidance primes for the two odd order parts at all ordinary prime positions. This includes the scalar-support primes and all local reference moduli. The compatible Gaussian ray conditions of Section~\ref{sec:new-realization} allow this without changing the squareclass graph. Choose the simultaneous prime pattern with all selected primes sufficiently large to avoid every fixed datum, all index-class moduli, $c_*$, all primes of bad reduction, all coefficient differences at the actual scale, and $(1+\Omega)(9+\Omega)$. In particular, since $u/w^2\to\infty$, all $p_i$ are in $(0.9u,1.1u)$ and all corresponding orders are even and less than $1.5u$ for sufficiently large patterns. No selected $p_i$ divides any of these orders: an even multiple of $p_i$ is at least $2p_i>1.8u$. At a selected prime the two reductions are assigned the primary and conjugate roles, with respective orders $d_{i,1}$ and $d_{i,2}$, one of each of $d_+(p_i),d_-(p_i)$. The odd parts of these orders, for all positions and both roles, are pairwise coprime and are coprime to the previously prescribed index moduli.

Choose the global anchor index $h_H$ in its prescribed local classes, with
\[
 h_H\equiv0\pmod{d_{i,1}},\qquad h_H\equiv c_*\pmod{d_{i,2}},\qquad h_H\equiv0\pmod{p_i}
\]
at every selected prime. The two-power conditions are compatible because $64\mid c_*$ and $64\mid h_H$ in the reference progression, while the order depths are two and six. The odd conditions are compatible by the preceding coprimality and avoidance. The generalized CRT gives an infinite progression; discard the finitely many undefined anchor values. At a selected genuine contact belonging to group $G$, choose each used comparison index $u_j$ locally sufficiently close to the integer $h_Hn_{G,j}$: impose this congruence modulo both reduction orders and an arbitrarily high power of $p_i$. These are compatible for all selected primes and with the reference moduli, again because $h_H$ is divisible by $64$ and the odd parts and selected primes are mutually coprime as required. At auxiliary contacts no raw group comparison is required. Use the generalized CRT to choose the actual $U_j=[u_j]P$ in the resulting classes. On the primary side $H$ is formal, on the conjugate side $H-[c_*]P$ is formal. The formal estimates of Lemma~\ref{lem:shared-two-anchor} and the chosen congruences ensure all ordinary tests of the designated group at a genuine selected contact meet their tags. In detail the exact logarithm ratios are $u_j/h_H$ and $(u_j-c_*n_{G,j})/(h_H-c_*)$, so the high $p_i$-adic congruence gives both quotients within precision $m$ of $n_{G,j}$, and their truncated analogues within the same precision. The exact source residual of that group vanishes. At the selected prime the normalizer is a unit and the other group gate has positive valuation; the raw depth is one, so this gives precisely the active and suppressed preparations of Lemma~\ref{lem:new-singleton}. At auxiliary contacts both group gates suppress the tests.

This completes every local preparation: for exterior scalar primes by the first index progression, for the core by the polar depth and the scalar progression, for the denominator and coefficient primes by the preceding suppression or normalization, and for the genuine and auxiliary contacts by the final CRT. The squareclass identities show that these are all exterior odd-source primes. Lemma~\ref{lem:new-real} produces the rational covering and terminal parameters. Thus every affirmative instance has a certificate satisfying all the required algebraic equalities and guards. The choices of primes, valuations and integer congruences occur only in the proof of existence and are not additional conditions or variables in the certificate.

\section{Two height arguments with one common detector}\label{sec:new-sound}
All variables of the certificate may now belong to $K$. The root points are required to be in $2E(L)$, so their indices are rational integers. Let $M_{jG}$ be the product of $N\mathcal P$ outside the fixed height-exception set at which any of the genuine shifted contacts of group $G$ at scale $j$ has positive odd depth. Each prime ideal is included only once. Also put
\[
 M_{\rm sc}=\prod_{\substack{\mathcal P\notin S\\v_{\mathcal P}(r_K(\beta))>0\ {\rm odd}}}N\mathcal P.
\]
At a genuine contact the appropriate ordinary source has odd valuation and the group gate and normalizer are units, by the fixed difference exclusions within the scale. This is true even when the other scale has an odd contact: the cover's local odd pattern is either a singleton or all three, and Proposition~\ref{prop:new-detection} works in both cases. The two-anchor comparison proves $v_{\mathcal P}(\Delta_G)\ge m$ for the exact residual of that group, with the monic integrality argument for $F_{2,a}$. At a scalar contact the scalar polynomial has an integral input $\beta$ locally by the construction of $r_K$. In the HTP case $Q_*$ is a unit there, since $\beta=Q_*^2$ and $r_K(0)$ is a unit at such a prime. Hence $a_0=-r$ is odd in both reductions and the scalar comparisons give $p^m\mid x$ or $p^m\mid z$, respectively. Each such prime ideal has norm $p$ or $p^2$. Consequently,
\begin{equation}\label{eq:new-radbound}
 \begin{gathered}
 M_{jG}^m\le N_{K/\Q}\mathfrak a(\Delta_G)\le \exp(2h(\Delta_G))\quad(\Delta_G\ne0),\\
 M_{\rm sc}^m\le |x|^2\quad\text{or}\quad M_{\rm sc}^m\le |z|^2.
 \end{gathered}
\end{equation}
Here $\mathfrak a(\alpha)$ denotes the numerator ideal of $(\alpha)$, and the bound follows from the product formula. In particular for a nonzero integral residual the right hand side can be replaced by $|\Delta_G|^2$. Using a field norm of a nonintegral residual in place of the ideal numerator would not be valid.

For every group at each scale, a surviving translate satisfies the relevant height theorem, so $h(\eta)\ll1+M_{0G}^c$ and $h(A\eta)\ll1+M_{1G}^c$. Since $A=(A\eta)/\eta$, this gives $h(A)\ll1+M_{0G}^c+M_{1G}^c$. With \eqref{eq:new-aggregate}, \eqref{eq:new-scalarheight} and \eqref{eq:new-scale-height} we have the two common inequalities
\begin{equation}\label{eq:new-height-common}
 B^2\ll1+M_{0G}^c+M_{1G}^c+h(\beta)\ \ \t{and}\ \  h(\beta)\ll1+M_{\rm sc}^c.
\end{equation}
All constants in these inequalities are uniform for the given field and fixed data. The precision $m$ is fixed so large that $2c/m<1$ and $2cd_*/m<1$, where $d_*$ bounds the index degrees of all source polynomials, including the degree $58$ of the Pell residual. In the following inequalities constants allowed to depend on the HTP input $a$ are denoted by $\ll_a$.

\subsection{The split integer source and the sign of its recovered scalar}
In the HTP case \eqref{eq:new-radbound} gives $h(\beta)\ll1+|x|^{2c/m}$. If $F_{1,a}\ne0$ then \eqref{eq:new-height-common} for group one gives
\[
 B^2\ll_a1+B^{2cd_*/m}+B^{2c/m},
\]
so $B$ is bounded for this input. If $F_{1,a}=0$, the growth lemma for the split source gives $|x|\ll_a1+\log B$ and hence
\[
 h(w_*)+h(\beta)\ll_a1+(\log B)^{2c/m}.
\]
If $F_{2,a}\ne0$, its monic quadratic form and bounded degree in the indices yield $h(F_{2,a})\le C_a+d_*\log B+2h(w_*)$. Substituting \eqref{eq:new-radbound} in the group-two inequality of \eqref{eq:new-height-common} gives
\[
 B^2\ll_a 1+B^{2cd_*/m}\exp((4c/m)h(w_*))+h(\beta).
\]
Since $2c/m<1$, the exponential factor is $B^{o(1)}$. This again bounds $B$ for the input. Neither argument uses an effective height constant.

It remains to recover the sign from a certificate with both residuals zero. If $O_0=0$ then \eqref{eq:new-F2} forces $w_*=0$, excluded. Otherwise $D_{\rm size}$ is a square in $K$. Since $D_{\rm size}/4$ is an integer congruent to one modulo four, it is a square already in $\Q$: if a rational number is a square in $\Q(\sqrt d)$, it is a rational square or $d$ times a rational square. In the latter alternative $D_{\rm size}/4=dv^2$ with $v\in\Q$ implies $v\in\Z$, as $d=2,3,6$ is squarefree. This contradicts the congruence modulo four. Thus $D_{\rm size}$ is an integer square, and \eqref{eq:new-F2} shows $w_*=O_0(\sqrt{D_{\rm size}}\pm1)^2\in\Z$. This formula allows either sign until the following argument.

Now $Q_*\in\Z\setminus\{0,1\}$, so $Q_*(Q_*-1)>0$. The definition \eqref{eq:new-scalar} gives $\rho>0$ at both real places, and \eqref{eq:new-A} gives $A>C_A^2$. At either real embedding, if $R_0<0$ then $\eta<-\min_i b_i$, because the numerator has even degree and its factors are all positive for $\eta>-\min b_i$. The positive scale and \eqref{eq:new-scale-bound} make every factor $A\eta+c_j$ negative. Since $|J|$ is odd and the normalizing denominator power even, $R_1<0$. If $r<0$, the second square equation of \eqref{eq:new-cover} forces $Z>1$, its first then forces $R_0<0$, and its third has a negative right hand side. This is impossible. Hence $r>0$, $Q_*>0$, and $w_*>1$ by the guards and $L_0\ge2$. The split-source lemma now implies $a\in\mathcal A$. No additional free parameter is used for this sign recovery.

\subsection{The three-index shifted Pell source}
In the integer-definition case let $B=\max(1,|n|,|t|,|z|)$. If $\Delta\ne0$, its degree bound and the scalar divisor bound give
\[
 B^2\ll1+B^{116c/m}+B^{2c/m}.
\]
Thus there is a uniform bound on $B$ in this case, independent of $X$ and of $M$. To handle the other case, put $a_n=1+r_K(n)^2$, $D_n=a_n^2-1$. The integer $r_K(n)$ is positive, so $a_n\ge2$ and $D_n$ is not a square. Every integer solution of $u^2-D_nz^2=1$ with $z\ne0$ is, up to signs, given by $(a_n+\sqrt{D_n})^j=u_j+z_j\sqrt{D_n}$ for some $j\ge1$. Indeed $a_n+\sqrt{D_n}$ is the least positive norm-one unit greater than one in $\Z[\sqrt{D_n}]$: any such unit has positive integer second coefficient, and the smallest possible coefficient is one. Dividing an arbitrary positive unit by a suitable power gives the claim. In particular $|z|=z_j\ge2^{j-1}$. If $p\ge5$ divides $D_n$, the binomial expansion gives
\begin{equation}\label{eq:new-pell-LTE}
 v_p(z_j)=v_p(j).
\end{equation}
In fact the leading summand $ja_n^{j-1}$ has that valuation, and the summand with $D_n^l$ for $l\ge1$ has valuation at least $v_p(j)+l-v_p(2l+1)>v_p(j)$ by $(2l+1)\binom j{2l+1}=j\binom{j-1}{2l}$ and $p^l>2l+1$. This also covers the sign choices for $z$.

At a scalar odd contact, the scalar tests give $p^m\mid z$, $t/z\equiv\beta\pmod{\mathcal P}$. Put $u=t-nz+1$, so $u\equiv1\pmod p$ since $t/z$ is integral there. From the exact Pell equation,
\[
 t/z=n+D_nz/(u+1)\equiv n\pmod{\mathcal P}.
\]
Hence $r_K(n)\equiv r_K(\beta)\equiv0\pmod{\mathcal P}$, and $p\mid D_n$. Formula \eqref{eq:new-pell-LTE} implies $p^m\mid j$. Thus the scalar estimate strengthens to $M_{\rm sc}^m\le j^2$ and $h(\beta)\ll1+j^{2c/m}$.

If $n-\beta\ne0$, group two now gives
\[
 B^2\ll1+B^{2c/m}\exp((2c/m)h(\beta))+h(\beta),
\]
using $h(n-\beta)\le\log B+h(\beta)+\log2$. Since $2c/m<1$, this bounds $\log B$ by $O(1+h(\beta))$ (take the logarithm after bounding the sum by a constant multiple of its largest term and absorb the sublinear $\log B$ term). Consequently
\[
 (j-1)\log2\le\log|z|\le\log B\ll1+j^{2c/m}.
\]
This bounds $j$, hence $h(\beta)$ and $B$, uniformly. Choose the fixed $M$ greater than the two uniform bounds on $B$ just proved; the constants do not depend on $M$ since it only restricts the root index $z$. Now $|z|\ge M$ excludes both alternatives. Every certificate therefore has $n=\beta\in\Z$. The input identity $\beta=3+\delta(X^3+X)$ makes $X^3+X$ an integer. As the function $s\mapsto s^3+s$ is strictly increasing on $\R$, the two real conjugates of $X$ are equal, so $X\in\Q$. The monic equation $X^3+X\in\Z$ gives $X\in\Z$ by the rational-root theorem. Together with the completeness of Section~\ref{sec:new-complete}, this proves the exact input equivalence.

\section{Compilation and proofs of the main theorems}\label{sec:new-finish}
The cyclic condition for every root, comparison point and anchor is imposed by representing a half point on $E(L)$ and taking its double. In the finite identity branch the half point is $O$; otherwise choose its two $K$-coordinates for the abscissa in the basis of $L/K$ and regard its ordinate as a finite algebraic choice. The curve equation is monic quadratic in the ordinate over $L$; over $K$ this gives a finite fibre, as can also be seen by splitting the two embeddings after passage to an algebraic closure. The guard on the $z$ index is imposed by $[M]$ of such a double. All multiplication and addition maps and parameter functions use finitely many algebraic charts. Identity branches and the guards on root indices are expressed by point identities, not by tests on unspecified integer indices. For HTP and each integer cutoff $N$, impose also that at least one root point has index of absolute value greater than $N$, by excluding the finitely many tuples of points with all indices in $[-N,N]$. These exclusions are effectively expressible, using the fixed generator $P$ and its multiples. They do not add any free coordinate.

For the cover, choose the coordinate $u_1\in K$ and compute $Z=u_1^2/a_1$. The remaining two roots $u_0,u_s$ are finite square-root choices. The terminal ordinate is finite over $s_b,t_b$. The multiplier choices and the input sign choice are finite. All remaining expressions are rational functions or finite chart branches of this data. Lemma~\ref{lem:finite-algebra-compiler} packs every finite choice, all chart conditions and the nonvanishing guards into one additional $K$-unknown. Thus the following are the two free-coordinate counts, excluding the input in the integer definition:
\[
\begin{array}{l|rr}
 &\text{HTP}&\text{integer definition}\\ \hline
 \text{root half-point abscissas}&16&6\\
 \text{comparison half-point abscissas}&14&6\\
 \text{anchor half-point abscissa}&2&2\\
 \eta&1&1\\ u_1&1&1\\ s_b,t_b&2&2\\ w_*&1&0\\
 \text{finite-choice compiler}&1&1\\ \hline
 \text{total}&38&19
\end{array}
\]
All coefficients before descent lie in the fixed number field $K$. Lemma~\ref{lem:coefficients}, equivalently clearing denominators and taking the product of the two coefficient-conjugates, replaces each compiled polynomial by an integer-coefficient polynomial. In the HTP decision problem the property of having a zero in $K$ is invariant under this product construction. For the input definition a zero of a conjugate factor at $X$ gives a zero of the original factor at the conjugate input, which is an integer if and only if $X$ is. Hence descent preserves the exact input equivalence.

\begin{proof}[Proof of Theorem~\ref{thm:unconditional}]
Assume an algorithm deciding solvability of arbitrary integer polynomials in $38$ variables over $K$. For $a\in\N$ and $N=1,2,\dots$, run it on the compiled HTP certificates with the indicated cutoff. The fixed constants can be hardwired, and the polynomials for each $a,N$ can be constructed effectively by the fixed algebraic charts and the compiler.
Recall that $\mathcal A$ is a nonrecursive r.e. subset of $\N$.
 If $a\in\mathcal A$, Section~\ref{sec:new-complete} gives solutions for every cutoff. If $a\notin\mathcal A$, Section~\ref{sec:new-sound} shows that all certificates have bounded root indices, since two exact residuals would imply membership. Some cutoff therefore returns insoluble. Running these tests alongside the enumeration of the recursively enumerable set $\mathcal A$ decides membership: accept on enumeration and reject on an insoluble cutoff. This contradicts nonrecursiveness. This parallel decision argument is why no effective height bound is necessary.
\end{proof}

\begin{proof}[Proof of Theorem~\ref{thm:integers25}]
For the chosen fixed $M$, Section~\ref{sec:new-sound} proves that every certificate implies $X\in\Z$, and Section~\ref{sec:new-complete} proves existence for every integral input. The common compiler and the count above give a polynomial in the input and nineteen quantified variables, over $\Z$ after descent. This establishes the asserted definition over each of the three fields.
\end{proof}

\end{document}